\documentclass[11pt]{article}

\usepackage[utf8]{inputenc}

\usepackage{amsmath}
\usepackage{CJK}
\usepackage{amsfonts}
\usepackage[dvips]{color}
\usepackage{amssymb}
\usepackage{ leftidx}
\usepackage{mathdots}
\usepackage{amssymb}
\usepackage{graphicx,epstopdf}
\usepackage{booktabs}
\usepackage{mathrsfs}
\usepackage{color}
\usepackage{indentfirst}
\usepackage{amsmath,url,cases,supertabular}
\usepackage{amssymb,mathtools,multirow}
\usepackage{enumerate,makecell}
\usepackage{amsfonts}
\usepackage{amsthm}
\usepackage{graphicx,epstopdf}
\usepackage{color,verbatim}
\usepackage{mathrsfs,subeqnarray,subfigure}
\usepackage{listings,array}
\usepackage{booktabs}
\usepackage{dashrule}
\usepackage{arydshln}
\usepackage[vlined,ruled]{algorithm2e}
\usepackage{amsmath}
\usepackage{stmaryrd}
\usepackage{txfonts}
\usepackage{amsfonts}
\usepackage{mathbbold}
\usepackage{amsfonts}
\usepackage{amssymb}
 \usepackage[table,xcdraw]{xcolor}
\RequirePackage{color,graphicx}

\usepackage{cite}

\makeatletter

\newcommand{\Rmnum}[1]{\expandafter\@slowromancap\romannumeral #1@}
\makeatother
\usepackage{fancyhdr}
\usepackage{amsfonts,latexsym,lscape}

\usepackage{epsfig}
\numberwithin{equation}{section}
\newtheorem{theorem}{Theorem}[section]
\newtheorem{proposition}{Proposition}[section]
\newtheorem{lemma}{Lemma}[section]
\newtheorem{corollary}{Corollary}[section]

\newtheorem{remark}{Remark}[section]
\newtheorem{assumption}{Assumption}[section]
\newtheorem{example}{Example}[section]

\newfont{\Bb}{msbm10 scaled\magstep1}

\newcommand{\Z}{\mbox{\Bb Z}}

\newcommand{\R}{{\mathbf R}}

\def\sqr#1#2{{\vcenter{\hrule height .#2pt
      \hbox{\vrule width .#2pt height#1pt \kern#1pt\vrule width.#2pt}
                       \hrule height.#2pt}}}

\def\L{{\cal L}}\def\S{{\cal S}}

\def \[{\begin{equation}}
\def \]{\end{equation}}
\def \dist{\hbox{dist}}

\def\B{{\cal B}}

\def\L{{\cal L}}
\def\S{{\cal S}}

\def\T{{\cal T}}

\def\bdes{\begin{description}}
\def\edes{\end{description}}
\def\benu{\begin{enumerate}}
\def\eenu{\end{enumerate}}
\def\bitm{\begin{itemize}}
\def\eitm{\end{itemize}}

\def\L{{\cal L}}

\def\R{{\sl I\kern-3.2pt R}}

\def\r{\right}

\def\S{{\cal S}}

\def\sqr#1#2{{\vcenter{\hrule height .#2pt
      \hbox{\vrule width .#2pt height#1pt \kern#1pt\vrule width.#2pt}
                       \hrule height.#2pt}}}

\def\Z{{\mathbf{Z}}}

\def\r1{{\mathbb{R}}}

\def\Z{{\cal Z}}

\def\W{{\cal W}}

\def\X{\mathcal{X}}
\def\Y{\mathcal{Y}}
\def\Z{\mathcal{Z}}

\title{Online Alternating Direction Method of Multipliers for Online Composite Optimization\thanks{Supported by National Key R\&D Program of China under project No. 2022YFA1004000.}}
\author{
Yule Zhang\footnote{School of Science, Dalian Maritime University, Dalian 116026, China. (ylzhang@dlmu.edu.cn) This author was supported by the Natural Science Foundation of China under No. 12201097.},\,\,
Zehao Xiao\footnote{Institute of Operations Research and Control Theory,
School of Mathematical Sciences, Dalian University of Technology,
Dalian 116024, China. (zehao xiao@mail.dlut.edu.cn)},
\,\, Jia Wu\footnote{Corresponding author. Institute of Operations Research and Control Theory,
School of Mathematical Sciences, Dalian University of Technology,
Dalian 116024, China. (wujia@dlut.edu.cn) This author was supported by the Natural Science Foundation of China under No.12071055.}\,\, and \,\,  Liwei Zhang\footnote{Institute of Operations Research and Control Theory, School of Mathematical Sciences, Dalian University
of Technology, Dalian 116024, China. (lwzhang@dlut.edu.cn) This author was supported by the Natural Science Foundation of China under No. 12371298 and partially supported by Dalian High-level Talent Innovation Project No. 2020RD09.}
}

\date{\today}
\begin{document}
\maketitle
 \vspace{2mm}
\begin{center}
\parbox{13.5cm}{\small
\textbf{Abstract.}
In this paper, we  investigate regrets of an online semi-proximal alternating direction method of multiplier (Online-spADMM) for  solving online linearly constrained convex composite optimization problems. Under mild  conditions, we establish ${\rm O}(\sqrt{N})$ objective regret and ${\rm O}(\sqrt{N})$ constraint violation regret at round $N$ when  the dual step-length is taken in $(0,(1 +\sqrt{5})/2)$ and penalty parameter $\sigma$ is taken as $\sqrt{N}$. We explain that the optimal value of parameter $\sigma$ is of order ${\rm O}(\sqrt{N})$.   Like the semi-proximal alternating direction method of multiplier (spADMM), Online-spADMM has the advantage to resolve the potentially non-solvability issue of the subproblems efficiently. We show the usefulness of the obtained results when applied to 
{different types of online optimization problems and verify the theoretical result by numerical experiments}. The inequalities established for Online-spADMM are also used to develop iteration complexity of the average update of spADMM  for solving linearly constrained convex composite optimization problems.
\\[10pt]
 }
\end{center}

\section{Introduction}
\setcounter{equation}{0}
In online  optimization, a decision maker (or a online player) makes decisions iteratively. At each round of decision, the outcomes associated with the decisions are unknown to the decision maker. After committing to a decision, the decision maker suffers a loss. These losses are unknown to the decision maker beforehand.

The Online Convex Optimization(OCO) framework models the feasible  set as a convex set $\Phi \subset {\cal U}$, where ${\cal U}$ is a linear space. The costs are modeled as  convex functions over ${\cal U}$.  A learning framework for OCO problems can described as: at round $t$, the online player chooses $u^t\in \Phi$. After the player has committed to this choice, a convex cost function $\psi_t \in \Psi: {\cal U} \rightarrow \overline \Re$
is revealed. Here $\Psi$ is the family of cost functions available to the adversary. The cost incurred by the online player is $\psi_t(u^t)$, the value of the cost function for the choice $u^t$. Let $T$ denote the total number of  rounds.

Let ${\cal A}$ be an algorithm for OCO, the regret of ${\cal A}$ after $T$ iterations is defined as:
\begin{equation}\label{eq:regret}
{\rm regret}_T({\cal A})=\displaystyle  \sup_{\{\psi_1,\ldots,\psi_T\}\subset \Psi} \left\{\sum_{t=1}^T \psi_t(u^t)-\min_{u\in \Phi}\sum_{t=1}^T \psi_t(u)\right\}
\end{equation}
 There are  a large number of algorithms for online convex optimization problems under different scenarios, among them the famous ones include Follow-the-leader \cite{Kalai2005},  Follow-the-Regularized-Leader  \cite{Shai2007a},\cite{Shai2007b},  Exponentiated Online Gradient \cite{Kiv1997}, Online Mirror Descent, Perceptron \cite{Rosenblatt1988} and Winnow \cite{Little1998}. There are a  lot of publications concerning algorithms for online convex optimization, see Chapter 7 of \cite{MRTal2012}, Chapter 21 of \cite{SSSBen2014}, and  survey papers Shalev-Shwartz \cite{Shai2011}, Hazan \cite{Hazan2015} and references cited in these two papers.

For most works in literature, as pointed out by \cite{Hazan2015}, there  are some restrictions for OCO:
 the losses determined by an adversary should not be allowed to be unbounded and the decision set must be somehow bounded and/or structured.  We know  from  \cite{Shai2011}  and \cite{Hazan2015} that $\psi_t$ is usually not allowed to take infinite values and $\Phi$ is only of simple structures. For examples, $\psi_t$ is required to be Lipschitz continuous or strongly convex, and/or $\Phi$ is the simplex set, the positive orthant, ball-shaped set, or box-shaped  set.

 This paper will eliminate the mentioned  restrictions by permitting  $\psi_t$ to take $+\infty$, this allows us to deal with complicated convex constraint sets. We will explain this point after we introduce the optimization model considered in this paper.

  In this paper, we consider the online composite optimization defined by
\begin{equation}\label{onlineP}
\psi_t(u)=f_t(x)+g(z), \quad \Phi=\{(x,z)\in {\cal U}: Ax+Bz=c\},
\end{equation}
with ${\cal U}={\cal X}\times {\cal Z}$, $u=(x,z)\in {\cal U}$, where ${\cal X}$ and ${\cal Z}$ are two finite-dimensional real Hilbert spaces each equipped with an inner product $\langle\cdot,\cdot \rangle$ and its induced norm $\|\cdot\|$, $f_t:{\cal X}\rightarrow (-\infty,+\infty)$ and $g:{\cal Z}\rightarrow (-\infty,+\infty]$ are proper closed convex functions, $A:{\cal X} \rightarrow {\cal Y}$ and $B:{\cal Z}\rightarrow {\cal Y}$ are two linear operators respectively, with ${\cal Y}$ being another finite-dimensional real Hilbert space equipped with an inner product $\langle\cdot,\cdot \rangle$ and its induced norm $\|\cdot\|$ and $c \in {\cal Y}$.
Namely, at round $N$, the online player is trying to solve
\begin{equation}\label{AonlineP}
\begin{array}{ll}
\min & F_N(x,z)=\displaystyle \sum_{t=1}^N[f_t(x)+g(z)]\\
{\rm s.t.} & Ax+Bz=c,x \in {\cal X},z\in {\cal Z}.\\
\end{array}
\end{equation}
For  online optimization problems  with no constraint in or  simple structured constraints embedded in function $g$ (for example, probability simplex is embedded in the entropy function), there are a large number of publications in machine learning filed for designing algorithms, among them see for example \cite{JohnSinger2009}, \cite{Xiao2010}, and references in \cite{Shai2011}  and \cite{Hazan2015}. Besides the mentioned literatures, there are some recent works related to online optimization (\ref{onlineP}). Mahdavi et. al \cite{MRYang2012}
designed a gradient based algorithm to achieve ${\rm O}(\sqrt{N})$ regret and ${\rm O}(N^{3/4})$  constraint violations for an online optimization problem whose constraint set is defined by a set of inequalities of smooth convex functions.  Recently
 Jenatton et. al \cite{JJArch2015} and  Yu and  Neely \cite{YNeedly2016} developed new algorithms to improve the performance  in comparison to prior works. However the online optimization model considered in these  papers does not cover model (\ref{onlineP}) as
 $\phi_t$ in their problem is required to be smooth and is not permitted to take $+\infty$ values.

 Now we explain that the online composite optimization model (\ref{onlineP}) covers many popular  online optimization problems.
 \begin{example}\label{ex-1}
 Consider a general  online optimization model in which the cost function at round $t$ is $\phi_t:{\cal X} \rightarrow \Re$   and the  constrained set is of the form:
 $$
 C:=\{x\in {\cal X}: Nx+b \in {\cal K}\},
 $$
 where  $N: {\cal X} \rightarrow {\cal V}$ is a linear mapping, $b \in {\cal V}$, ${\cal K} \subset {\cal V}$ is closed convex set and ${\cal V}$ is a Hilbert space.
  Define
  $$
  f_t(x)=\phi_t(x), \quad  g(z)=\delta_{{\cal K}}(z), \quad A=N, \quad B=-{\cal I}
    $$
 where $\delta$  is the indicator function, ${\cal I}$ is the identity in ${\cal V}$ and ${\cal Z}={\cal V}$.  Then the online optimization problem is  expressed as  the form
(\ref{onlineP}).
\end{example}
 \begin{example}\label{ex-2}
 Consider an online optimization model in which the cost function at round $t$ is $\phi_t:{\cal X} \rightarrow \Re$   and the  constrained set is a simple convex set $X$. For avoiding decision jumping, online player introduces a regularizer $R:{\cal X}\rightarrow \Re$.
   Define
  $$
  f_t(x)=\phi_t(x), \quad  g(z)=R(z)+\delta_X(z), \quad A={\cal I}, \quad B=-{\cal I}
    $$
 where  ${\cal I}$ is the identity in ${\cal X}$.  Then the online optimization problem is  expressed as  the form
(\ref{onlineP}).
\end{example}

The off-line problem, in which all losses are known to the decision maker beforehand, corresponds to the case
where $f_k(x)=f(x)$,$\forall k=1,\ldots, N$.
In this case, Problem (\ref{onlineP}) or Problem (\ref{AonlineP}) is reduced to
\begin{equation}\label{sconvex0}
\begin{array}{ll}
\min & f(x)+g(z)\\
{\rm s.t.} & Ax+Bz=c,x \in {\cal X},z\in {\cal Z}.\\
\end{array}
\end{equation}
The convex composite   optimization problem (\ref{sconvex0}) is an important optimization model widely distributed in scientific and engineering fields, see  examples considered in  \cite{boyd11}. Alternating direction methods of multipliers  for solving Problem (\ref{sconvex0}) are an important class of numerical algorithms, which are  extensively studied in recent twenty  years. The classic  ADMM  was   designed   by Glowinski and Marroco \cite{glo75} and  Gabay and Mercier \cite{Gabay:1976ff} and its construction  was much influenced by Rockafellar's works on proximal point algorithms   for solving the more general  maximal monotone inclusion problems \cite{rockafellar76A,rockafellar76B}.

An important progress in the ADMM  field is
the  semi-proximal ADMM (in short, spADMM)   proposed  by Fazel et al. \cite{Fazel}. This method has several advantages. First, it allows  the dual step-length to be at least as large as  the golden ratio of $1.618$. Second,  spADMM not only covers  the classic ADMM but also  resolves  the potentially non-solvability issue of the subproblems in the classic  ADMM.  Third,  perhaps more important one,  it possesses the abilities of handling   multi-block convex optimization problems.
For example, it has been shown most recently that the spADMM  is quite efficient in solving multi-block convex composite semidefinite programming problems \cite{STYang2015,LSToh2014, CSToh2015} with a low to medium accuracy. Importantly, under the calmness of  the inverse Karush-Kuhn-Tucker mapping, spADMM has the linear rate of convergence, this result was established by Han, Sun and Zhang \cite{HSZhang2018}.

Inspired by spADMM for Problem (\ref{sconvex0}), we construct the following  Online-spADMM  for the online problem (\ref{onlineP}). For a Hilbert space $\Z$, for any   self-adjoint   positive semidefinite linear operator $\B:\Z \to \Z$,  denote  $\|z\|_\B:=\sqrt{\langle z , \B z\rangle }$ and   $$\dist_\B(z,D)=\inf_{z'\in D}\|z'-z\|_{\B}$$ for any $z\in \Z$ and any set $D\subseteq \Z$. If $\B$ is the identity mapping in $\Z$, namely $\B=\I$, we use $\dist (z,D)$ to denote the distance of $z$ from $D$.
At round $k\in \textbf{N}$, for problem
\begin{equation}\label{ad}
\begin{array}{ll}
\min & f_k(x)+g(z)\\
{\rm s.t.} & Ax+Bz=c,x \in {\cal X},z\in {\cal Z},\\
\end{array}
\end{equation}
the augmented Lagrangian function
is defined by
\[\label{augL}
{\cal L}^k_{\sigma}(x,z;y):=f_k (x)+g (z)+\langle y,Ax+Bz-c\rangle+\displaystyle \frac{\sigma}{2}\|Ax+Bz-c\|^2,\ \forall\, (x,z,y)\in \X\times \Z\times \Y.
\]
Then   Online-spADMM may be described as follows.\mbox{}\\[0pt]
{\bf Online-spADMM}: An online  semi-proximal alternating direction method of multipliers for solving the online convex optimization problem (\ref{onlineP}).

\begin{description}
\item[Step 0 ] Input $(x^1,z^1,y^1)\in {\cal X} \times {\cal Z}\times \mathcal{Y}.$ Let  $\tau\in (0, +\infty)$ be a positive parameter   (e.g., $\tau \in ( 0, (1+\sqrt{5} )/2)$\,),  $\S_1:\X\to \X$ and   $\T:\Z\to \Z$ be a self-adjoint positive semidefinite, not necessarily positive definite, linear operators. Set $k:=1$.
\item[ Step 1]
Set \begin{equation}\label{xna}
\begin{array}{l}
x^{k+1}\in \hbox{arg}\min \, \L^k_{\sigma }(x,z^k; y^k) +\displaystyle\frac{\sigma}{2}\|x-x^k\|^2_{\S_k}\, ,\\[2mm]
z^{k+1}\in \hbox{arg}\min\, \L^k_{\sigma  }(x^{k+1},z; y^k) +\displaystyle\frac{1}{2}\|z-z^k\|^2_\T\, ,\\[2mm]
y^{k+1}=y^k+\tau\sigma (Ax^{k+1}+Bz^{k+1}-c).
\end{array}
\end{equation}
\item[ Step 2] Receive a cost function $f_{k+1}$ and incur loss $f_{k+1}(x^{k+1})+g(z^{k+1})$ and constraint violation
$\|Ax^{k+1}+Bz^{k+1}-c\|$.
\item[ Step 3] Choose a self-adjoint positive semidefinite linear operator $\S_{k+1}:\X\to \X$ .
\item[ Step 4] Set $k:=k+1$  and go to Step 1.
\end{description}
\quad \, To our knowledge,  the first online alternating direction method perhaps is the one proposed by
Wang and  Banerjee \cite{WangBan2013}. Their algorithm corresponds  to  a modified version of the above Online spADMM where $\T=0$ and the term $\displaystyle\frac{\sigma}{2}\|x-x^k\|^2_{\S_k}$ is replaced by $\sigma B_{\phi}(x,x^k)$, where $B_{\phi}$ is the Bergman of a smooth convex function $\phi$. 
{Based on the Wang and  Banerjee \cite{WangBan2013}, Hosseini et. al \cite{hosseini2014online} extend the online ADMM algorithm to a distributed setting based on dual-averaging. This algorithm is applicable to solve the online convex optimization over a
network of agents and attain a sub-linear regret bound of ${\rm O}(\sqrt{N})$ for the objective function and linear local constraints violation. Liu et. al \cite{liu2018zeroth} design and analyze a new zeroth-order online algorithm. It can not only perform gradient-free operation, but also employing the ADMM to accommodate complex structured regularizers. Compared with the first-order gradient based online algorithm, it requires $m$ times more iterations, where $m$ is the number of optimization variables.}

For Online-spADMM, we define objective and constraint violation  regret by
\begin{equation}\label{eq:ob-regret}
{\rm regret}^{\rm obj}_N=\displaystyle \sum_{t=1}^N  \left[ f_t(x^t)+g(z^t)\right]-\min_{Ax+Bz=c} \left\{\sum_{t=1}^N f_t(x)+g(z)\right\}
\end{equation}
and
\begin{equation}\label{eq:con-regret}
{\rm regret}^{\rm ctr}_N=\displaystyle \sum_{t=1}^N  \left[ \|Ax^t+Bz^t-c\|^2\right],
\end{equation}
respectively.

As far as we are concerned, the main contributions of this paper can be summarized as follows.
\begin{itemize}
\item Cost functions are proper lower semi-continuous convex extended real-valued functions and this makes the optimization model (\ref{onlineP}) include more online optimization problems. For instance, Problem (\ref{ad}) includes linear semi-definite programming, quadratic semi-definite programming and convex composite programming.
\item The proposed Online-spADMM  allows  the dual step-length to be at least as large as  the golden ratio of $1.618$, which is independent on the time horizon and other parameters.
\item When $\sigma=\sqrt{N}$ and $\S_k$ is chosen in a smart way, under mild assumptions (these assumptions are quite weaker than those in \cite{WangBan2013}), it is proved that the regret of objective function of $N$ iterations is  of order ${\rm O}(\sqrt{N})$, and the regret of constraints of $N$ iterations is  of order ${\rm O}(\sqrt{N})$. It is proved that the solution regret  is  of order ${\rm O}(\sqrt{N})$ under strong assumptions.
\item It is proved that, for the average of the first $N$ iterations by spADMM for solving linearly constrained convex composite optimization problems, the iteration complexity of objective function  is  of order ${\rm O}(\sqrt{N})$ and  the iteration complexity of constraint violation is  of order ${\rm O}(\sqrt{N})$ when $\sigma=\sqrt{N}$.
{\item[$\bullet$] We apply the Online-spADMM to solve several examples of online linear constrained optimization problems. The theoretical results are verified by numerical experiments. Compared with other numerical results, Online-spADMM performs very well in all respects, especially the running time}.
\end{itemize}

The remaining parts of this paper are organized as follows. In Section 2, we develop two important inequalities, which play a key role in the analysis for objective regret, constraint violation regret and solution error regret. Section 3 establishes bounds of objective regret, constraint violation regret and solution error regret of Online-spADMM for the online optimization. Section 4 is about the  complexity of the average iteration of spADMM for solving Problem (\ref{sconvex0}) and the recovery of an important inequality in \cite{HSZhang2018}. 
{Section 5 applies Online-spADMM to online quadratic optimization, lasso and generalized Lasso and present experimental results.} We make our final conclusions in 
{Section 6}.

\section{Key Inequalities of Online spADMM}\label{Mainresluts}
\setcounter{equation}{0}
\quad \, In this section, we demonstrate important inequalities for upper bounds of $[f_k(x^{k+1})+g(z^{k+1})]-[f_k(\widehat x)+g(\widehat z)]$, where $(\widehat x,\widehat z)$ is a feasible point of the set $\Phi$. These  upper bounds play crucial roles in the analysis for constraint violation regret and objective regret of Online spADMM.
\begin{theorem}\label{inequality}Let $\{(x^k,z^k,y^k):k=1,\ldots, N+1\}$ be generated by Online-spADMM. Then, for any $(\widehat x,\widehat z) \in \Phi$, any $k=1,\ldots, N$,
$$
\begin{array}{l}
[f_k(x^{k+1})+g(z^{k+1})]-[f_k(\widehat x)+g(\widehat z)]\\[6pt]
\leq \left [ (2\sigma \tau)^{-1}\|y^k\|^2+\displaystyle \frac{1}{2}\|x^k- \widehat x\|_{\sigma\S_k}^2+\displaystyle \frac{1}{2}\|z^k -\widehat z\|_{\T}^2+
\displaystyle \frac{\sigma}{2}\|B(z^k-\widehat z)\|^2\right.\\[6pt]
\quad +\displaystyle \frac{1}{2}\left.\|z^k-z^{k-1}\|^2_{\T}+\left(1-\min\{\tau,\tau^{-1}\}\right)\displaystyle \frac{\sigma}{2}\|Ax^k+Bz^k-c\|^2 \right ]\\[10pt]
-\left [ (2\sigma \tau)^{-1}\|y^{k+1}\|^2+\displaystyle \frac{1}{2}\|x^{k+1}-\widehat x\|_{\sigma\S_k}^2+\displaystyle \frac{1}{2}\|z^{k+1}-\widehat z\|_{\T}^2+
\displaystyle \frac{\sigma}{2}\|B(z^{k+1}-\widehat z\|^2\right.\\[6pt]
\end{array}
$$
\begin{equation}\label{eq:mainineq}
\begin{array}{l}
\quad +\displaystyle \left. \frac{1}{2}\|z^{k+1}-z^{k}\|^2_{\T}+\left(1-\min\{\tau,\tau^{-1}\}\right)\displaystyle \frac{\sigma}{2}\|Ax^{k+1}+Bz^{k+1}-c\|^2 \right ]\\[10pt]
-\left[\displaystyle \frac{1}{2}\tau(1-\tau+\min (\tau,\tau^{-1}))\sigma \|B(z^{k+1}-z^{k})\|^2+\displaystyle \frac{1}{2}\|z^{k+1}-z^{k}\|_{\T}^2\right.\\[16pt]
\quad +\left. \displaystyle \frac{1}{2}\|x^{k+1}-x^{k}\|_{\sigma\S_{k}}^2+\displaystyle \frac{1}{2}\|x^{k+1}-\widehat x\|^2_{\Sigma_{f_k}}+\displaystyle \frac{1}{2}\|z^{k+1}-\widehat z\|^2_{\Sigma_g}.\right.\\[10pt]
\quad +\left.\left(1-\tau+\min (\tau,\tau^{-1})\right)\displaystyle \frac{\sigma}{2}\|Ax^{k+1}+Bz^{k+1}-c\|^2\right].
\end{array}
\end{equation}
\end{theorem}
{\bf Proof.} The proof is quite lengthy. We put it in Appendix A. \hfill $\Box$\\
We define
$$
s_{\tau}:=\displaystyle \frac{1}{4}\Big[5-\tau-3\min \{\tau,\tau^{-1}\}\Big], t_{\tau}:=\displaystyle \frac{1}{2}\Big[1-\tau+\min \{\tau,\tau^{-1}\}\Big].
$$
Let $\overline E:{\cal X}\times {\cal Z} \rightarrow {\cal Y}$ be linear operator defined by $\overline E(x,z)=Ax+Bz$ for $(x,z)\in {\cal X}\times {\cal Z}$ and
\begin{equation}\label{eq:notationMH}
\begin{array}{l}
\overline {\cal M}_k=\mbox{Diag}\left(\sigma\S_k+\Sigma_{f_k}, \T+\Sigma_g+\sigma B^*B\right)+s_{\tau}\overline E^*\overline E,\\[6pt]
\overline {\cal H}_k=\mbox{Diag}\left(\sigma\S_k+\displaystyle \frac{1}{2}\Sigma_{f_k}, \T+\Sigma_g+2t_{\tau}\tau \sigma B^*B\right)+\displaystyle \frac{1}{4}t_{\tau}\sigma \overline E^*\overline E.
\end{array}
\end{equation}
\begin{proposition}\label{prop-equivalence}
Let $\tau \in (0,(1+\sqrt{5})/2)$, then
\begin{equation}\label{eq:a13}
\Sigma_{f_k}+\sigma\S_k+\sigma A^*A\succ 0\quad  \& \quad \Sigma_g+\T+\sigma B^*B\succ 0 \Longleftrightarrow \overline {\cal M}_k \succ 0 \Longleftrightarrow \overline  {\cal H}_k \succ 0.
\end{equation}
\end{proposition}
{\bf Proof}. Similar to the proof of Proposition 3 of \cite{HSZhang2018}. \hfill $\Box$
\begin{theorem}\label{inequalitya}Let $\{(x^k,z^k,y^k):k=1,\ldots, N+1\}$ be generated by Online-spADMM. Then, for any $(\widehat x,\widehat z) \in \Phi$, any $k=1,\ldots, N$,
\begin{equation}\label{eq:mainineqa}
\begin{array}{l}
[f_k(x^{k+1})+g(z^{k+1})]-[f_k(\widehat x)+g(\widehat z)]\\[6pt]
\leq \displaystyle \frac{1}{2}\left [ (\tau\sigma)^{-1}\|y^k\|^2+\|x^k -\widehat x\|_{\sigma\S_k}^2+\|z^k-\widehat z\|_{\T}^2+\sigma \|B(z^k-\widehat z)\|^2+\|z^k-z^{k-1}\|_{\T}^2\right.\\[6pt]
\quad \quad \left.+s_{\tau}\sigma\|A(x^k-\widehat x)+B(z^k-\widehat z)\|^2+\|x^k-\widehat x\|_{\Sigma_{f_k}}^2+\|z^k-\widehat z\|_{\Sigma_g}^2
 \right ]\\[10pt]
 -\displaystyle \frac{1}{2}\left [ (\tau\sigma)^{-1}\|y^{k+1}\|^2+\|x^{k+1} -\widehat x\|_{\sigma\S_k}^2+\|z^{k+1}-\widehat z\|_{\T}^2\right.\\[10pt]
\quad \quad \quad\quad \quad  \left.+\sigma \|B(z^{k+1}-\widehat z)\|^2+\|z^{k+1}-z^{k-1}\|_{\T}^2\right.\\[6pt]
\quad \quad \left.+s_{\tau}\sigma\|A(x^{k+1}-\widehat x)+B(z^{k+1}-\widehat z)\|^2+\|x^{k+1}-\widehat x\|_{\Sigma_{f_k}}^2+\|z^{k+1}-\widehat z\|_{\Sigma_g}^2
 \right ]\\[10pt]
-\displaystyle \frac{1}{2}\left[2t_{\tau}\sigma \tau \|B(z^{k+1}-z^k)\|^2+\|z^{k+1}-z^k\|_{\T}^2+\|x^{k+1}-x^k\|_{\sigma\S_k}^2\right.\\[10pt]
\quad \quad \quad  \left. +\displaystyle \frac{1}{2}\|x^{k+1}-x^k\|_{\Sigma_{f_k}}^2+\displaystyle \frac{1}{2}\|z^{k+1}-z^k\|_{\Sigma_g}^2+t_{\tau}\sigma \|Ax^{k+1}+Bz^{k+1}-c\|^2\right.\\[10pt]
\quad \quad \quad \quad \quad \quad \quad  \quad\quad \quad \quad \quad\left. +\displaystyle \frac{1}{4} t_{\tau} \sigma \|A(x^{k+1}-x^k)+B(z^{k+1}-z^k)\|^2
\right].
\end{array}
\end{equation}
\end{theorem}
{\bf Proof}. We have from Theorem \ref{inequality} that

\begin{equation}\label{eq:a14}
\begin{array}{l}
[f_k(x^{k+1})+g(z^{k+1})]-[f_k(\widehat x)+g(\widehat z)]\\[6pt]
\leq \displaystyle \frac{1}{2}\left [ (\tau\sigma)^{-1}\|y^k\|^2+\|x^k -\widehat x\|_{\sigma\S_k}^2+\|z^k-\widehat z\|_{\T}^2+\sigma \|B(z^k-\widehat z)\|^2+\|z^k-z^{k-1}\|_{\T}^2\right.\\[6pt]
\quad \quad \quad \quad\quad \quad  \left.+ (1-\min \{\tau,\tau^{-1}\})\sigma\|Ax^k+Bz^k-c\|^2
 \right ]\\[10pt]
-\displaystyle \frac{1}{2}\left [ (\tau\sigma)^{-1}\|y^{k+1}\|^2+\|x^{k+1} -\widehat x\|_{\sigma\S_k}^2+\|z^{k+1}-\widehat z\|_{\T}^2+\sigma \|B(z^{k+1}-\widehat z)\|^2+\|z^{k+1}-z^{k}\|_{\T}^2\right.\\[10pt]
\quad \quad \quad\quad \quad  \left.+ (1-\min \{\tau,\tau^{-1}\})\sigma\|Ax^{k+1}+Bz^{k+1}-c\|^2
 \right ]\\[10pt]
-\displaystyle \frac{1}{2}\left[(1-\min \{\tau,\tau^{-1}\})\sigma \tau \|B(z^{k+1}-z^k)\|^2+\|z^{k+1}-z^k\|_{\T}^2+\|x^{k+1}-x^k\|_{\sigma\S_k}^2\right.\\[10pt]
\quad \quad \quad  \left. +\|x^{k+1}-\widehat x\|_{\Sigma_{f_k}}^2+\|z^{k+1}-\widehat z\|_{\Sigma_g}^2\right.\\[10pt]
\quad \quad \quad \left. + \left(1-\tau+\min (\tau,\tau^{-1})\right)\sigma \|A(x^{k+1}-x^k)+B(z^{k+1}-z^k)\|^2
\right].
\end{array}
\end{equation}
Reorganizing the terms on (\ref{eq:a14}), we obtain
$$
\begin{array}{l}
2[f_k(x^{k+1})+g(z^{k+1})]-[f_k(\widehat x)+g(\widehat z)]\\[6pt]
\leq \left [ (\tau\sigma)^{-1}\|y^k\|^2+\|x^k -\widehat x\|_{\sigma\S_k}^2+\|z^k-\widehat z\|_{\T}^2+\sigma \|B(z^k-\widehat z)\|^2+\|z^k-z^{k-1}\|_{\T}^2\right.\\[6pt]
 \left.+\displaystyle \frac{1}{4} (5-\tau -3\min \{\tau,\tau^{-1}\})\sigma\|Ax^k+Bz^k-c\|^2+\|x^k-\widehat x\|_{\Sigma_{f_k}}^2+\|z^k-\widehat z\|_{\Sigma_g}^2
 \right ]\\[10pt]
-\left [ (\tau\sigma)^{-1}\|y^{k+1}\|^2+\|x^{k+1} -\widehat x\|_{\sigma\S_k}^2+\|z^{k+1}-\widehat z\|_{\T}^2+\sigma \|B(z^{k+1}-\widehat z)\|^2+\|z^{k+1}-z^{k}\|_{\T}^2\right.\\[6pt]
 \left.+\displaystyle \frac{1}{4} (5-\tau -3\min \{\tau,\tau^{-1}\})\sigma\|Ax^{k+1}+Bz^{k+1}-c\|^2+\|x^{k+1}-\widehat x\|_{\Sigma_{f_k}}^2+\|z^{k+1}-\widehat z\|_{\Sigma_g}^2
 \right ]\\[10pt]
-\left[2t_{\tau}\sigma \tau \|B(z^{k+1}-z^k)\|^2+\|z^{k+1}-z^k\|_{\T}^2+\|x^{k+1}-x^k\|_{\sigma\S_k}^2+\|x^{k+1}-\widehat x\|_{\Sigma_{f_k}}^2\right.\\[10pt]
\quad \quad \quad \left.  +\|z^{k+1}-\widehat z\|_{\Sigma_g}^2+ \displaystyle \frac{1}{2}\left(1-\tau+\min (\tau,\tau^{-1})\right)\sigma \|Ax^{k+1}+Bz^{k+1}-c\|^2 \right.\\[10pt]
\quad \left. + \displaystyle \frac{1}{4}\left(1-\tau+\min (\tau,\tau^{-1})\right)\sigma \left [ \|Ax^{k+1}+Bz^{k+1}-c\|^2+\|Ax^{k}+Bz^{k}-c\|^2 \right]
\right].
\end{array}
$$
Or equivalently
$$
\begin{array}{l}
2[f_k(x^{k+1})+g(z^{k+1})]-[f_k(\widehat x)+g(\widehat z)]\\[6pt]
\leq \left [ (\tau\sigma)^{-1}\|y^k\|^2+\|x^k -\widehat x\|_{\sigma\S_k}^2+\|z^k-\widehat z\|_{\T}^2+\sigma \|B(z^k-\widehat z)\|^2+\|z^k-z^{k-1}\|_{\T}^2\right.\\[6pt]
 \left.+\displaystyle \frac{1}{4} (5-\tau -3\min \{\tau,\tau^{-1}\})\sigma\|Ax^k+Bz^k-c\|^2+\|x^k-\widehat x\|_{\Sigma_{f_k}}^2+\|z^k-\widehat z\|_{\Sigma_g}^2
 \right ]\\[10pt]
-\left [ (\tau\sigma)^{-1}\|y^{k+1}\|^2+\|x^{k+1} -\widehat x\|_{\sigma\S_k}^2+\|z^{k+1}-\widehat z\|_{\T}^2+\sigma \|B(z^{k+1}-\widehat z)\|^2+\|z^{k+1}-z^{k}\|_{\T}^2\right.\\[6pt]
 \left.+s_{\tau}\sigma\|Ax^{k+1}+Bz^{k+1}-c\|^2+\|x^{k+1}-\widehat x\|_{\Sigma_{f_k}}^2+\|z^{k+1}-\widehat z\|_{\Sigma_g}^2
 \right ]\\[10pt]
 \end{array}
 $$
 \begin{equation}\label{eq:a15a}
\begin{array}{l}
-\left[2t_{\tau}\sigma \tau \|B(z^{k+1}-z^k)\|^2+\|z^{k+1}-z^k\|_{\T}^2+\|x^{k+1}-x^k\|_{\sigma\S_k}^2+\|x^{k+1}-\widehat x\|_{\Sigma_{f_k}}^2\right.\\[10pt]
\quad \quad \quad \left.  +\|z^{k+1}-\widehat z\|_{\Sigma_g}^2+ t_{\tau}\sigma \|Ax^{k+1}+Bz^{k+1}-c\|^2 \right.\\[10pt]
\quad \left. + \displaystyle \frac{1}{2}t_{\tau}\sigma \left [ \|Ax^{k+1}+Bz^{k+1}-c\|^2+\|Ax^{k}+Bz^{k}-c\|^2 \right]
\right].
\end{array}
\end{equation}
Using equalities
$$
Ax^{k+1}+Bz^{k+1}-c=A(x^{k+1}-\widehat x)+B(z^{k+1}-\widehat z),\quad Ax^{k}+Bz^{k}-c=A(x^{k}-\widehat x)+B(z^{k}-\widehat z)
$$
and inequalities
$$
\begin{array}{l}
\|x^{k+1}-\widehat x\|_{\Sigma_{f_k}}^2 +\|x^{k+1}-\widehat x\|_{\Sigma_{f_k}}^2\geq \displaystyle \frac{1}{2}\|x^{k+1}-x^k\|_{\Sigma_{f_k}}^2,\\[6pt]
\|z^{k+1}-\widehat z\|_{\Sigma_g}^2 +\|z^{k+1}-\widehat z\|_{\Sigma_g}^2\geq \displaystyle \frac{1}{2}\|z^{k+1}- z^k\|_{\Sigma_g}^2,\\[6pt]
\|Ax^{k+1}+Bz^{k+1}-c\|^2+\|Ax^{k}+Bz^{k}-c\|^2\geq \displaystyle \frac{1}{2}\|A(x^{k+1}-x^k)+B(z^{k+1}-z^k)\|^2,
\end{array}
$$
we obtain from (\ref{eq:a15a}),
\begin{equation}\label{eq:mainineqb}
\begin{array}{l}
2\left\{[f_k(x^{k+1})+g(z^{k+1})]-[f_k(\widehat x)+g(\widehat z)]\right\}\\[6pt]
\leq \left [ (\tau\sigma)^{-1}\|y^k\|^2+\|x^k -\widehat x\|_{\sigma\S_k}^2+\|z^k-\widehat z\|_{\T}^2+\sigma \|B(z^k-\widehat z)\|^2+\|z^k-z^{k-1}\|_{\T}^2\right.\\[6pt]
\quad \quad \left.+s_{\tau}\sigma\|A(x^k-\widehat x)+B(z^k-\widehat z)\|^2+\|x^k-\widehat x\|_{\Sigma_{f_k}}^2+\|z^k-\widehat z\|_{\Sigma_g}^2
 \right ]\\[10pt]
-\left [ (\tau\sigma)^{-1}\|y^{k+1}\|^2+\|x^{k+1} -\widehat x\|_{\sigma\S_k}^2+\|z^{k+1}-\widehat z\|_{\T}^2\right.\\[10pt]
\quad \quad \quad\quad \quad  \left.+\sigma \|B(z^{k+1}-\widehat z)\|^2+\|z^{k+1}-z^{k-1}\|_{\T}^2\right.\\[6pt]
\quad \quad \left.+s_{\tau}\sigma\|A(x^{k+1}-\widehat x)+B(z^{k+1}-\widehat z)\|^2+\|x^{k+1}-\widehat x\|_{\Sigma_{f_k}}^2+\|z^{k+1}-\widehat z\|_{\Sigma_g}^2
 \right ]\\[10pt]
-\left[2t_{\tau}\sigma \tau \|B(z^{k+1}-z^k)\|^2+\|z^{k+1}-z^k\|_{\T}^2+\|x^{k+1}-x^k\|_{\sigma\S_k}^2\right.\\[10pt]
\quad \quad \quad  \left. +\displaystyle \frac{1}{2}\|x^{k+1}-x^k\|_{\Sigma_{f_k}}^2+\displaystyle \frac{1}{2}\|z^{k+1}-z^k\|_{\Sigma_g}^2+t_{\tau}\sigma \|Ax^{k+1}+Bz^{k+1}-c\|^2\right.\\[10pt]
\quad \quad \quad \quad \quad \quad \quad  \quad\quad \quad \quad \quad\left. +\displaystyle \frac{1}{4} t_{\tau} \sigma \|A(x^{k+1}-x^k)+B(z^{k+1}-z^k)\|^2
\right].
\end{array}
\end{equation}
This is just inequality (\ref{eq:mainineqa}). \hfill $\Box$
\begin{corollary}\label{cora}
Let $\{(x^k,z^k,y^k):k=1,\ldots, N+1\}$ be generated by Online-spADMM. Then, for any $(\widehat x,\widehat z) \in \Phi$, any $k=1,\ldots, N$,
\begin{equation}\label{eq:ineqc}
\begin{array}{l}
[f_k(x^{k+1})+g(z^{k+1})]-[f_k(\widehat x)+g(\widehat z)]+\displaystyle \frac{1}{2}\sigma t_{\tau} \|Ax^{k+1}+Bz^{k+1}-c\|^2\\[10pt]
\leq \displaystyle \frac{1}{2}\left [ (\tau\sigma)^{-1}\|y^k\|^2+\|z^k-z^{k-1}\|_{\T}^2\right]- \displaystyle \frac{1}{2}\left [ (\tau\sigma)^{-1}\|y^{k+1}\|^2+\|z^{k+1}-z^{k}\|_{\T}^2\right]\\[10pt]
\quad \quad+\displaystyle \frac{1}{2}\left [ \|(x^k,z^k)-(\widehat x,\widehat z)\|_{\overline {\cal M}_k}^2-\|(x^{k+1},z^{k+1})-(\widehat x,\widehat z)\|_{\overline {\cal M}_k}^2\right]\\[10pt]
\quad \quad -\displaystyle \frac{1}{2}\left[\|(x^{k+1},z^{k+1})-(x^k,z^k)\|_{\overline {\cal H}_k}^2
\right].
\end{array}
\end{equation}
\end{corollary}
{\bf Proof}. The inequality (\ref{eq:ineqc}) follows from (\ref{eq:mainineqa}) and the definitions of $\overline{\cal M}_k$ and $\overline{\cal H}_k$. \hfill $\Box$
\section{Regret  Analysis}
\setcounter{equation}{0}
Let $S^*_N$ denote the solution set of Problem (\ref{AonlineP}) and for $(\overline x,\overline z) \in S^*_N$,
 $$\nu^*_N=\displaystyle \frac{1}{N}\sum_{k=1}^N \left[ f_k(\overline x)+g(\overline z)\right].$$
 In this section, we discuss the iteration complexity of Online-spADMM for solving Problem (\ref{onlineP}). For establishing the constraint violation regret and the objective regret of Online-spADMM in terms of round number $N$, we make the following assumptions.
\begin{assumption}\label{ass:value-bound}
Suppose that the sequence $\{(x^k,z^k):k=1,\ldots,N\}$ generated by Online-spADMM satisfies
$$
[f_k(\widehat x)+g(\widehat z)] -[f_k(x^{k})+g(z^{k})]\leq \gamma_0, \forall k=1,\ldots, N.
$$
for some $\gamma_0>0$ and $(\widehat x,\widehat z) \in S^*_N$.
\end{assumption}
\begin{assumption}\label{ass:matrixthta}
For any $k=1,\ldots, N$, assume that $\S_k$ satisfies
\begin{equation}\label{eq:matrixthtaSigma}
\sigma S_k+\displaystyle \frac{2\tau}{1+8\tau}t_{\tau}\sigma A^*A \succeq \displaystyle \frac{2\tau}{1+8\tau}t_{\tau}\sigma {\cal I},
\end{equation}
where ${\cal I}$ is the identity mapping of ${\cal Y}$.
\end{assumption}
\begin{assumption}\label{ass:matrixthtac}
For any $k=1,\ldots, N$ and  $g^k \in \partial f_k(x^k)$, suppose that there exists $L>0$ such that
\begin{equation}\label{eq:subgradN}
\displaystyle \frac{1}{N} \sum_{k=1}^N \|g^k\|^2\leq L^2.
\end{equation}
\end{assumption}
When $f_t$ is Lipschitz continuous with Lipschitz constant $L_t$ and
$$
\displaystyle \frac{1}{N} \sum_{k=1}^N L_k^2\leq L^2,
$$
then Assumption \ref{ass:matrixthtac} is satisfied.
 \subsection{Constraint and objective regrets}
 Define
\begin{equation}\label{notation:sigmamatrix}
\Theta_k(\sigma)=\displaystyle \frac{1}{2}\Sigma_{f_k}+\sigma S_k+\displaystyle \frac{2\tau}{1+8\tau}t_{\tau}\sigma A^*A.
\end{equation}
For this purpose, we give the following lemma.
\begin{lemma}\label{lem:Hkinequ}
For $\overline {\cal H}_k$ defined by (\ref{eq:notationMH}) and $(x,z) \in {\cal X}\times {\cal Z}$, one has
\begin{equation}\label{Mkinequality}
\|(x,z)\|^2_{\overline {\cal H}_k}\geq \|z\|^2_{\Sigma_g+\T} +\|x\|^2_{\Theta_k(\sigma)}.
\end{equation}
\end{lemma}
{\bf Proof}. We can express $\|(x,z)\|^2_{\overline {\cal H}_k}$ as
\begin{equation}\label{eq:norm2}
\begin{array}{ll}
\|(x,z)\|^2_{\overline {\cal H}_k} &= \left\langle \left(\sigma \S_k+\displaystyle \frac{1}{2}\Sigma_{f_k}\right)x,x\right\rangle
+\displaystyle \frac{1}{4}t_{\tau}\sigma \|Ax+Bz\|^2\\[12pt]
& \quad \,  +
\Big\langle (\T+\Sigma_g+2t_{\tau}\tau \sigma B^*B)z, z\Big\rangle\\[12pt]
&=\left\langle \left(\displaystyle \frac{1}{2}\Sigma_{f_k}+\sigma \S_k+\displaystyle \frac{1}{4}t_{\tau}\sigma A^*A\right)x,x\right\rangle\\[12pt]
& \quad \, + \Big\langle \left(\Sigma_g+\T+ \displaystyle \frac{1}{4}t_{\tau}\sigma (1+8\tau)B^*B\right)z, z\Big\rangle+
\displaystyle \frac{1}{2} t_{\tau}\sigma \Big\langle Ax, Bz \Big \rangle.
\end{array}
\end{equation}
In view of the inequality
$$
\displaystyle \frac{1}{2} t_{\tau}\sigma \Big\langle Ax, Bz \Big \rangle\geq
-\displaystyle \frac{1}{4}t_{\tau}\sigma\left[(1+8\tau)^{-1}\|Ax\|^2+(1+8\tau)\|Bz\|^2\right],
$$
we have from (\ref{eq:norm2}) that
$$
\|(x,z)\|^2_{\overline {\cal H}_k}\geq  \Big\langle \left(\T+\Sigma_g\right)z, z\Big\rangle
 +\left\langle \left(\displaystyle \frac{1}{2}\Sigma_{f_k}+\sigma \S_k+\displaystyle \frac{2\tau}{1+8\tau}t_{\tau}\sigma A^*A\right)x,x\right\rangle.
$$
The proof is completed. \hfill $\Box$

Making a summation of (\ref{eq:ineqc}) for $k=1,\ldots, N$, we obtain
\begin{equation}\label{eq:a16}
\begin{array}{l}
\displaystyle \sum_{k=1}^N\left\{[f_k(x^{k+1})+g(z^{k+1})]-[f_k(\widehat x)+g(\widehat z)]+\displaystyle \frac{1}{2} t_{\tau}\sigma \|Ax^{k+1}+Bz^{k+1}-c\|^2\right\}\\[16pt]
\leq \displaystyle \frac{1}{2}\left [ (\tau\sigma)^{-1}\|y^1\|^2+\|z^1-z^{0}\|_{\T}^2\right]- \displaystyle \frac{1}{2}\left [ (\tau\sigma)^{-1}\|y^{N+1}\|^2+\|z^{N+1}-z^{N}\|_{\T}^2\right]\\[10pt]
\quad \quad+\displaystyle \sum_{k=1}^N\left\{\displaystyle \frac{1}{2}\left [ \|(x^k,z^k)-(\widehat x,\widehat z)\|_{\overline {\cal M}_k}^2-\|(x^{k+1},z^{k+1})-(\widehat x,\widehat z)\|_{\overline {\cal M}_k}^2\right]\right\}\\[10pt]
\quad \quad -\displaystyle \sum_{k=1}^N\left\{\displaystyle \frac{1}{2}\left[\|(x^{k+1},z^{k+1})-(x^k,z^k)\|_{\overline {\cal H}_k}^2
\right]\right\}.
\end{array}
\end{equation}
We first give a proposition about a bound for the constraint violation regret by Online-spADMM.
\begin{proposition}\label{Citer-complexity}
Let $\{(x^k,z^k,y^k):k=1,\ldots, N+1\}$ be generated by Online-spADMM. Let the following matrix orders hold:
\begin{equation}\label{eq:orders}
\S_1 \succeq \S_2 \succeq \cdots \succeq \S_{N+1}.
\end{equation}
Then
\begin{equation}\label{eq:c-compexity}
\begin{array}{ll}
\displaystyle  \frac{1}{N} {\rm regret}^{\rm ctr}_N & =\displaystyle  \frac{1}{N}\sum_{k=1}^N\left\{[\displaystyle   \|Ax^{k+1}+Bz^{k+1}-c\|^2\right\}\\[16pt]
& \leq \displaystyle \frac{1}{Nt_{\tau}\sigma}
\left [ (\tau\sigma)^{-1}\|y^1\|^2+\|z^1-z^{0}\|_{\T}^2\right]
+\displaystyle \displaystyle \frac{1}{Nt_{\tau}\sigma}\left[{\rm Dist}_{\overline{\cal M}_1}((x^1,z^1),S^*_N)\right]^2 \\[12pt]
& \quad +\displaystyle \frac{2}{t_{\tau}\sigma}\left[\nu^*_N -\displaystyle\frac{1}{N} \sum_{k=1}^N[f_k(x^{k+1})+g(z^{k+1})]\right].
\end{array}
\end{equation}
\end{proposition}
{\bf Proof}. For any $(\widehat x,\widehat z) \in S^*_N$, we have
$$
\displaystyle \frac{1}{N}\displaystyle \sum_{k=1}^N \left[f_k(\widehat x)+g(\widehat z)\right]=\nu^*_N.
$$
Since $\S_1 \succeq \S_2 \succeq \cdots \succeq \S_N$, we have
$$
\overline {\cal M}_1 \succeq \overline {\cal M}_2 \succeq \cdots \succeq \overline {\cal M}_{N+1}.
$$
Then, for $\{w^1, w^2,\ldots, w^{N+1}\}\subset {\cal X} \times {\cal Z}$, one has
$$
\|w^2\|_{\overline{\cal M}_1} \geq \|w^2\|_{\overline{\cal M}_2}\geq  \ldots \geq\|w^{N+1}\|_{\overline{\cal M}_N} \geq \|w^{N+1}\|_{\overline{\cal M}_{N+1}}.
$$
Thus we obtain
$$
\begin{array}{l}
\displaystyle \sum_{k=1}^N \left [ \|(x^k,z^k)-(\widehat x,\widehat z)\|_{\overline{\cal M}_k}^2-\|(x^{k+1},z^{k+1})-(\widehat x,\widehat z)\|_{\overline{\cal M}_k}^2\right]\\[16pt]
\quad \quad =\|(x^1,z^1)-(\widehat x,\widehat z)\|_{\overline{\cal M}_1}^2- \left [ \|(x^2,z^2)-(\widehat x,\widehat z)\|_{\overline{\cal M}_1}^2-\|(x^2,z^2-(\widehat x,\widehat z)\|_{\overline{\cal M}_2}^2\right]\\[20pt]
\quad \quad - \left [ \|(x^3,z^3)-(\widehat x,\widehat z)\|_{\overline{\cal M}_2}^2-\|(x^3,z^3)-(\widehat x,\widehat z)\|_{\overline{\cal M}_3}^2\right]\\[6pt]
\quad \quad - \quad\quad\quad  \cdots \\[6pt]
\quad \quad - \left [ \|(x^N,z^N)-(\widehat x,\widehat z)\|_{\overline{\cal M}_N}^2-\|(x^N,z^N)-(\widehat x,\widehat z)\|_{\overline{\cal M}_{N+1}}^2\right]\\[10pt]
\quad \quad -\|(x^{N+1},z^{N+1})-(\widehat x,\widehat z)\|_{\overline{\cal M}_{N+1}}^2\\[10pt]
\leq \|(x^1,z^1)-(\widehat x,\widehat z)\|_{\overline{\cal M}_1}^2-\|(x^{N+1},z^{N+1})-(\widehat x,\widehat z)\|_{\overline{\cal M}_{N+1}}^2.
\end{array}
$$
Therefore we  have from (\ref{eq:a16}) that
$$
 \begin{array}{l}
\displaystyle \frac{1}{N}\sum_{k=1}^N\left\{[f_k(x^{k+1})+g(z^{k+1})]\right\}-\nu^*_N+
\displaystyle \frac{t_{\tau}\sigma}{2N}\displaystyle \sum_{k=1}^N  \|Ax^{k+1}+Bz^{k+1}-c\|^2\\[12pt]
\leq \displaystyle \frac{1}{2N}\left [ (\tau\sigma)^{-1}\|y^1\|^2+\|z^1-z^{0}\|_{\T}^2\right]
+\displaystyle \displaystyle \frac{1}{2N}\|(x^1,z^1)-(\widehat x,\widehat z)\|_{\overline{\cal M}_1}^2.
\end{array}
$$
Therefore,  inequality (\ref{eq:c-compexity}) follows directly from inequality (\ref{eq:a16}) and Assumption \ref{ass:value-bound}. \hfill $\Box$\\

We now  give  the following proposition about a bound for the sum of objective regret and constraint violation regret by Online-spADMM.
\begin{proposition}\label{Oiter-complexity}
Let $\{(x^k,z^k,y^k):k=1,\ldots, N+1\}$ be generated by Online-spADMM. Let Assumption \ref{ass:value-bound} be satisfied  and the following matrix orders hold:
\begin{equation}\label{eq:ordersa}
\S_1 \succeq \S_2 \succeq \cdots \succeq \S_{N+1}.
\end{equation}
 Then the following  property holds:
\begin{equation}\label{eq:objective-compexity}
\begin{array}{l}
\displaystyle \frac{1}{N}\sum_{k=}^{N}\left\{[f_k(x^{k})+g(z^{k})]\right\}-\nu^*_N+
\displaystyle \frac{t_{\tau}\sigma}{2N}\displaystyle \sum_{k=1}^N  \|Ax^{k+1}+Bz^{k+1}-c\|^2\\[12pt]
\leq \displaystyle \frac{1}{2N}\left [ (\tau\sigma)^{-1}\|y^1\|^2+\|z^1-z^{0}\|_{\T}^2+2(g(z^1)-g(z^{N+1}))\right]\\[12pt]
\quad+\displaystyle \displaystyle \frac{1}{2N}\|(x^1,z^1)-(\widehat x,\widehat z)\|_{\overline{\cal M}_1}^2
 + \displaystyle \frac{1}{N}\displaystyle \sum_{k=1}^N \|g^k\|^2_{\Theta_k(\sigma)^{-1}}
\end{array}
\end{equation}
for $g^k \in \partial f_k(x^k), k=1,\ldots, N$.
\end{proposition}
{\bf Proof}. For $k=1,\ldots, N$, we have for any $g^k \in \partial f_k(x^k)$,
\begin{equation}\label{eq:b1}
f_k(x^k)-f_k(x^{k+1}) \leq \left \langle g^k, x^k-x^{k+1} \right \rangle \leq \displaystyle \frac{1}{2}
\left ( \|g^k\|^2_{\Theta_k(\sigma)^{-1}}+ \|x^k-x^{k+1}\|^2_{\Theta_k(\sigma)}\right).
\end{equation}
It follows from Lemma \ref{lem:Hkinequ} that
\begin{equation}\label{eq:b2}
\|x^{k+1}-x^k\|^2_{\Theta_k(\sigma)} \leq \|(x^{k+1}-x^k,z^{k+1}-z^k)\|^2_{\overline {\cal H}_k}- \|z^{k+1}-z^k\|^2_{\Sigma_g+\T}.
\end{equation}
We have from (\ref{eq:a16}),(\ref{eq:b1}) and  (\ref{eq:b2}) that
\begin{equation}\label{eq:a161}
\begin{array}{l}
\displaystyle \sum_{k=1}^N
\left\{[f_{k}(x^{k})+g(z^{k})]-[f_{k}(\widehat x)+g(\widehat z)]\right\}+\displaystyle \frac{1}{2} t_{\tau}\sigma \displaystyle \sum_{k=1}^N\left\{ \|Ax^{k+1}+Bz^{k+1}-c\|^2\right\}\\[16pt]
\leq \displaystyle \frac{1}{2}\left [ (\tau\sigma)^{-1}\|y^1\|^2+\|z^1-z^{0}\|_{\T}^2\right]- \displaystyle \frac{1}{2}\left [ (\tau\sigma)^{-1}\|y^{N+1}\|^2+\|z^{N+1}-z^{N}\|_{\T}^2\right]\\[10pt]
\quad \quad+\displaystyle \sum_{k=1}^N\left\{\displaystyle \frac{1}{2}\left [ \|(x^k,z^k)-(\widehat x,\widehat z)\|_{\overline {\cal M}_k}^2-\|(x^{k+1},z^{k+1})-(\widehat x,\widehat z)\|^2_{\overline {\cal M}_k}\right]\right\}\\[10pt]
\quad \quad -\displaystyle \sum_{k=1}^N \left\{ \displaystyle \frac{1}{2}
\left[\|(x^{k+1},z^{k+1})-(x^k,z^k)\|^2_{\overline {\cal H}_k}
\right]\right\}\\[10pt]
\quad \quad + \displaystyle \sum_{k=1}^N\left\{[f_k(x^{k})-f_k(x^{k+1})]\right\}+g(z^1)-g(z^{N+1})\\[10pt]
\leq \displaystyle \frac{1}{2}\left [ (\tau\sigma)^{-1}\|y^1\|^2+\|z^1-z^{0}\|_{\T}^2\right]- \displaystyle \frac{1}{2}\left [ (\tau\sigma)^{-1}\|y^{N+1}\|^2+\|z^{N+1}-z^{N}\|_{\T}^2\right]\\[10pt]
\quad \quad+\displaystyle \sum_{k=1}^N\left\{\displaystyle \frac{1}{2}\left [ \|(x^k,z^k)-(\widehat x,\widehat z)\|_{\overline {\cal M}_k}^2-\|(x^{k+1},z^{k+1})-(\widehat x,\widehat z)\|^2_{\overline {\cal M}_k}\right]\right\}\\[10pt]
\quad \quad + \displaystyle \sum_{k=1}^N\displaystyle \frac{1}{2}
\left ( \|g^k\|^2_{\Theta_k(\sigma)^{-1}}- \|z^{k+1}-z^k\|^2_{\Sigma_g+\T} \right)+g(z^1)-g(z^{N+1}).
\end{array}
\end{equation}
Since $\S_1 \succeq \S_2 \succeq \cdots \succeq \S_N$, we have $\overline {\cal M}_1 \succeq \cdots \succeq \overline {\cal M}_N$ and from (\ref{eq:a161}) we obtain
\begin{equation}\label{eq:objective-compexity}
\begin{array}{l}
\displaystyle \sum_{k=1}^{N}\left\{[f_k(x^{k})+g(z^{k})]\right\}-N\nu^*_N+
\displaystyle \frac{t_{\tau}\sigma}{2}\displaystyle \sum_{k=1}^N  \|Ax^{k+1}+Bz^{k+1}-c\|^2\\[12pt]
\leq \displaystyle \frac{1}{2}\left [ (\tau\sigma)^{-1}\|y^1\|^2+\|z^1-z^{0}\|_{\T}^2\right]
+\displaystyle \displaystyle \frac{1}{2}\|(x^1,z^1)-(\widehat x,\widehat z)\|_{\overline{\cal M}_1}^2\\[12pt]
\quad + \displaystyle \sum_{k=1}^N \|g^k\|^2_{\Theta_k(\sigma)^{-1}}+g(z^1)-g(z^{N+1}).
\end{array}
\end{equation}
The proof is completed. \hfill $\Box$\\

Now we are in a position to state the main theorem about  the bounds of the constraint violation regret and the objective regret by Online-spADMM, respectively.
\begin{theorem}\label{iter-complexity}
Let $\{(x^k,z^k,y^k):k=1,\ldots, N+1\}$ be generated by Online-spADMM. Suppose that
the following matrix orders hold:
\begin{equation}\label{eq:ordersb}
\S_1 \succeq \S_2 \succeq \cdots \succeq \S_{N+1}.
\end{equation}
Then the following properties hold:
\begin{itemize}
\item[(i)]  Let Assumptions \ref{ass:matrixthta} and \ref{ass:matrixthtac} be satisfied, then
\begin{equation}\label{eq:objective-compexityN}
\begin{array}{ll}
\displaystyle  \frac{1}{N} {\rm regret}^{\rm obj}_N &=\displaystyle \frac{1}{N}\displaystyle \sum_{k=1}^{N}\left\{[f_k(x^{k})+g(z^{k})]\right\}-\nu^*_N\\[12pt]
&\leq \displaystyle \frac{1}{2N}\left [ (\tau\sigma)^{-1}\|y^1\|^2+\|z^1-z^{0}\|_{\T}^2\right]
+\displaystyle \displaystyle \frac{1}{2N}\left[{\rm Dist}_{\overline{\cal M}_1}((x^1,z^1),S^*_N)\right]^2\\[12pt]
& \quad + \displaystyle \frac{(1+8\tau)}{2t_{\tau}\tau\sigma }L^2+\displaystyle \frac{1}{N}\Big(g(z^1)-g(z^{N+1})\Big).
\end{array}
\end{equation}
\item[(ii)]
Let Assumptions \ref{ass:value-bound}, \ref{ass:matrixthta}   and \ref{ass:matrixthtac} be satisfied, then
\begin{equation}\label{eq:constraint-compexity}
\begin{array}{ll}
\displaystyle  \frac{1}{N} {\rm regret}^{\rm ctr}_N & =\displaystyle  \frac{1}{N}\sum_{k=1}^N\left\{\displaystyle   \|Ax^{k+1}+Bz^{k+1}-c\|^2\right\}\\[10pt]
& \leq \displaystyle \frac{1}{Nt_{\tau}\sigma}
\left [ (\tau\sigma)^{-1}\|y^1\|^2+\|z^1-z^{0}\|_{\T}^2\right]\\[10pt]
& \quad \,
+\displaystyle \displaystyle \frac{1}{Nt_{\tau}\sigma}\left[{\rm Dist}_{\overline{\cal M}_1}((x^1,z^1),S^*_N)\right]^2 +\displaystyle \frac{2\gamma_0}{t_{\tau}\sigma}\\[10pt]
& \quad + \displaystyle \frac{(1+8\tau)}{t_{\tau}\tau\sigma }L^2+\displaystyle \frac{2}{t_{\tau}\sigma N}\Big(g(z^1)-g(z^{N+1})\Big).
\end{array}
\end{equation}
\end{itemize}
\end{theorem}
{\bf Proof}. From (\ref{eq:matrixthtaSigma}) in Assumption \ref{ass:matrixthta}, we know that
\begin{equation}\label{thetaK}
\Theta_k(\sigma)^{-1} \preceq \displaystyle \frac{(1+8\tau)}{2t_{\tau}\tau\sigma }{\cal I},
\end{equation}
and Conclusion (i) follows from Proposition \ref{Oiter-complexity} and Assumption \ref{ass:matrixthtac}. In view of (\ref{thetaK}), Assumption \ref{ass:value-bound} and  Assumption \ref{ass:matrixthtac}, we obtain (ii) from Proposition \ref{Oiter-complexity}. \hfill $\Box$\\
Define
$$
\begin{array}{l}
\eta(\tau)=\displaystyle \frac{1+8\tau}{2t_{\tau}\tau},\\[8pt]
\kappa_1=\displaystyle \frac{1}{2}(\|\Sigma_{f_1}\|+\|\T+\Sigma_g\|+s_{\tau}\|\overline E^*\overline E\|),\\[8pt]
\kappa_2=\displaystyle \frac{1}{2}(\|\S_1\|+\|B^*B\|),\quad \kappa_3=\displaystyle \frac{1}{2}[\tau^{-1}\|y^1\|^2+\|z^1-z^0\|^2]+g(z^1).\\[8pt]
\end{array}
$$
We obtain the following result from Theorem \ref{iter-complexity}, which provides upper bounds in terms of $N$ for objective regret and constraint violation regret.
\begin{corollary}\label{coro:N-complexity}
Let $\{(x^k,z^k,y^k):k=1,\ldots, N+1\}$ be generated by Online-spADMM. Suppose that
the following matrix orders hold:
\begin{equation}\label{eq:ordersbz}
\S_1 \succeq \S_2 \succeq \cdots \succeq \S_{N+1}.
\end{equation}
Let $\sigma=\sqrt{N}$.
Then the following properties hold:
\begin{itemize}
\item[(i)]  Let Assumptions \ref{ass:matrixthta} and \ref{ass:matrixthtac} be satisfied, then
\begin{equation}\label{eq:objective-compexityNN}
\begin{array}{ll}
\displaystyle  \frac{1}{N} {\rm regret}^{\rm obj}_N &=\displaystyle \frac{1}{N}\displaystyle \sum_{k=1}^{N}\left\{[f_k(x^{k})+g(z^{k})]\right\}-\nu^*_N\\[10pt]
& \leq \left ( \displaystyle \frac{\kappa_1}{N}+\displaystyle \frac{\kappa_2}{\sqrt{N}}\right )\left[{\rm Dist}((x^1,z^1),S^*_N)\right]^2\\[12pt]
& \quad \quad + \displaystyle \frac{1}{\sqrt{N}}\eta(\tau) L^2+\displaystyle \frac{1}{N}\Big(\kappa_3-g(z^{N+1})\Big).
\end{array}
\end{equation}
\item[(ii)]
Let Assumptions \ref{ass:value-bound}, \ref{ass:matrixthta}   and \ref{ass:matrixthtac} be satisfied, then
\begin{equation}\label{eq:constraint-compexityNN}
\begin{array}{ll}
\displaystyle  \frac{1}{N} {\rm regret}^{\rm ctr}_N &=\displaystyle  \frac{1}{N}\sum_{k=1}^N\left\{\displaystyle   \|Ax^{k+1}+Bz^{k+1}-c\|^2\right\}\\[10pt]
& \leq \displaystyle \frac{1}{t_{\tau}\sqrt{N}} \left ( \displaystyle \frac{\kappa_1}{N}+\displaystyle \frac{\kappa_2}{\sqrt{N}}\right )\left[{\rm Dist}((x^1,z^1),S^*_N)\right]^2\\[12pt]
& \quad \quad + \displaystyle \frac{1}{t_{\tau}N}\eta(\tau) L^2+\displaystyle \frac{1}{N^{3/2}t_{\tau}}\Big(\kappa_3-g(z^{N+1})\Big)
+\displaystyle \frac{2\gamma_0}{t_{\tau}\sqrt{N}}
\end{array}
\end{equation}
\end{itemize}
\end{corollary}
{\bf Proof}. From the definitions of $\overline{\cal M}_1$,$\kappa_1$ and $\kappa_2$, one has that
$$
\begin{array}{l}
\displaystyle \frac{1}{2N} \left[{\rm Dist}_{\overline{\cal M}_1}((x^1,z^1),S^*_N)\right]^2\\[10pt]
\leq\displaystyle \frac{1}{2N} \left [ \|\Sigma_{f_1}\|+\|\T+\Sigma_g\|+s_{\tau} \|\overline E^*\overline E\|+\sqrt{N}(\|\S_1\|+\|B^*B\|) \right]\left[{\rm Dist}((x^1,z^1),S^*_N)\right]^2\\[10pt]
\leq  \left ( \displaystyle \frac{\kappa_1}{N}+\displaystyle \frac{\kappa_2}{\sqrt{N}}\right )\left[{\rm Dist}((x^1,z^1),S^*_N)\right]^2.
\end{array}
$$
Thus we obtain  conclusion (i) from (i) of Theorem \ref{iter-complexity} and the definition of $\kappa_3$ and $\kappa_4$.
Conclusion (ii) follows from  (ii) of Theorem \ref{iter-complexity}. \hfill $\Box$\\
\begin{remark}\label{rem:severalpoints}
We have the following observations:
\begin{itemize}
\item[a)]
It follows from Corollary \ref{coro:N-complexity} that the objective regret and the constraint violation regret are of order $\sqrt{N}$ whose coefficients are dependent on the distance of the initial point $(x^1,z^1)$ from $S^*_N$ if $g$ is nonnegative-valued. \item[b)] Since $g$ is often a regularizer,  the assumption  $g$ being nonnegative-valued is quite natural. Assumption \ref{ass:matrixthtac} is a natural assumption, which was adopted by \cite{Shai2011}.
\item[c)] Assumption \ref{ass:value-bound} is satisfied in many circumstances.  For example if $g$ has the following form
$$
g(z)=\theta (z)+\delta_{Z}(z),
$$
where $Z \subset {\cal Z}$ is a nonempty convex compact set and $\theta: Z:\rightarrow \Re$ is a continuous function and $A^*$ is an onto linear operator, then  Assumption \ref{ass:value-bound} is satisfied.
\item[d)] Since we can not neglect the term
$$
\left[{\rm Dist}_{\overline{\cal M}_1}((x^1,z^1),S^*_N)\right]^2,
$$
which is of order $\sigma$, in (\ref{eq:objective-compexityN}) and (\ref{eq:constraint-compexity}) for objective regret and constraint violation regret, the optimal choice of $\sigma$ is of order ${\rm O}(\sqrt{N})$.
\end{itemize}
\end{remark}
\subsection{Solution regret}
In this subsection, we discuss the possibility for deriving solution regret. For this purpose, we define
  $$
\begin{array}{l}
{\cal W}_k=\mbox{Diag} (\Sigma_{f_k}, \Sigma_{g})+\displaystyle \frac{1}{2} t_{\tau}\sigma E^*E,\\[6pt]
\widehat {\cal M}_k=\mbox{Diag} (\sigma\S_k, \T +\sigma BB^*)+(1-\min (\tau,\tau^{-1}))\sigma E^*E,\\[6pt]
\widehat {\cal H}_k=\mbox{Diag} (\sigma\S_k, \T +2\sigma t_{\tau} BB^*)
\end{array}
$$
 and  derive new inequalities from (\ref{eq:mainineq}). From (\ref{eq:mainineq}), we have
\begin{equation}\label{eq:a15}
\begin{array}{l}
[f_k(x^{k+1})+g(z^{k+1})]-[f_k(\widehat x)+g(\widehat z)]\\[6pt]
\quad +\displaystyle \frac{1}{2}\|x^{k+1}-\widehat x\|^2_{\Sigma_{f_k}}+\displaystyle \frac{1}{2}\|z^{k+1}-\widehat z\|^2_{\Sigma_g}
+\displaystyle \frac{1}{2}\sigma t_{\tau}\|Ax^{k+1}+Bz^{k+1}-c\|^2\\[6pt]
\leq \left [ (2\sigma \tau)^{-1}\|y^k\|^2+\displaystyle \frac{1}{2}\|x^k- \widehat x\|_{\S_k}^2+\displaystyle \frac{1}{2}\|z^k -\widehat z\|_{\T}^2+
\displaystyle \frac{\sigma}{2}\|B(z^k-\widehat z)\|^2\right.\\[6pt]
\quad +\displaystyle \frac{1}{2}\left.\|z^k-z^{k-1}\|^2_{\T}+\left(1-\min\{\tau,\tau^{-1}\}\right)\displaystyle \frac{\sigma}{2}\|Ax^k+Bz^k-c\|^2 \right ]\\[10pt]
-\left [ (2\sigma \tau)^{-1}\|y^{k+1}\|^2+\displaystyle \frac{1}{2}\|x^{k+1}-\widehat x\|_{\S_k}^2+\displaystyle \frac{1}{2}\|z^{k+1}-\widehat z\|_{\T}^2+
\displaystyle \frac{\sigma}{2}\|B(z^{k+1}-\widehat z\|^2\right.\\[6pt]
\quad +\displaystyle \left. \frac{1}{2}\|z^{k+1}-z^{k}\|^2_{\T}+\left(1-\min\{\tau,\tau^{-1}\}\right)\displaystyle \frac{\sigma}{2}\|Ax^{k+1}+Bz^{k+1}-c\|^2 \right ]\\[10pt]
-\left[\displaystyle \frac{1}{2}\tau(1-\tau+\min (\tau,\tau^{-1}))\sigma \|B(z^{k+1}-z^{k})\|^2+\displaystyle \frac{1}{2}\|z^{k+1}-z^{k}\|_{\T}^2\right.\\[16pt]
\quad +\left. \displaystyle \frac{1}{2}\|x^{k+1}-x^{k}\|_{\S_{k}}^2+\displaystyle \frac{1}{2}\sigma t_{\tau}\|Ax^{k+1}+Bz^{k+1}-c\|^2\right].
\end{array}
\end{equation}

Then inequality (\ref{eq:a15}) can be equivalently written as
\begin{equation}\label{eq:a16a}
\begin{array}{l}
[f_k(x^{k+1})+g(z^{k+1})]-[f_k(\widehat x)+g(\widehat z)]\\[6pt]
\quad +\displaystyle \frac{1}{2}\|(x^{k+1},z^{k+1})-(\widehat x,\widehat z)\|_{{\cal W}_k}^2
+\displaystyle \frac{1}{2}\sigma t_{\tau}\|Ax^{k+1}+Bz^{k+1}-c\|^2\\[6pt]
\leq \left [ (2\sigma \tau)^{-1}\|y^k\|^2+\displaystyle \frac{1}{2}\|x^k- \widehat x\|_{\sigma\S_k}^2+\displaystyle \frac{1}{2}\|z^k -\widehat z\|_{\T}^2+
\displaystyle \frac{\sigma}{2}\|B(z^k-\widehat z)\|^2\right.\\[6pt]
\quad +\displaystyle \frac{1}{2}\left.\|z^k-z^{k-1}\|^2_{\T}+\left(1-\min\{\tau,\tau^{-1}\}\right)\displaystyle \frac{\sigma}{2}\|Ax^k+Bz^k-c\|^2 \right ]\\[10pt]
\quad-\left [ (2\sigma \tau)^{-1}\|y^{k+1}\|^2+\displaystyle \frac{1}{2}\|x^{k+1}-\widehat x\|_{\sigma\S_k}^2+\displaystyle \frac{1}{2}\|z^{k+1}-\widehat z\|_{\T}^2+
\displaystyle \frac{\sigma}{2}\|B(z^{k+1}-\widehat z\|^2\right.\\[6pt]
\quad +\displaystyle \left. \frac{1}{2}\|z^{k+1}-z^{k}\|^2_{\T}+\left(1-\min\{\tau,\tau^{-1}\}\right)\displaystyle \frac{\sigma}{2}\|Ax^{k+1}+Bz^{k+1}-c\|^2 \right ]\\[10pt]
\quad-\left[\tau\sigma t_{\tau} \|B(z^{k+1}-z^{k})\|^2+\displaystyle \frac{1}{2}\|z^{k+1}-z^{k}\|_{\T}^2
+\displaystyle \frac{1}{2}\|x^{k+1}-x^{k}\|_{\S_{k}}^2\right]\\[10pt]
=\left [ (2\sigma \tau)^{-1}\|y^k\|^2+\displaystyle \frac{1}{2}\|(x^k,z^k)- (\widehat x,\widehat z)\|_{\widehat {\cal M}_k}^2+\displaystyle \frac{1}{2}\|z^k-z^{k-1}\|^2_{\T} \right]\\[10pt]
\quad-\left [ (2\sigma \tau)^{-1}\|y^{k+1}\|^2+\displaystyle \frac{1}{2}\|(x^{k+1},z^{k+1})- (\widehat x,\widehat z)\|_{\widehat {\cal M}_k}^2+\displaystyle \frac{1}{2}\|z^{k+1}-z^{k}\|^2_{\T} \right]\\[10pt]
\quad-\displaystyle \frac{1}{2}\|(x^{k+1}, z^{k+1})-(x^{k},z^k)\|_{\widehat {\cal H}_k}^2.
\end{array}
\end{equation}
Like Proposition \ref{Citer-complexity},  basing on (\ref{eq:a16a}), we may prove a bound for both constraint violation regret and solution regret by Online-spADMM.
\begin{proposition}\label{Citer-complexityaaa}
Let $\{(x^k,z^k,y^k):k=1,\ldots, N+1\}$ be generated by Online-spADMM. Let the following matrix orders hold:
\begin{equation}\label{eq:orders}
\S_1 \succeq \S_2 \succeq \cdots \succeq \S_{N+1}.
\end{equation}
Then
\begin{equation}\label{eq:objective-compexityaa}
\begin{array}{l}
\displaystyle  \frac{t_{\tau}\sigma}{N}\sum_{k=1}^N\left\{[\displaystyle   \|Ax^{k+1}+Bz^{k+1}-c\|^2\right\}+\displaystyle \frac{1}{N}\sum_{k=1}^N\|(x^{k+1},z^{k+1})-(\widehat x,\widehat z)\|_{{\cal W}_k}^2  \\[16pt]
\leq \displaystyle \frac{1}{N}
\left [ (\tau\sigma)^{-1}\|y^1\|^2+\|z^1-z^{0}\|_{\T}^2\right]
+\displaystyle \displaystyle \frac{1}{N}\left[{\rm Dist}_{\widehat{\cal M}_1}((x^1,z^1),S^*_N)\right]^2 \\[12pt]
\quad +\displaystyle 2\left[\nu^*_N -\displaystyle\frac{1}{N} \sum_{k=1}^N[f_k(x^{k+1})+g(z^{k+1})]\right].
\end{array}
\end{equation}
\end{proposition}
Like Proposition \ref{Oiter-complexity}, we may prove another  bound of the objective regret by Online-spADMM.
\begin{proposition}\label{Oiter-complexitybbbb}
Let $\{(x^k,z^k,y^k):k=1,\ldots, N+1\}$ be generated by Online-spADMM. Let Assumption \ref{ass:value-bound} be satisfied  and the following matrix orders hold:
\begin{equation}\label{eq:ordersabbb}
\S_1 \succeq \S_2 \succeq \cdots \succeq \S_{N+1}\succ 0.
\end{equation}
 Then the following two property holds:
\begin{equation}\label{eq:objective-compexityaaa}
\begin{array}{l}
\displaystyle \frac{1}{N}\sum_{k=}^{N}\left\{[f_k(x^{k})+g(z^{k})]\right\}-\nu^*_N+
\displaystyle \frac{t_{\tau}\sigma}{2N}\displaystyle \sum_{k=1}^N  \|Ax^{k+1}+Bz^{k+1}-c\|^2\\[12pt]
\quad \,  +\displaystyle \frac{1}{2N}\sum_{k=1}^N\|(x^{k+1},z^{k+1})-(\widehat x,\widehat z)\|_{{\cal W}_k}^2  \\[16pt]
\leq \displaystyle \frac{1}{2N}\left [ (\tau\sigma)^{-1}\|y^1\|^2+\|z^1-z^{0}\|_{\T}^2+2(g(z^1)-g(z^{N+1})\right]\\[12pt]
\quad+\displaystyle \displaystyle \frac{1}{2N}\|(x^1,z^1)-(\widehat x,\widehat z)\|_{\widehat{\cal M}_1}^2
 + \displaystyle \frac{1}{\sigma N}\displaystyle \sum_{k=1}^N \|g^k\|^2_{\S_k^{-1}}
\end{array}
\end{equation}
for $g^k \in \partial f_k(x^k), k=1,\ldots, N$.
\end{proposition}
We obtain the following result from Proposition \ref{Oiter-complexitybbbb}.
\begin{theorem}\label{iter-complexityanother}
Let $\{(x^k,z^k,y^k):k=1,\ldots, N+1\}$ be generated by Online-spADMM. Suppose that
the following matrix orders hold:
\begin{equation}\label{eq:ordersbanother}
\S_1 \succeq \S_2 \succeq \cdots \succeq \S_{N+1}\succeq \S,
\end{equation}
where $\S$ is a positively definite self-adjoint operator.
Then the following properties hold:
\begin{itemize}
\item[(i)]  Let  Assumption  \ref{ass:matrixthtac}  be satisfied, then
\begin{equation}\label{eq:objective-compexityNanother}
\begin{array}{ll}
\displaystyle  \frac{1}{N} {\rm regret}^{\rm obj}_N &=\displaystyle \frac{1}{N}\displaystyle \sum_{k=1}^{N}\left\{[f_k(x^{k})+g(z^{k})]\right\}-\nu^*_N\\[12pt]
&\leq \displaystyle \frac{1}{2N}\left [ (\tau\sigma)^{-1}\|y^1\|^2+\|z^1-z^{0}\|_{\T}^2\right]
+\displaystyle \displaystyle \frac{1}{2N}\|(x^1,z^1)-(\widehat x,\widehat z)\|_{\widehat{\cal M}_1}^2\\[12pt]
& \quad + [\sigma\lambda_{max}(\S)]^{-1}L^2+\displaystyle \frac{1}{N}\Big(g(z^1)-g(z^{N+1})\Big).
\end{array}
\end{equation}

\item[(ii)]Let  Assumptions \ref{ass:value-bound} and \ref{ass:matrixthtac} be satisfied, then
\begin{equation}\label{eq:constraint-compexityanother}
\begin{array}{ll}
\displaystyle  \frac{1}{N} {\rm regret}^{\rm ctr}_N & =\displaystyle  \frac{1}{N}\sum_{k=1}^N\left\{\displaystyle   \|Ax^{k+1}+Bz^{k+1}-c\|^2\right\}\\[16pt]
& \leq \displaystyle \frac{1}{Nt_{\tau}\sigma}
\left [ (\tau\sigma)^{-1}\|y^1\|^2+\|z^1-z^{0}\|_{\T}^2\right]\\[10pt]
& \quad \,
+\displaystyle \displaystyle \frac{1}{Nt_{\tau}\sigma}\left[{\rm Dist}_{\widehat{\cal M}_1}((x^1,z^1),S^*_N)\right]^2 +\displaystyle \frac{2\gamma_0}{t_{\tau}\sigma}\\[10pt]
& \quad + \displaystyle \frac{2}{[t_{\tau}\sigma^2\lambda_{max}(\S)]}L^2+\displaystyle \frac{2}{t_{\tau}\sigma N}\Big(g(z^1)-g(z^{N+1})\Big).
\end{array}
\end{equation}
\item[(iii)]
Let   Assumption \ref{ass:matrixthtac} be satisfied and
\begin{equation}\label{ass:regretassanother}
\displaystyle \frac{1}{N}\sum_{k=1}^{N}\left\{[f_k(x^{k})+g(z^{k})]\right\}-\nu^*_N \geq 0,
\end{equation}
then the solution regret has the following bound
\begin{equation}\label{eq:constraint-compexityanother}
\begin{array}{l}
\displaystyle \frac{1}{N}\sum_{k=1}^N\|(x^{k+1},z^{k+1})-(\widehat x,\widehat z)\|_{{\cal W}_k}^2 \\[16pt]
\leq  \displaystyle \frac{1}{N}
\left [ (\tau\sigma)^{-1}\|y^1\|^2+\|z^1-z^{0}\|_{\T}^2\right]
+ \displaystyle \frac{1}{N}\left[{\rm Dist}_{\widehat{\cal M}_1}((x^1,z^1),S^*_N)\right]^2 \\[10pt]
\quad + \displaystyle \frac{2}{[\sigma\lambda_{max}(\S)]}L^2+\displaystyle \frac{2}{N}\Big(g(z^1)-g(z^{N+1})\Big).
\end{array}
\end{equation}
\end{itemize}
\end{theorem}
Define
$$
\begin{array}{l}
\mu_1(\tau)=\displaystyle \frac{1}{2}\left[ \|\S_1\| +\|BB^*\|+(1-\min(\tau,\tau^{-1}))\|E^*E\|\right],\\[10pt]
\mu_2(\tau)=\displaystyle \frac{1}{2}\left[ \tau^{-1}\|y^1\|^2+\|z^1-z^0\|_{\T}^2+2g(z^1)\right].
\end{array}
$$
We obtain the following result from Theorem \ref{iter-complexityanother}, which provides upper bounds in terms of $N$ for objective regret, constraint violation regret and solution regret.
\begin{corollary}\label{iter-complexityanotherCoro}
Let $\{(x^k,z^k,y^k):k=1,\ldots, N+1\}$ be generated by Online-spADMM. Suppose that
the following matrix orders hold:
\begin{equation}\label{eq:ordersbanothercoro}
\S_1 \succeq \S_2 \succeq \cdots \succeq \S_{N+1}\succeq \S,
\end{equation}
where $\S$ is a positively definite self-adjoint operator.
Then the following properties hold:
\begin{itemize}
\item[(i)]  Let  Assumption  \ref{ass:matrixthtac}  be satisfied, then
\begin{equation}\label{eq:objective-compexityNanothercoro}
\begin{array}{ll}
\displaystyle  \frac{1}{N} {\rm regret}^{\rm obj}_N &=\displaystyle \frac{1}{N}\displaystyle \sum_{k=1}^{N}\left\{[f_k(x^{k})+g(z^{k})]\right\}-\nu^*_N\\[16pt]
& \leq \displaystyle \left ( \displaystyle \frac{\|\T\|}{2}\displaystyle \frac{1}{N}+\mu_1(\tau)\displaystyle \frac{1}{\sqrt{N}}\right ) {\rm dist} ((x^1,z^1),S^*_N)^2\\[12pt]
& \quad + \Big(\mu_2(\tau)-g(z^{N+1}\Big) \displaystyle \frac{1}{N}+\displaystyle \frac{L^2}{\lambda_{max}(\S)}\displaystyle \frac{1}{\sqrt{N}}.
\end{array}
\end{equation}
\item[(ii)]Let  Assumptions \ref{ass:value-bound} and \ref{ass:matrixthtac} be satisfied, then
\begin{equation}\label{eq:constraint-compexityanothercoro}
\begin{array}{ll}
\displaystyle  \frac{1}{N} {\rm regret}^{\rm ctr}_N &=\displaystyle  \frac{1}{N}\sum_{k=1}^N\left\{\displaystyle   \|Ax^{k+1}+Bz^{k+1}-c\|^2\right\}\\[16pt]
& \leq \displaystyle \frac{2}{t_{\tau}}
\left [  \displaystyle \left ( \displaystyle \frac{\|\T\|}{2}\displaystyle \frac{1}{N^{3/2}}+\mu_1(\tau)\frac{1}{N}\displaystyle \right ) {\rm dist} ((x^1,z^1),S^*_N)^2\right.\\[12pt]
& \quad + \left.\Big(\mu_2(\tau)-g(z^{N+1}\Big) \displaystyle \frac{1}{N^{3/2}}+\displaystyle \frac{L^2}{\lambda_{max}(\S)}\displaystyle \frac{1}{N}
\right]
+\displaystyle \frac{2\gamma_0}{t_{\tau}}\displaystyle \frac{1}{\sqrt{N}}.\\[10pt]
\end{array}
\end{equation}
\item[(iii)]
Let   Assumption \ref{ass:matrixthtac} be satisfied and
\begin{equation}\label{ass:regretassanotherc}
\displaystyle \frac{1}{N}\sum_{k=}^{N}\left\{[f_k(x^{k})+g(z^{k})]\right\}-\nu^*_N \geq 0,
\end{equation}
then the solution regret has the following bound
\begin{equation}\label{eq:constraint-compexityanotherc}
\begin{array}{l}
\displaystyle \frac{1}{N}\sum_{k=1}^N\|(x^{k+1},z^{k+1})-(\widehat x,\widehat z)\|_{{\cal W}_k}^2 \\[16pt]
\leq \displaystyle \left ( \displaystyle \|\T\|\displaystyle \frac{1}{N}+2\mu_1(\tau)\displaystyle \frac{1}{\sqrt{N}}\right ) {\rm dist} ((x^1,z^1),S^*_N)^2\\[12pt]
\quad + 2\Big(\mu_2(\tau)-g(z^{N+1}\Big) \displaystyle \frac{1}{N}+\displaystyle \frac{2L^2}{\lambda_{max}(\S)}\displaystyle \frac{1}{\sqrt{N}}+\displaystyle \frac{2\gamma_0}{t_{\tau}}\displaystyle \frac{1}{\sqrt{N}}.
\end{array}
\end{equation}
\end{itemize}
\end{corollary}
We should point out that Assumption \ref{ass:matrixthta} is weaker than condition (\ref{eq:ordersbanothercoro}). However, if $\W_k$ satisfies
$$
\W_k \succeq \W \succ 0
$$
for some self-adjoint operator $\W$, then under  condition (\ref{ass:regretassanotherc}), (\ref{eq:constraint-compexityanotherc}) provides a solution regret. But whether  condition (\ref{ass:regretassanotherc})
always holds is a problem left.

\section{Averaging in spADMM}
\setcounter{equation}{0}

In this section, we consider the off-line problem, namely the case where $f_k(x)=f(x)$,$\forall k=1,\ldots, N$.
In this case, Problem (\ref{onlineP}) is reduced to
\begin{equation}\label{sconvex}
\begin{array}{ll}
\min & f(x)+g(z)\\
{\rm s.t.} & Ax+Bz=c,x \in {\cal X},z\in {\cal Z}.\\
\end{array}
\end{equation}
The augmented Lagrangian function of problem (\ref{sconvex})
is defined by
\[\label{augL}
{\cal L}_{\sigma}(x,z;y):=f (x)+g (z)+\langle y,Ax+Bz-c\rangle+\displaystyle \frac{\sigma}{2}\|Ax+Bz-c\|^2,\ \forall\, (x,z,y)\in \X\times \Z\times \Y.
\]
For this probem,   Online-spADMM is reduced to the semi-proximal alternating direction method of multipliers, namely spADMM proposed by \cite{Fazel}. \\
{\bf spADMM}: A   semi-proximal alternating direction method of multipliers for solving the convex optimization problem (\ref{sconvex}).
\mbox{}\vspace{-2mm}
\begin{description}
\item[Step 0 ] Input $(x^1,z^1,y^1)\in\hbox{dom}\; f \times \hbox{dom}\; g\times \mathcal{Y}$.
 Let  $\tau\in (0, +\infty)$ be a positive parameter   (e.g., $\tau \in ( 0, (1+\sqrt{5} )/2)$\,),  $\S_1:\X\to \X$ and   $\T:\Z\to \Z$ be a self-adjoint positive semidefinite, not necessarily positive definite, linear operators. Set $k:=1$.
\item[ Step 1]
Set \begin{equation}\label{xna}
\begin{array}{l}
x^{k+1}\in \hbox{arg}\min \, \L_{\sigma }(x,z^k; y^k) +\frac{1}{2}\|x-x^k\|^2_{\S_k}\, ,\\[2mm]
z^{k+1}\in \hbox{arg}\min\, \L_{\sigma  }(x^{k+1},z; y^k) +\frac{1}{2}\|z-z^k\|^2_\T\, ,\\[2mm]
y^{k+1}=y^k+\tau\sigma (Ax^{k+1}+Bz^{k+1}-c).
\end{array}
\end{equation}
\item[ Step 2] Set $k:=k+1$  and go to Step 1.
\end{description}

There is a slight difference between Online-spADMM and the above spADMM   proposed in \cite{Fazel}. In Online-spADMM, in stead of using $\S_k$ in Step 1, we use $\sigma \S_k$ in this step, because this is more convenient for regret bound analysis.

In this section, we first discuss regrets of spADMM iterations and the iteration complexity of the averaging of  generated iterative points. After that we will  recover an important inequality in \cite{HSZhang2018}, which is the key for establishing the linear rate of convergence for spADMM  under calmness condition of the inverse KKT mapping.
\subsection{Regrets of spADMM Iterations}
  Let $S^*$ and $\nu^*$  denote the solution set and the optimal value  of Problem (\ref{sconvex}), respectively.
We introduce the following notations:
\begin{equation}\label{eq:notationMH4}
\begin{array}{l}
\overline {\cal M}=\mbox{Diag}\left(\S+\Sigma_{f}, \T+\Sigma_g+\sigma B^*B\right)+s_{\tau}\overline E^*\overline E,\\[6pt]
\overline {\cal H}=\mbox{Diag}\left(\S+\displaystyle \frac{1}{2}\Sigma_{f}, \T+\Sigma_g+2t_{\tau}\tau \sigma B^*B\right)+\displaystyle \frac{1}{4}t_{\tau}\sigma \overline E^*\overline E,
\end{array}
\end{equation}
where $\overline E$ is a linear operator defined by $\overline E(x,z)=Ax+Bz$.

   In order to derive the constraint violation regret of spADMM for solving Problem (\ref{sconvex}), we need the following assumption similar to Assumption \ref{ass:value-bound}.
  \begin{assumption}\label{ass:value-bound4}
Suppose that the sequence $\{(x^k,z^k)\}$ generated by spADMM satisfies
$$
[f(\widehat x)+g(\widehat z)] -[f(x^{k})+g(z^{k})]\leq \gamma_0, \forall k=1,\ldots, N.
$$
for some $\gamma_0>0$ and $(\widehat x,\widehat z) \in S^*$.
\end{assumption}
 From Corollary \ref{cora}, we obtain the following result directly.
  \begin{proposition}\label{cora4-prop}
Let $\{(x^k,z^k,y^k)\}$ be generated by spADMM. Then, for any $(\widehat x,\widehat z) \in \Phi$, any $k=1,\ldots,$
\begin{equation}\label{eq:mainineqc4p}
\begin{array}{l}
[f(x^{k+1})+g(z^{k+1})]-[f(\widehat x)+g(\widehat z)]+\displaystyle \frac{1}{2}\sigma t_{\tau} \|Ax^{k+1}+Bz^{k+1}-c\|^2\\[10pt]
\leq \displaystyle \frac{1}{2}\left [ (\tau\sigma)^{-1}\|y^k\|^2+\|z^k-z^{k-1}\|_{\T}^2\right]- \displaystyle \frac{1}{2}\left [ (\tau\sigma)^{-1}\|y^{k+1}\|^2+\|z^{k+1}-z^{k}\|_{\T}^2\right]\\[10pt]
\quad \quad+\displaystyle \frac{1}{2}\left [ \|(x^k,z^k)-(\widehat x,\widehat z)\|_{\overline {\cal M}}^2-\|(x^{k+1},z^{k+1})-(\widehat x,\widehat z)\|_{\overline {\cal M}}^2\right]\\[10pt]
\quad \quad -\displaystyle \frac{1}{2}\left[\|(x^{k+1},z^{k+1})-(x^k,z^k)\|_{\overline {\cal H}}^2
\right].
\end{array}
\end{equation}
\end{proposition}
From (\ref{eq:mainineqc4p}), we obtain the following inequality:
\begin{equation}\label{eq:value-compexityp4}
 \begin{array}{l}
\displaystyle \frac{1}{N}\sum_{k=1}^{N}\left\{f(x^{k})+g(z^{k})\right\}-\nu^*+\displaystyle \frac{\sigma t_{\tau}}{2N} \sum_{k=1}^N\|Ax^k+Bz^k-c\|^2\\[12pt]
\leq \displaystyle \frac{1}{2N}\left [ (\tau\sigma)^{-1}\|y^1\|^2+\|z^1-z^{0}\|_{\T}^2\right]
+\displaystyle \displaystyle \frac{1}{2N}\left[{\rm Dist}_{\overline {\cal M}}((x^1,z^1),S^*)\right]^2\\[12pt]
+\displaystyle  \frac{1}{N}\left[f(x^1)+g(z^1)\right]+\displaystyle \frac{\sigma t_{\tau}}{2N} \|Ax^{1}+Bz^{1}-c\|^2.
\end{array}
\end{equation}
Define
where $$\begin{array}{l}
\eta_1=\displaystyle \frac{1}{t_{\tau}}\|B^*B\|,\quad \eta_2=\displaystyle \frac{1}{t_{\tau}} \left [ \tau^{-1}\|y^1\|^2+\|z^1-z^{0}\|_{\T}^2+2(f(x^1)+g(z^1))\right],\\[10pt]
\eta_3=\displaystyle \frac{1}{t_{\tau}}\left[\|\S+\Sigma_f\|+\|\T+\Sigma_g\|+s_{\tau}\|\overline E^*\overline E\|\right].
\end{array}$$
Then, using (\ref{eq:value-compexityp4}) and the definition of $\overline {\cal M}$, we obtain the following result about regrets of spADMM.
  \begin{theorem}\label{iter-complexityp4}
  Let $N$ be a  positive integer.
Let $\{(x^k,z^k,y^k)\}$ be generated by spADMM with $\sigma =\sqrt{N}$. Suppose that  Assumption \ref{ass:value-bound4}  is satisfied. Then  the following properties hold.
\begin{itemize}
\item[(i)] The objective regret satisfies the following bound:
 \begin{equation}\label{eq:value-compexityp4}
 \begin{array}{ll}
\displaystyle  \frac{1}{N} {\rm regret}^{\rm obj}_N &=\displaystyle \frac{1}{N}\sum_{k=1}^{N}\left\{f(x^{k})+g(z^{k})\right\}-\nu^*\\[12pt]
& \leq   \displaystyle \frac{t_{\tau}}{2}\displaystyle \frac{1}{\sqrt{N}} \left [ \|Ax^1+Bz^1-c\|^2+\eta_1\left[{\rm Dist}((x^1,z^1),S^*)\right]^2 \right]\\[16pt]
& \quad \quad\quad \quad\quad \, +  \displaystyle \frac{t_{\tau}}{2}\displaystyle\frac{1}{N}\left[\eta_2+\eta_3\left[{\rm Dist}((x^1,z^1),S^*)\right]^2 \right]
\end{array}
\end{equation}
\item[(ii)] The constraint violation regret has the following bound:
\begin{equation}\label{eq:constraint-compexityp4}
\begin{array}{ll}
\displaystyle  \frac{1}{N} {\rm regret}^{\rm ctr}_N &=\displaystyle  \frac{1}{N}\sum_{k=1}^N\left\{[\displaystyle   \|Ax^k+Bz^k-c\|^2\right\}\\[16pt]
& \leq \displaystyle  \frac{1}{N}  \left [ \|Ax^1+Bz^1-c\|^2+\eta_1\left[{\rm Dist}((x^1,z^1),S^*)\right]^2 \right]\\[16pt]
& \quad \quad  \quad  \quad  \quad \, +  \displaystyle\frac{1}{N^{3/2}}\left[\eta_2+\eta_3\left[{\rm Dist}((x^1,z^1),S^*)\right]^2 \right]
+\displaystyle \frac{2\gamma_0}{t_{\tau}}\displaystyle \frac{1}{\sqrt{N}}.
\end{array}
\end{equation}
\end{itemize}
\end{theorem}
Define for $t>1$,
\begin{equation}\label{eq:averageN}
(\widehat x^t, \widehat z^t,\widehat y^t)=\displaystyle\frac{1}{t}\sum_{j=1}^{t}(x^j,z^j,y^j).
\end{equation}
Then we can easily to obtain the following estimates from  Theorem \ref{iter-complexityp4}.

\begin{theorem}\label{iter-complexityp4a}
 Let $N$ be a  positive integer.
Let $\{(x^k,z^k,y^k)\}$ be generated by spADMM with $\sigma =\sqrt{N}$. Suppose that  Assumption \ref{ass:value-bound4}  is satisfied. Then  the following properties hold.
\begin{itemize}
\item[(i)] The error in  objective  at $(\widehat x^N,\widehat z^N)$ satisfies
 \begin{equation}\label{eq:value-compexityp4a}
 \begin{array}{l}
f(\widehat x^N )+g(\widehat z^N)-\nu^*\\[12pt]
\quad \quad \leq   \displaystyle \frac{t_{\tau}}{2}\displaystyle \frac{1}{\sqrt{N}} \left [ \|Ax^1+Bz^1-c\|^2+\eta_1\left[{\rm Dist}((x^1,z^1),S^*)\right]^2 \right]\\[16pt]
\quad \quad\quad \quad\quad \, +  \displaystyle \frac{t_{\tau}}{2}\displaystyle\frac{1}{N}\left[\eta_2+\eta_3\left[{\rm Dist}((x^1,z^1),S^*)\right]^2 \right]
\end{array}
\end{equation}
\item[(ii)] The error in  constraint at $(\widehat x^N,\widehat z^N)$ satisfies
\begin{equation}\label{eq:constraint-compexityp4a}
\begin{array}{l}
  \|A\widehat x^N+B\widehat z^N-c\|^2\leq \displaystyle  \frac{1}{N}  \left [ \|Ax^1+Bz^1-c\|^2+\eta_1\left[{\rm Dist}((x^1,z^1),S^*)\right]^2 \right]\\[16pt]
\quad \quad  \quad  \quad  \quad \, +  \displaystyle\frac{1}{N^{3/2}}\left[\eta_2+\eta_3\left[{\rm Dist}((x^1,z^1),S^*)\right]^2 \right]
+\displaystyle \frac{2\gamma_0}{t_{\tau}}\displaystyle \frac{1}{\sqrt{N}}.
\end{array}
\end{equation}
\end{itemize}
\end{theorem}

\subsection{Recovery of an important inequality in \cite{HSZhang2018}}

Let $(\overline x,\overline z)\in S^*$ be a solution to Problem  (\ref{sconvex}) and $\overline y \in {\cal Y}$  be a vector  such that $(\overline x,\overline z,\overline y)$ satisfies the following Karush-Kuhn-Tucker system
\begin{equation}\label{SCKKT}
 0\in \partial f(\overline x)+A^* \overline y, \,\, 0\in \partial g(\overline z)+B^* \overline y,\,\, c-A\overline x-B\overline z=0.
\end{equation}
From equalities
$$
\begin{array}{ll}
\left \langle y^{k+1},y^k-y^{k+1} \right \rangle &= \displaystyle \frac{1}{2}\left[ \|y^k\|^2-\|y^k-y^{k+1}\|^2-\|y^{k+1}\|^2\right]\\[6pt]
&= \displaystyle \frac{1}{2}\left[ \|y^k\|^2-\|y^{k+1}\|^2\right]-\displaystyle \frac{1}{2}[\sigma \tau]^2\|Ax^{k+1}+Bz^{k+1}-c\|^2,
\end{array}
$$
and
$$
\begin{array}{ll}
\left \langle y^{k+1},y^k-y^{k+1} \right \rangle &= \displaystyle \left\langle y^{k+1}-\overline y, y^k-\overline y-(y^{k+1}-\overline y) \right \rangle+\displaystyle \left\langle \overline y,y^k-y^{k+1}\right \rangle\\[6pt]
&=\displaystyle \frac{1}{2}\left[ \|y^k-\overline y\|^2-\|y^k-y^{k+1}\|^2-\|y^{k+1}-\overline y\|^2\right]+\displaystyle \left\langle \overline y,y^k-y^{k+1}\right \rangle\\[6pt]
&= \displaystyle \frac{1}{2}\left[ \|y^k-\overline y\|^2-\|y^{k+1}-\overline y\|^2\right]-\displaystyle \frac{1}{2}[\sigma \tau]^2\|Ax^{k+1}+Bz^{k+1}-c\|^2\\[6pt]
& \quad \quad -(\tau\sigma)\displaystyle \left\langle \overline y,Ax^{k+1}+Bz^{k+1}-c\right \rangle,
\end{array}
$$
 we obtain
 \begin{equation}\label{eq:yestimates}
 \|y^k\|^2-\|y^{k+1}\|^2=\|y^k-\overline y\|^2-\|y^{k+1}-\overline y\|^2-2(\tau\sigma)\displaystyle \left\langle \overline y,Ax^{k+1}+Bz^{k+1}-c\right \rangle.
 \end{equation}
Noting from (\ref{SCKKT}) that $(\overline x,\overline z)\in \Phi$, we have from Corollary \ref{cora}, for
 $\{(x^k,z^k,y^k)\}$  generated by spADMM, that
\begin{equation}\label{eq:mainineqc}
\begin{array}{l}
[f(x^{k+1})+g(z^{k+1})]-[f(\widehat x)+g(\widehat z)]+\displaystyle \frac{1}{2}\sigma t_{\tau} \|Ax^{k+1}+Bz^{k+1}-c\|^2\\[10pt]
\leq \displaystyle \frac{1}{2}\left [ \|(x^k,z^k)-(\widehat x,\widehat z)\|_{\overline{\cal M}}^2+\displaystyle  \frac{1}{\tau\sigma} \|y^k-\overline y\|^2+\|z^k-z^{k-1}\|_{\T}^2\right]\\[10pt]
\quad \quad- \displaystyle \frac{1}{2}\left [ \|(x^{k+1},z^{k+1})-(\widehat x,\widehat z)\|_{\overline {\cal M}}^2+\displaystyle  \frac{1}{\tau\sigma}\|y^{k+1}-\overline y\|^2+\|z^{k+1}-z^{k}\|_{\T}^2\right]\\[10pt]
\quad \quad -\displaystyle \frac{1}{2}\left[\|(x^{k+1},z^{k+1})-(x^k,z^k)\|_{\overline {\cal H}}^2
\right]\\[10pt]
\quad \quad -\displaystyle \frac{1}{2}\sigma t_{\tau} \|Ax^{k+1}+Bz^{k+1}-c\|^2\\[10pt]
\end{array}
\end{equation}
In terms of (\ref{SCKKT}), we know that $-A^* \overline y\in \partial f(\overline x)$ and $-B^* \overline y\in \partial g(\overline z)$. Then we obtain form the convexity of $f$ and $g$ that
$$
[f(x^{k+1})+g(z^{k+1})]-[f(\widehat x)+g(\widehat z)]+\displaystyle \left\langle \overline y,Ax^{k+1}+Bz^{k+1}-c\right \rangle\geq 0.
$$
From this and (\ref{eq:mainineqc}), we get
\begin{equation}\label{eq:mainineqd}
\begin{array}{l}
0
\leq \displaystyle \frac{1}{2}\left [ \|(x^k,z^k)-(\widehat x,\widehat z)\|_{\overline{\cal M}}^2+\displaystyle  \frac{1}{\tau\sigma} \|y^k-\overline y\|^2+\|z^k-z^{k-1}\|_{\T}^2\right]\\[10pt]
\quad \quad- \displaystyle \frac{1}{2}\left [ \|(x^{k+1},z^{k+1})-(\widehat x,\widehat z)\|_{\overline {\cal M}}^2+\displaystyle  \frac{1}{\tau\sigma}\|y^{k+1}-\overline y\|^2+\|z^{k+1}-z^{k}\|_{\T}^2\right]\\[10pt]
\quad \quad -\displaystyle \frac{1}{2}\left[\|(x^{k+1},z^{k+1})-(x^k,z^k)\|_{\overline {\cal H}}^2
\right]\\[10pt]
\quad \quad -\displaystyle \frac{1}{2}\sigma t_{\tau} \|Ax^{k+1}+Bz^{k+1}-c\|^2\\[10pt]
\end{array}
\end{equation}
or
\begin{equation}\label{eq:mainineqd}
\begin{array}{l}
0
\leq \displaystyle \left [ \|(x^k,z^k)-(\widehat x,\widehat z)\|_{\overline{\cal M}}^2+\displaystyle  \frac{1}{\tau\sigma} \|y^k-\overline y\|^2+\|z^k-z^{k-1}\|_{\T}^2\right]\\[10pt]
\quad \quad- \displaystyle \left [ \|(x^{k+1},z^{k+1})-(\widehat x,\widehat z)\|_{\overline {\cal M}}^2+\displaystyle  \frac{1}{\tau\sigma}\|y^{k+1}-\overline y\|^2+\|z^{k+1}-z^{k}\|_{\T}^2\right]\\[10pt]
\quad \quad -\displaystyle \left[\|(x^{k+1},z^{k+1})-(x^k,z^k)\|_{\overline {\cal H}}^2
\right]\\[10pt]
\quad \quad -\displaystyle (\sigma\tau^2)^{-1} t_{\tau} \|y^{k+1}-y^k\|^2\\[10pt]
\end{array}
\end{equation}
We introduce the following linear operators, the same notations as in \cite{HSZhang2018},
$$
{\cal M}={\rm Diag}\Big(\overline{\cal M},\, (\sigma \tau)^{-1}{\cal I}\Big),\quad  {\cal H}={\rm Diag}\Big(\overline{\cal H}, \, (\sigma\tau^2)^{-1} t_{\tau} {\cal I}\Big),
$$
where ${\cal I}$ is the identity operator in ${\cal Y}$.
Then inequality (\ref{eq:mainineqd}) is equivalent to
\begin{equation}\label{eq:recovery}
\begin{array}{l}
\displaystyle \left [ \|(x^{k+1},z^{k+1},y^{k+1})-(\widehat x,\widehat z,\widehat y)\|_{ {\cal M}}^2+\|z^{k+1}-z^{k}\|_{\T}^2\right]\\[12pt]
\quad \quad \leq \displaystyle \left [ \|(x^{k},z^{k},y^{k})-(\widehat x,\widehat z,\widehat y)\|_{ {\cal M}}^2+\|z^{k}-z^{k-1}\|_{\T}^2\right]\\[10pt]
\quad \quad \quad -\displaystyle \left[\|(x^{k+1},z^{k+1},y^{k+1})-(x^k,z^k,y^k)\|_{{\cal H}}^2
\right],
\end{array}
\end{equation}
which coincides with the important formula (26) in \cite{HSZhang2018}. Thus inequality  (\ref{eq:mainineqc})  in Corollary \ref{cora} is not so strong as  formula (26) in \cite{HSZhang2018}.

\section{Examples and Numerical Evaluations}

In this section, numerical results of the proposed algorithms are presented. We apply the Online-spADMM to several specific questions and conduct numerical experiments to validate the theoretical performance of our algorithm acting on synthetic data. We get the time-averaged objective regret and time-average constraint violations of our algorithm. In this experiment, we intend to compare the performance of Online-spADMM with the following well-studied algorithms. We will test the following numerical cases provided by \cite{boyd11}.

In this part, we evaluate the performance of Online-spADMM for solving online quadratic optimization, Lasso and total variation (TV), respectively. All computational results are obtained by running Matlab R2020a on Windows 10 (Intel Core i5-10400 CPU @ 2.90GHz 16GB RAM).

\subsection{Application to online quadratic optimization}
Consider the online quadratic optimization (OQO) problem with

\begin{equation}\label{eq:qd}
    f_t(x)=\displaystyle \frac{1}{2}x^TG_tx+c_t^Tx, \, x \in \Theta,
\end{equation}
where
\begin{equation}\label{eq:Ac}
\Theta=\{x\in X: Ax=b\},
\end{equation}
and $X \subset \Re^n$ is a closed convex compact set, $G_{t} \in \Re^{n \times n}$, $c_t \in \Re^n$, $A \in \Re^{m\times n}$ and $b \in \Re^m$.

We will test the following numerical case provided in Boyd's website (http://www.stanford.edu/
$\thicksim$boyd/papers/admm/quadprog/quadprog\_example.html). For better numerical result, we generate a well-conditioned symmetrical positive definite matrix $G_{t} \in \mathbb S^{n \times n}$ and vector $c_t \in \Re^n$ at time $t$, where $G_{t}$ is uniformly distributed over $(0, 1)$ and $c_t$ generated by the standard normal distribution.

Define
$$
g(z)=\delta_X(z)
$$
and
$$
\Phi=\{(x,z)\in \Re^n \times \Re^m:Ax-b=0,z-x=0\}.
$$
We can reformulate the OQO problem into the framework of Problem (\ref{onlineP}), i.e.,

\begin{equation} \label{OQO}
	\left\{
		\begin{array}{rl}
			\textup{min} & \displaystyle\frac{1}{2}x^TG_tx+c_t^Tx + \delta_X(z) \\[4pt]
			\textup{s.t.} & Ax - b = 0 \\
                            &~~ x - z = 0.
		\end{array}
	\right.
\end{equation}


Then the online quadratic optimization problem is reformulated as Problem (\ref{OQO}). The augmented Lagrangian is defined as
$$
{\cal L}^t_{\sigma}(x,z;\mu,\lambda)=\displaystyle \frac{1}{2}x^TG_tx+c_t^Tx+\delta_X(z)+
\langle \mu, Ax-b\rangle +\langle \lambda, x-z \rangle+\displaystyle \frac{\sigma}{2} \|Ax-b\|^2+\displaystyle \frac{\sigma}{2} \|x-z\|^2.
$$
Then Online-spADMM may be described as follows.\\
{\bf Online-spADMM} for  online convex quadratic optimization problem (\ref{OQO}).

\begin{description}
\item[Step 0 ] Input $(x^1,z^1,\mu^1,\lambda^1)\in \Re^n \times \Re^n \times \Re^m\times \Re^n .$ Let  $\tau\in (0, +\infty)$ be a positive parameter   (e.g., $\tau \in ( 0, (1+\sqrt{5} )/2)$\,),  $S_1\in \mathbb S^n_+$ be a symmetric positive semidefinite matrix.  Set $k:=1$.
\item[ Step 1]
Set \begin{equation}\label{xnaq}
\begin{array}{l}
x^{k+1}\in \hbox{arg}\min \, \frac{1}{2}x^TG_tx+c_t^Tx+
\langle A^{T}\mu + \lambda, x \rangle + \frac{\sigma}{2} \|Ax-b\|^2+\frac{\sigma}{2}\|x-z\|^2 +\frac{\sigma}{2}\|x-x^k\|^2_{S_k}\, ,\\[2mm]

z^{k+1}\in \hbox{arg}\min\, \delta_X(z) - \langle \lambda^{k}, z \rangle + \frac{\sigma}{2}\|x^{k+1}-z\|^2 ,\\[2mm]

\mu^{k+1}=\mu^k+\tau\sigma (Ax^{k+1}-b),\\[2mm]

\lambda^{k+1}=\lambda^k+\tau\sigma (x^{k+1}-z^{k+1}).
\end{array}
\end{equation}
\item[ Step 2] Receive a cost function $f_{k+1}$ and incur loss $f_{k+1}(x^{k+1})$ and constraint violation
$\|Ax^{k+1}-b\|$.
\item[ Step 3] Choose a symmetric positive semidefinite matrix $S_{k+1}\in \mathbb S^n_+$.
\item[ Step 4] Set $k:=k+1$  and go to Step 1.
\end{description}
For integer $N$, choose $\alpha >0$ large enough such that
$$
\alpha I-\displaystyle \frac{1}{\sigma}G_t-A^TA\succeq 0, \forall t=1,\ldots, N.
$$
Define
\begin{equation}\label{k:Sk}
S_t=\alpha I-\displaystyle \frac{1}{\sigma}G_t-A^TA\succeq 0, \forall t=1,\ldots, N.
\end{equation}
 Let $\sigma=\sqrt{N}$. Then subproblems for $x^{k+1}$ and $z^{k+1}$ in (\ref{xnaq}) have the following explicit solutions:
\begin{equation}\label{eq:Solq}
\begin{array}{l}
x^{k+1}=\displaystyle \frac{1}{(\alpha+1)\sqrt{N}}[-c_k-A^T\mu^k-\lambda^k]+\displaystyle \frac{1}{(\alpha+1)}[A^Tb+z^k+S_kx^k],\\[10pt]
z^{k+1}=\Pi_X\left(x^{k+1}+\displaystyle \frac{\lambda^k}{\sqrt{N}}\right).
\end{array}
\end{equation}
Therefore, when $X$ is a simple convex set and $\Pi_X$ is easy to calculate, then Online-spADMM with $S_t$ defined by (\ref{k:Sk}) is quite effective as subproblems in (\ref{xnaq}) have explicit expressions.

In our numerical tests, we choose $n = [10, 20, 50, 100]$, total number of iterations $T = 5000$, $\tau = [0.1, 0.3, 1.618]$ and $\sigma = \sqrt{N}$. In the existing literature, dimensions of numerical examples are generally very
low, usually 2 or 5, for instances see \cite{chaudhary2021safe} , \cite{hosseini2014online} and \cite{liu2018online}.


In this experiment, we intend to compare the performance
of online-mspADMM with the Online Alternating Direction Method (OADM) algorithm proposed in \cite{WangBan2013} is in the following form:
\begin{equation*}
	\left\{
		\begin{array}{rl}
			x^{k + 1} = &(G_k + \eta_1 A^{T}A + (\eta_1 + \eta_2)I_n)^{-1}[\eta_1(z^k + A^{T}b) + \eta_2 x^k - (A^T \mu_k + \lambda_k + c_k)], \\[6pt]
			z^{k + 1} = & \textup{max}\{0, x^{k+1}+\lambda^k / \eta_1\}, \\[6pt]
			\mu^{k+1}=&\mu^k+\eta_1 (Ax^{k+1}-b),\\[2mm]
            \lambda^{k+1}=&\lambda^k+\eta_1 (x^{k+1}-z^{k+1}).
		\end{array}
	\right.
\end{equation*}
We set the parameters to $\eta_1 = \sqrt{N}$ and $\eta_2 = T$.

\begin{figure}[htbp]
	\centering
	\subfigure[time-averaged objective regret (n = 10)]{
	\includegraphics[height =4cm, width=5.5cm]{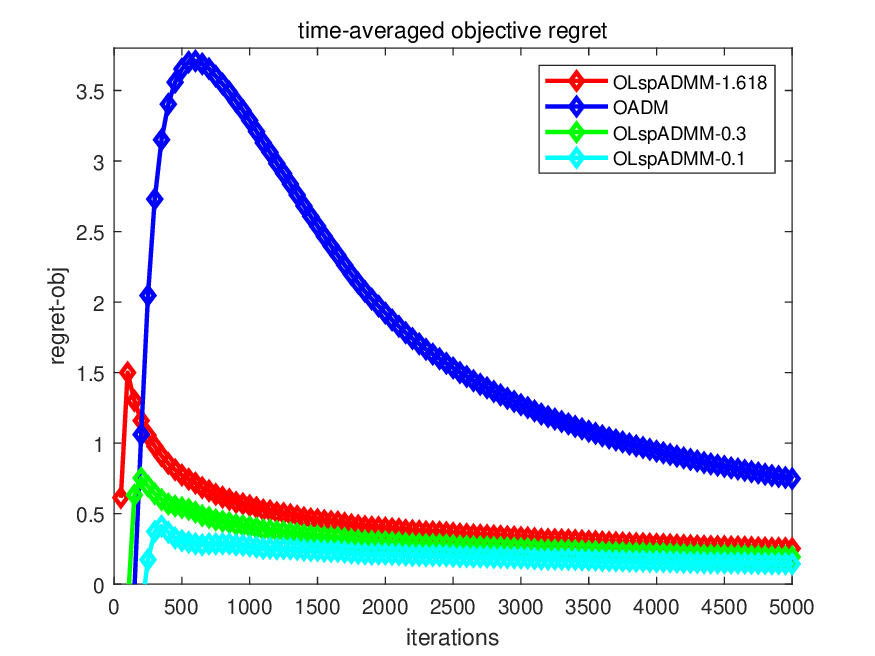}\label{figure1a}
	}
	\quad
	\subfigure[time-average constraint violation (n = 10)]{
	\includegraphics[height =4cm, width=5.5cm]{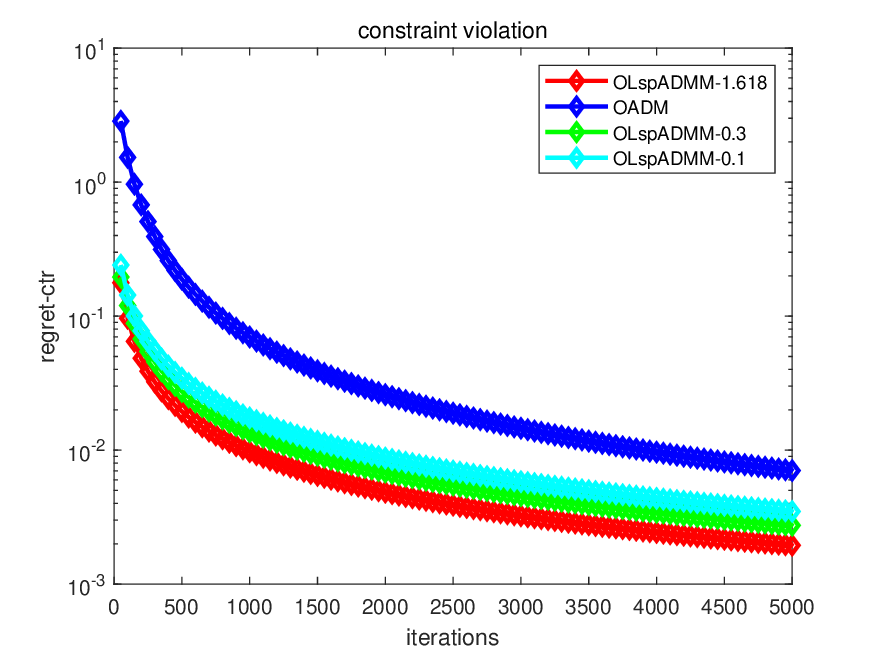} \label{figure1b}
	}
    \subfigure[time-averaged objective regret (n = 20)]{
	\includegraphics[height =4cm, width=5.5cm]{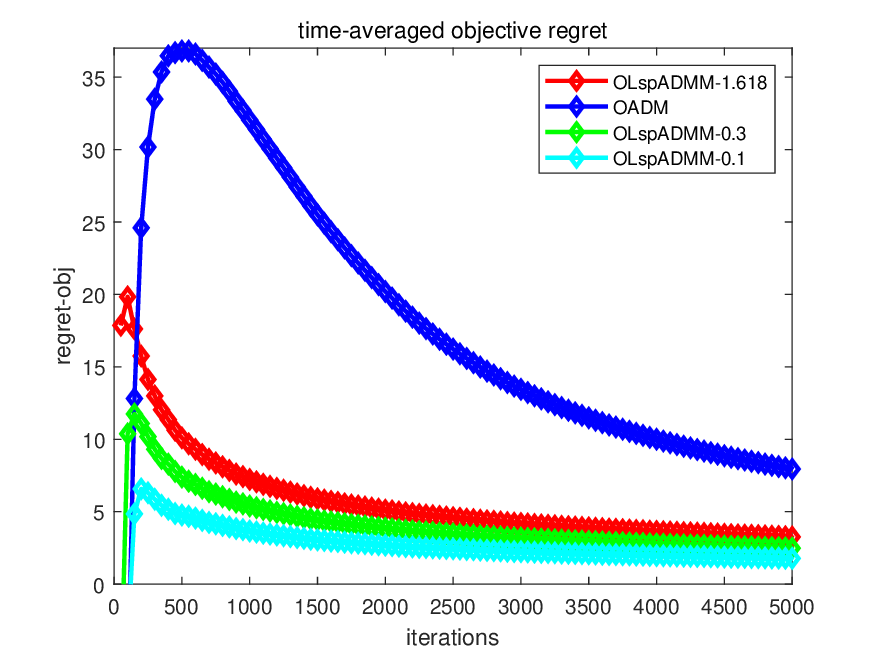} \label{figure1b}
	}
    \quad
    \subfigure[time-average constraint violation (n = 20)]{
	\includegraphics[height =4cm, width=5.5cm]{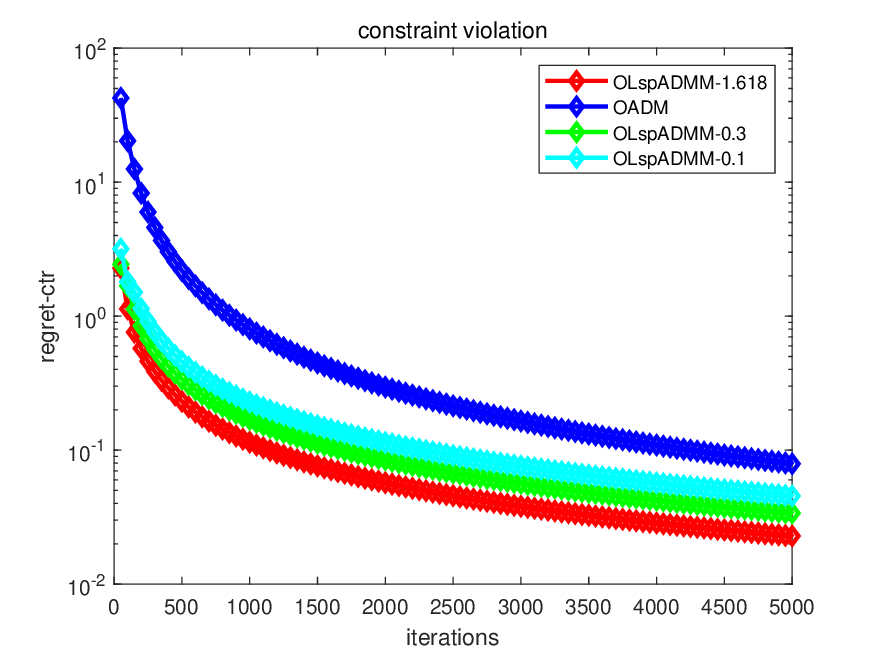} \label{figure1b}
	}
\subfigure[time-averaged objective regret (n = 50)]{
	\includegraphics[height =4cm, width=5.5cm]{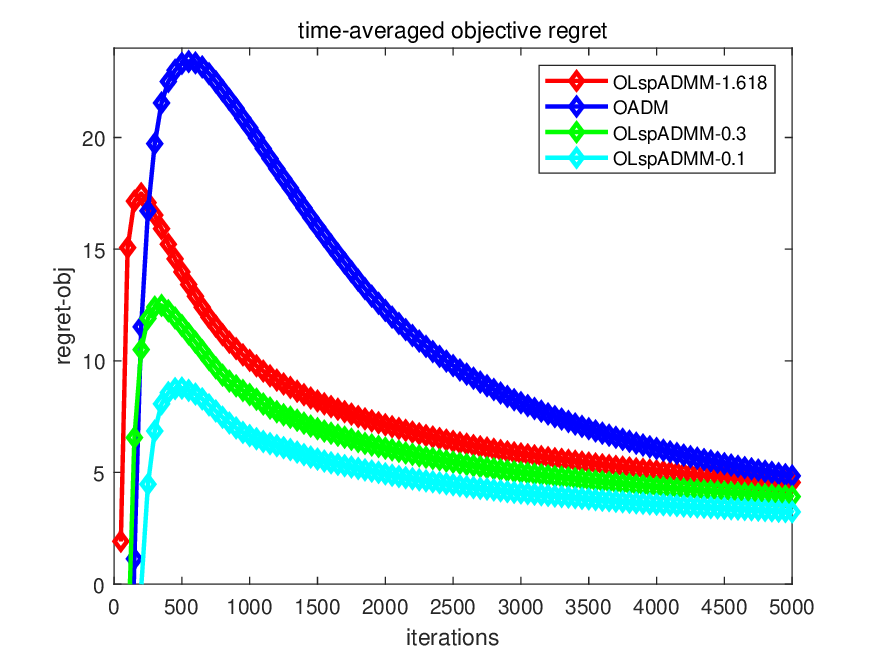} \label{figure1b}
	}
    \quad
    \subfigure[time-average constraint violation (n = 50)]{
	\includegraphics[height =4cm, width=5.5cm]{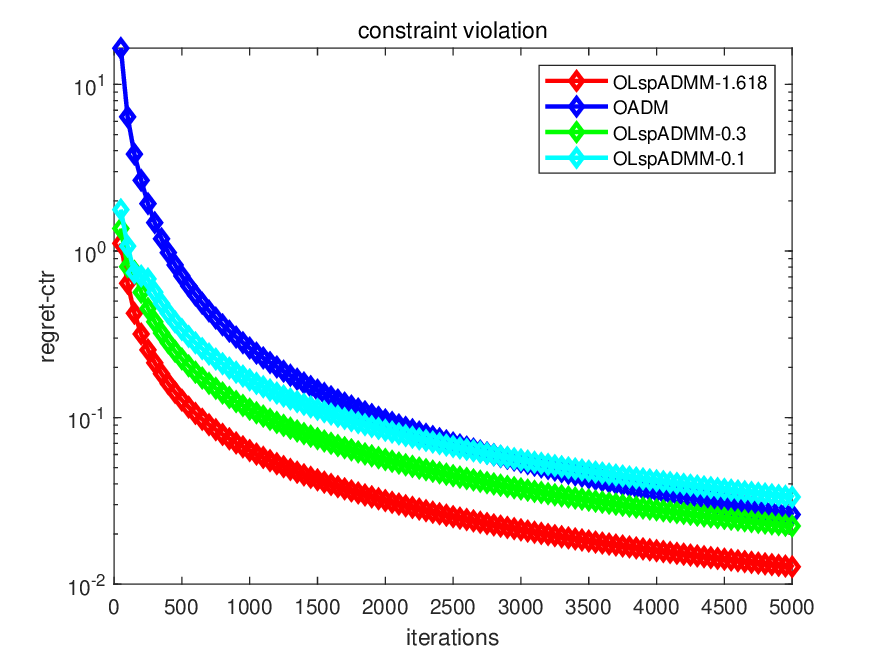} \label{figure1b}
	}
\subfigure[time-averaged objective regret (n = 100)]{
	\includegraphics[height =4cm, width=5.5cm]{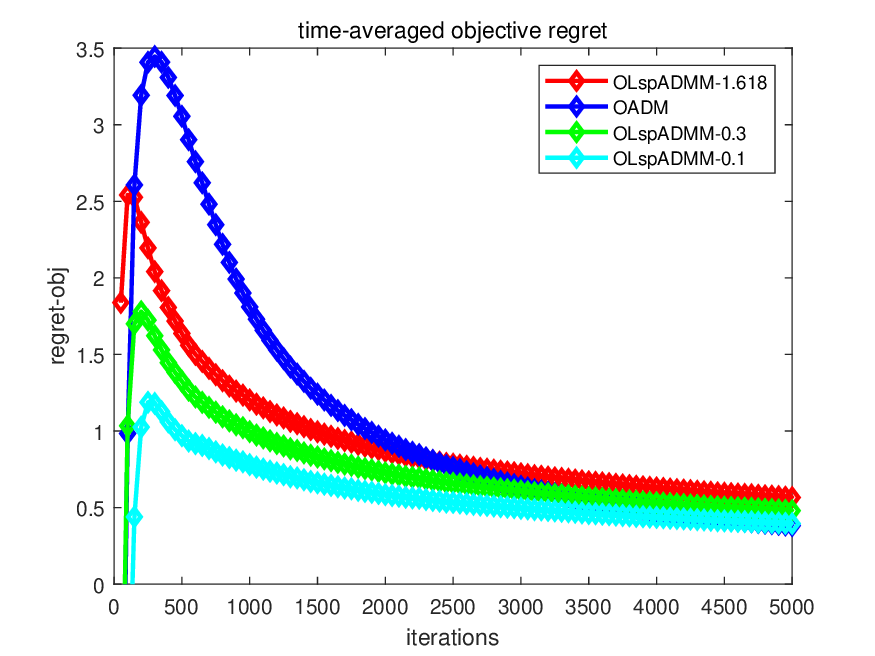} \label{figure1b}
	}
    \quad
    \subfigure[time-average constraint violation (n = 100)]{
	\includegraphics[height =4cm, width=5.5cm]{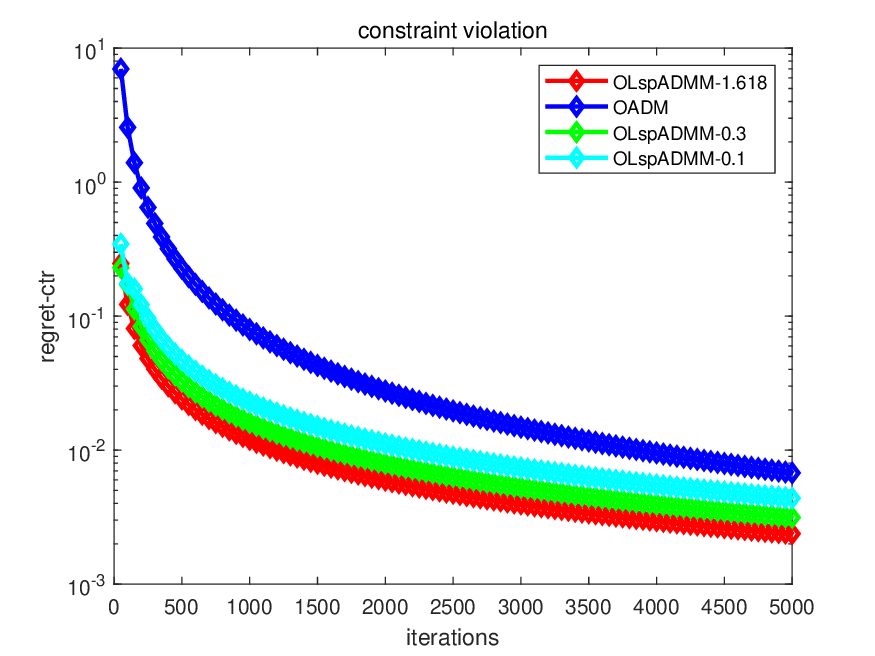} \label{figure1b}
	}
	\caption{Comparison of algorithms with respect to the time-average objective regret and time-average constraint violation for OQO.} \label{figure1}
\end{figure}

The results are shown in Figure \ref{figure1}. The images show the numerical results of the OADM and Online-spADMM algorithms with the different selections of parameter $\tau$. We can see that Online-spADMM achieves $\mathcal{O}(\sqrt{N})$ objective regret and $\mathcal{O}(\sqrt{N})$ constraint violation regret at round $N$, which validates our theoretical results. It is obvious that the Online-spADMM algorithm performs slightly better than the OADM algorithm. Combine Figure \ref{figure1} and Table \ref{table1}, we find that the smaller $\tau$, the smaller time-averaged objective regret, the time-average constraint violation is the opposite. Table \ref{table2} reports the comparison between the running time of the Online-spADMM and OADM for all cases. From this table, we can observe that Online-spADMM has a shorter running time.

\begin{table}[h]
\caption{Comparison between the time-averaged objective regret of Online semi-proximal ADMM with the different selections of parameter $\tau$ (OLspADMM-$\tau$) and Online Alternating Direction Method(OADM), where the total number of iterations $T = 5000$.} \label{table1}
\begin{tabular*}{\textwidth}{@{\extracolsep\fill}lcccc}
\toprule%
\multicolumn{1}{@{}c@{}}{Dimension} &
\multicolumn{1}{@{}c@{}}{OLspADMM-1.618} & \multicolumn{1}{@{}c@{}}{OLspADMM-0.3} & \multicolumn{1}{@{}c@{}}{OLspADMM-0.1} &
\multicolumn{1}{@{}c@{}}{OADM}\\
\midrule
10  & 0.253 & 0.193 & 0.143 & 0.748 \\
20  & 3.258 & 2.483 & 1.777 & 7.946 \\
50  & 4.555 & 3.918 & 3.232 & 4.847 \\
100 & 0.566 & 0.480 & 0.395 & 0.450 \\
\bottomrule
\end{tabular*}
\end{table}

\begin{table}[h]
\caption{Comparison between the running time (sec) of Online semi-proximal ADMM (Online-spADMM) and Online Alternating Direction Method (OADM).} \label{table2}
\begin{tabular*}{\textwidth}{@{\extracolsep\fill}lcc}
\toprule%
&
 \multicolumn{1}{@{}c@{}}{Online-spADMM} & \multicolumn{1}{@{}c@{}}{OADM} \\\cmidrule{2-3}%
IterNum & n = 10$|$n = 20$|$n = 50$|$n = 100 & n = 10$|$n = 20$|$n = 50$|$n = 100 \\
\midrule
5000  & 0.023 $|$ 0.028 $|$ 0.097 $|$ 0.264 & 0.029 $|$ 0.087 $|$ 0.222 $|$ 0.628 \\
10000 & 0.046 $|$ 0.055 $|$ 0.212 $|$ 0.536 & 0.055 $|$ 0.174 $|$ 0.441 $|$ 1.261 \\
20000 & 0.090 $|$ 0.108 $|$ 0.420 $|$ 1.095 & 0.220 $|$ 0.358 $|$ 0.879 $|$ 2.757 \\
50000 & 0.225 $|$ 0.272 $|$ 1.086 $|$ 2.933 & 0.562 $|$ 0.782 $|$ 2.249 $|$ 6.515 \\
\bottomrule
\end{tabular*}
\end{table}

\subsection{Lasso}

In this subsection, we study a numerical example of Lasso. Lasso \cite{tibshirani1996regression} is an important special case of $l_{1}$ regularized linear regression. For a given scalar $\lambda \geq 0$, the Lasso regularizer is chosen as $\varphi(x) = \lambda \|x\|_{1}$. This involves solving

\begin{equation} \label{Lasso}
	\mathop{\textup{min}}\limits_{x} ~~ \displaystyle\frac{1}{2}\|A_{t}x - b_{t}\|^{2} + \lambda \|x\|_{1},
\end{equation}
where $A_{t} \in \Re^{m\times n}$ and $b_{t} \in \Re^m$. By introducing an auxiliary variable $z \in \Re^{n}$, we can reformulate the Lasso problem \eqref{Lasso} into the framework of Problem \eqref{onlineP}, i.e.,

\begin{equation} \label{OCO-Lasso}
	\left\{
		\begin{array}{rl}
			\textup{min} & \displaystyle\frac{1}{2}\|A_{t}x - b_{t}\|^{2} + \lambda \|z\|_{1} \\[4pt]
			\textup{s.t.} & x - z = 0,
		\end{array}
	\right.
\end{equation}
where $A_{t}$ and $b_t$ are generated by the standard normal distribution. The specific generating code can be found in Boyd's website (http://www.stanford.edu/$\thicksim$boyd/papers/admm/lasso/lasso$\_$example.html).


Then the Lasso is reformulated as Problem (\ref{OCO-Lasso}). The augmented Lagrangian is defined as
$$
{\cal L}^t_{\sigma}(x,z;y)=\displaystyle \frac{1}{2}\|A_{t}x - b_{t}\|^{2} + \lambda \|z\|_{1} +
\langle y, x - z \rangle + \displaystyle \frac{\sigma}{2} \|x - z\|^2.
$$
Then Online-spADMM may be described as follows.\\
{\bf Online-spADMM} for Lasso (\ref{Lasso}).

\begin{description}
\item[Step 0 ] Input $(x^1,z^1,\mu^1,\lambda^1)\in \Re^n \times \Re^n \times \Re^m\times \Re^n .$ Let  $\tau\in (0, +\infty)$ be a positive parameter   (e.g., $\tau \in ( 0, (1+\sqrt{5} )/2)$\,),  $S_1\in \mathbb S^n_+$ be a symmetric positive semidefinite matrix.  Set $k:=1$.
\item[ Step 1]
Set \begin{equation}\label{xnaq_lasso}
\begin{array}{l}
x^{k+1}\in \hbox{arg}\min \, \frac{1}{2}\|A_k x - b_k\|^2 +
\langle y^k, x \rangle + \frac{\sigma}{2} \|x - z^k\|^2 + \frac{\sigma}{2}\|x-x^k\|^2_{S_k}\, ,\\[2mm]

z^{k+1}\in \hbox{arg}\min\, \lambda \|z\|_1 + \frac{\sigma}{2}\|z-(x^{k+1} + \frac{y^k}{\sigma})\|^2 ,\\[2mm]

y^{k+1}=y^k+\tau\sigma (x^{k+1}-z^{k+1}).
\end{array}
\end{equation}
\item[ Step 2] Receive a cost function $f_{k+1}$ and incur loss $f_{k+1}(x^{k+1})$ and constraint violation
$\|x^{k+1} - z^{k+1}\|$.
\item[ Step 3] Choose a symmetric positive semidefinite matrix $S_{k+1}\in \mathbb S^n_+$.
\item[ Step 4] Set $k:=k+1$  and go to Step 1.
\end{description}
For integer $N$, choose $\alpha >0$ large enough such that
$$
\alpha I-\displaystyle \frac{1}{\sigma} A^{T}_{k}A_{k} \succeq 0, \forall t=1,\ldots, N.
$$
Define
\begin{equation}\label{k:Sk}
S_t=\alpha I-\displaystyle \frac{1}{\sigma} A^{T}_{k}A_{k} \succeq 0, \forall t=1,\ldots, N.
\end{equation}
 Let $\sigma=\sqrt{N}$. Then subproblems for $x^{k+1}$ and $z^{k+1}$ in (\ref{xnaq_lasso}) have the following explicit solutions:
\begin{equation}\label{eq:Solq}
\begin{array}{rl}
x^{k+1}= &\displaystyle \frac{1}{(\alpha+1)}(z^k+S_kx^k)-\displaystyle \frac{1}{(\alpha+1)\sqrt{N}}(y^{k}-A^{T}_{k}b_{k}),\\[10pt]
z^{k+1}= &\hbox{arg}\min \displaystyle\frac{\lambda}{\sqrt{N}}\|z\|_{1} + \frac{1}{2}\|z-(x^{k+1}+\displaystyle\frac{y^k}{\sqrt{N}})\|^2\\[10pt]
= &\textup{sgn}(x^{k + 1} + y^{k} / \sqrt{N}) \cdot \textup{max}(|x^{k + 1} + y^{k} / \sqrt{N}| - \lambda / \sqrt{N}, 0).
\end{array}
\end{equation}

In our numerical tests, we choose $n = [10, 20, 50]$, total number of iterations $T = 5000$, $\tau = [0.1, 0.3, 1.618]$ and $\sigma = \sqrt{N}$. We intend to compare the performance
of Online-spADMM with the following well-studied algorithms:
\begin{itemize}
\item[$\bullet$] OADM. The Online Alternating Direction Method (OADM) algorithm
 proposed in \cite{WangBan2013} is in the following form:

 \begin{equation*}
	\left\{
		\begin{array}{rl}
			x^{k + 1} = & ((\eta_1 + \eta_2)I + A^T_k A_k)^{-1} [(\eta_2 x_k + \eta_1 z^k + A^T_k b_k) - y^k ], \\[6pt]
			z^{k + 1} = & \textup{argmin} ~ \lambda \|z\|_{1} + \frac{\eta_1}{2} \|z - (x^{k + 1} + \frac{y^{k}}{\eta_1})\|^{2} \\[4pt]
			= & \textup{sgn}(x^{k + 1} + y^{k} / \eta_1) \cdot \textup{max}(|x^{k + 1} + y^{k} / \eta_1| - \lambda / \eta_1, 0), \\[6pt]
			y^{k + 1} = & y^{k} + \eta_1 (x^{k + 1} - z^{k + 1}).
		\end{array}
	\right.
\end{equation*}
We set the parameters to $\eta_1 = \sqrt{N}$ and $\eta_2 = T / 2$.
\item[$\bullet$] FOBOS. The Online Forward-Backward Splitting Method (FOBOS) algorithm
 proposed in \cite{JohnSinger2009} is in the following form:

 \begin{equation*}
	\left\{
		\begin{array}{rl}
			w^{k + 1} = & x_k - \rho_k A^T_k (A_k x_k - b_k), \\[6pt]
			x^{k + 1} = & \textup{argmin} ~ \lambda \rho_{k+1} \|x\|_{1} + \frac{1}{2} \|x - w^{k + 1}\|^{2} \\[4pt]
			= & \textup{sgn}(w_{k+1}) \cdot \textup{max}(|w_{k+1}| - \lambda \rho_{k+1}, 0),
		\end{array}
	\right.
\end{equation*}
where parameters $\rho_k \propto 1 / k$.
\item[$\bullet$] RDA. The Regularized Dual Averaging Method (RDA) algorithm
 proposed in \cite{Xiao2010} is in the following form:

 \begin{equation*}
	\left\{
		\begin{array}{rl}
			g^{k} = & A^T_k (A_kx_k - b_k), \\[6pt]
			\bar{g}^{k} = & \frac{k-1}{k} \bar{g}^{k} + \frac{1}{k}g_k, \\[4pt]
			x_{k+1} = & \textup{argmin} ~ (\lambda + \frac{\eta \beta_k}{k}) \|x\|_{1} + \frac{\beta_k}{2k} \|x + \frac{k}{\beta_k}\bar{g}_k\|^{2} \\[4pt]
			= & \textup{sgn}(-\frac{k}{\beta_k}\bar{g}_k) \cdot \textup{max}(|-\frac{k}{\beta_k}\bar{g}_k| - (\frac{\lambda k}{\beta_k} + \eta), 0).
		\end{array}
	\right.
\end{equation*}
We set the parameters to $\eta = 0.005$ and $\beta_k = \gamma \sqrt{k}$, where $\gamma = 5000$.
\end{itemize}

\begin{figure}[htbp]
	\centering
	\subfigure[time-averaged objective regret (n = 10)]{
	\includegraphics[height =4cm, width=5.5cm]{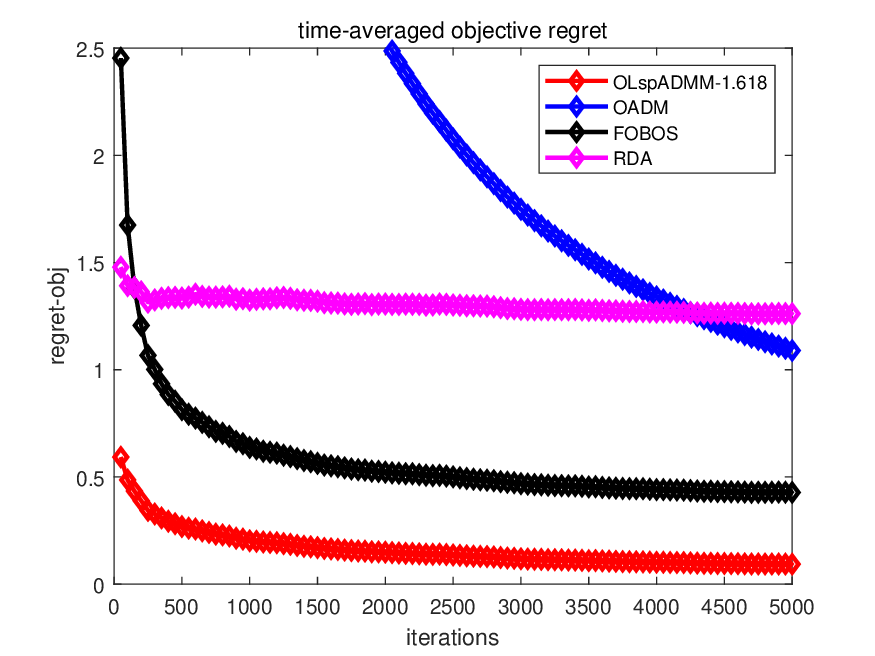}\label{figure1a}
	}
	\quad
	\subfigure[time-averaged constraint violation (n = 10)]{
	\includegraphics[height =4cm, width=5.5cm]{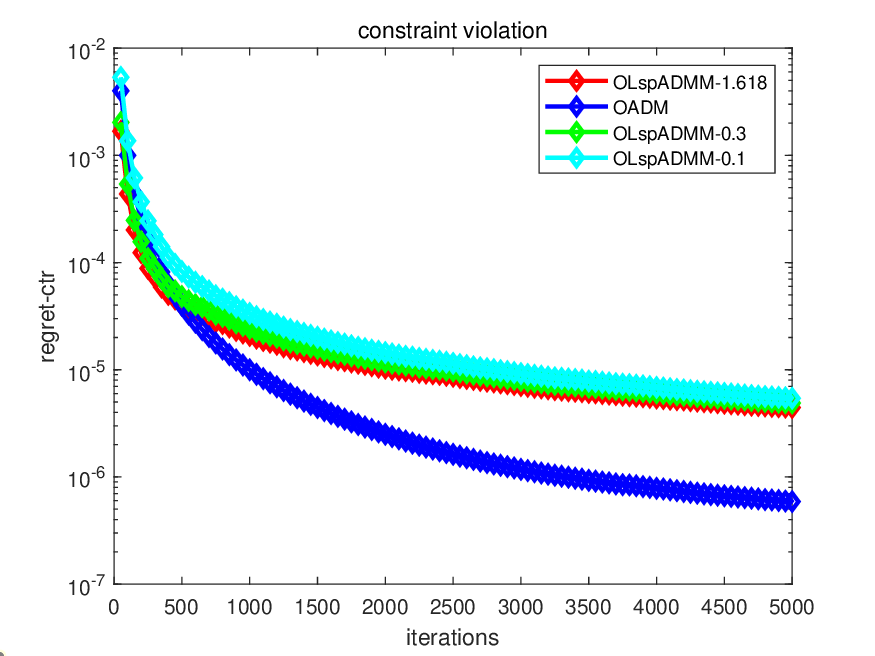} \label{figure1b}
	}
    \subfigure[time-averaged objective regret (n = 20)]{
	\includegraphics[height =4cm, width=5.5cm]{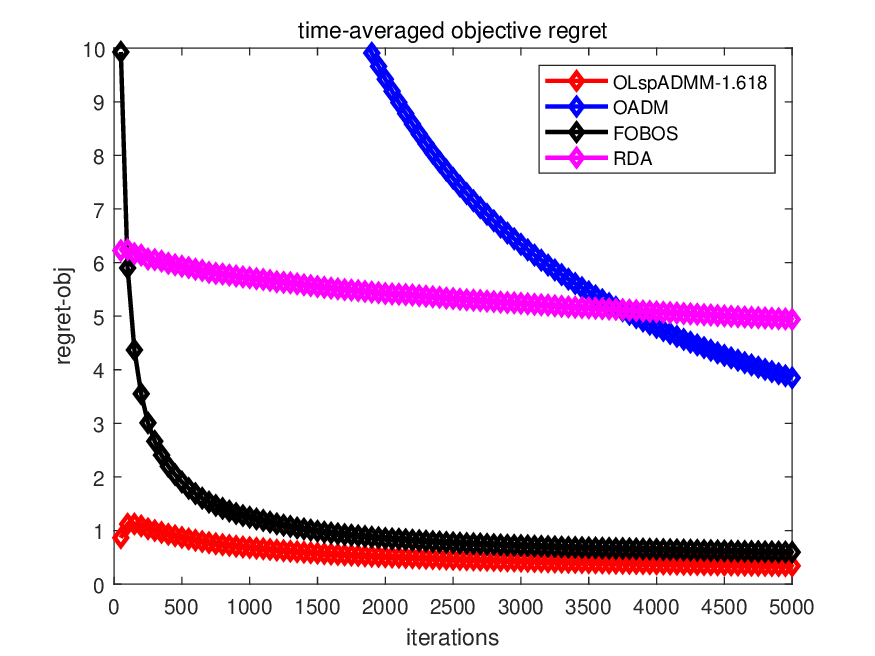} \label{figure1b}
	}
    \quad
    \subfigure[time-averaged constraint violation (n = 20)]{
	\includegraphics[height =4cm, width=5.5cm]{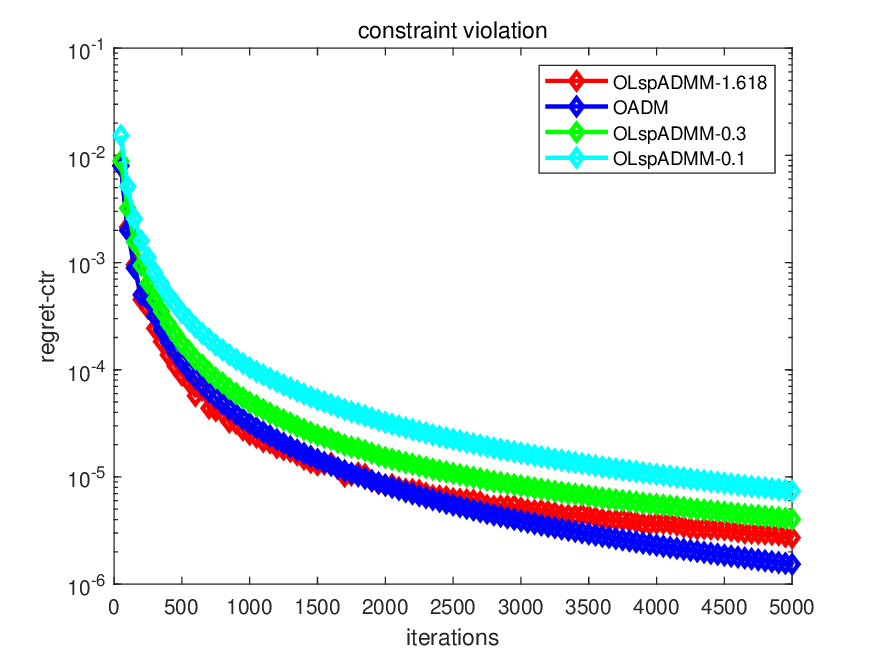} \label{figure1b}
	}
\subfigure[time-averaged objective regret (n = 50)]{
	\includegraphics[height =4cm, width=5.5cm]{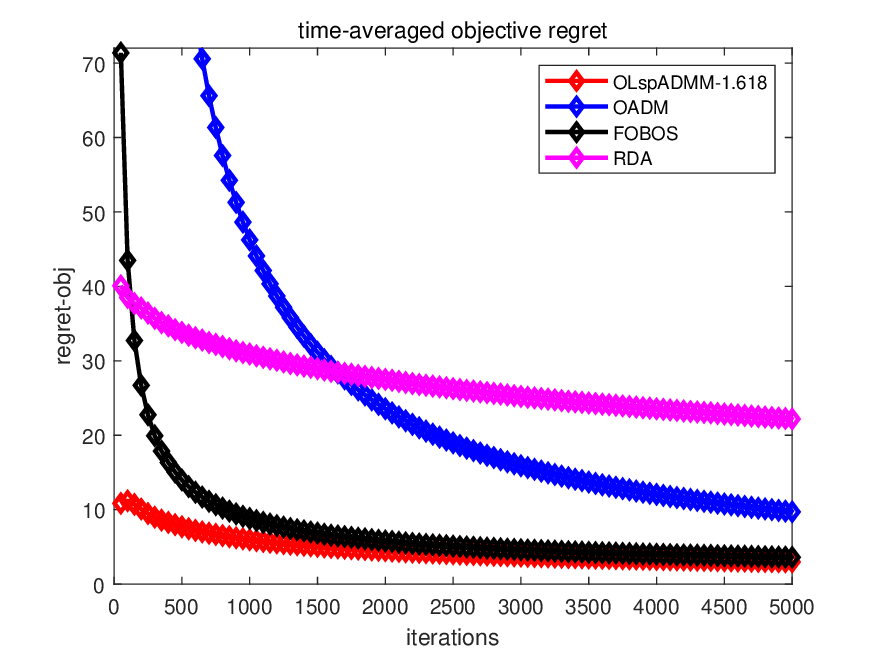} \label{figure1b}
	}
    \quad
    \subfigure[time-averaged constraint violation (n = 50)]{
	\includegraphics[height =4cm, width=5.5cm]{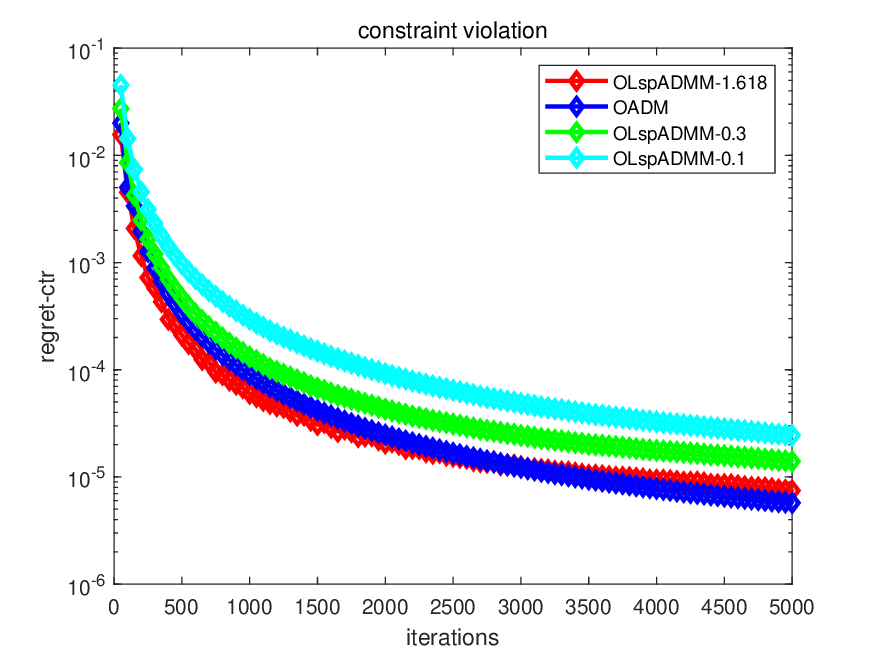} \label{figure1b}
	}
	\caption{Comparison of algorithms with respect to the time-average objective regret and constraint violation for Lasso.} \label{figure2}
\end{figure}

The results are shown in Figure \ref{figure2}. Since parameter $\tau$ changes have little effect on time-averaged objective regret, we don't consider the effect of $\tau$. The images show the numerical results of the OADM and Online-spADMM algorithms with the different selections of parameter $\tau$. Combine Figure \ref{figure2} and Table \ref{table3}, we can see that Lasso and OQO have the similar numerical performance about the time-averaged objective regret. From Table \ref{table3} with different dimensions, we observe that Online-spADMM always performs better than the other three algorithms. Since FOBOS and RDA are used to solve unconstrained online optimization problems, they don't have time-average constraint violation as a performance metric. We just have to care about the Online-spADMM and OADM. It's easy to observe that the smaller $\tau$, the larger constraint violation and the performance of Online-spADMM not much worse than OADM. From table \ref{table4}, we can observe that Online-spADMM has a shorter running time.

\begin{table}[h]
\caption{Comparison between the time-averaged objective regret of Online semi-proximal ADMM (Online-spADMM), Online Alternating Direction Method(OADM), Online Forward-Backward Splitting Method (FOBOS) and Regularized Dual Averaging Method (RDA), where the total number of iterations $T = 5000$.} \label{table3}
\begin{tabular*}{\textwidth}{@{\extracolsep\fill}lcccc}
\toprule%
\multicolumn{1}{@{}c@{}}{Dimension} &
\multicolumn{1}{@{}c@{}}{Online-spADMM} & \multicolumn{1}{@{}c@{}}{OADM} &
\multicolumn{1}{@{}c@{}}{FOBOS} &
\multicolumn{1}{@{}c@{}}{RDA}\\
\midrule
10  & 0.472 & 2.086 & 1.081 & 1.262 \\
20  & 0.341 & 3.850 & 0.596 & 4.939 \\
50  & 2.981 & 9.708 & 3.606 & 22.180 \\
\bottomrule
\end{tabular*}
\end{table}

\begin{table}[h]
\caption{Same as Table \ref{table2} but for Lasso.} \label{table4}
\begin{tabular*}{\textwidth}{@{\extracolsep\fill}lcc}
\toprule%
&
 \multicolumn{1}{@{}c@{}}{Online-spADMM} & \multicolumn{1}{@{}c@{}}{OADM} \\\cmidrule{2-3}%
IterNum & n = 10$|$n = 20$|$n = 50 & n = 10$|$n = 20$|$n = 50 \\
\midrule
5000  & 0.026 $|$ 0.031 $|$ 0.056 & 0.059 $|$ 0.072 $|$ 0.155 \\
10000 & 0.052 $|$ 0.061 $|$ 0.112 & 0.113 $|$ 0.141 $|$ 0.310 \\
20000 & 0.102 $|$ 0.125 $|$ 0.223 & 0.225 $|$ 0.282 $|$ 0.619 \\
50000 & 0.257 $|$ 0.308 $|$ 0.567 & 0.556 $|$ 0.705 $|$ 1.551 \\
\bottomrule
\end{tabular*}
\end{table}

\subsection{Generalized Lasso}

Online-spADMM is more powerful for problems with complex objective function. Such as Problem (\ref{OQO}), FOBOS and RDA are no longer applicable, since there will be no closed-form for it. The other important lasso problem can be generalized to
\begin{equation}
	\mathop{\textup{min}}\limits_{x} ~~ \displaystyle\frac{1}{2}\|A_{t}x - b_{t}\|^{2} + \lambda \|Fx\|_{1},
\end{equation}
where $F$ is an arbitrary linear transformation. An important special case is when $F \in \Re^{(n-1) \times n}$ is the difference matrix,
\begin{equation*}
	Fx = [x_{1} - x_{2}, x_{2} - x_{3}, \dots, x_{n - 1} - x_{n}]^{T}, \forall x \in \Re^{n}.
\end{equation*}
When $A_t = I$, generalized lasso can be expressed as
\begin{equation} \label{TV}
	\mathop{\textup{min}}\limits_{x} ~~ \displaystyle\frac{1}{2}\|x - b_{t}\|^{2} + \lambda \sum_{i = 1}^{n-1} |x_{i+1} - x_i|.
\end{equation}
This problem is the TV model of the removing noise from images problem \cite{rudin1992nonlinear}.

By introducing an auxiliary variable $z \in \Re^{n}$, we can reformulate the Lasso problem \eqref{TV} into the framework of Problem \eqref{onlineP}, i.e.,

\begin{equation} \label{OCO-TV}
	\left\{
		\begin{array}{rl}
			\textup{min} & \frac{1}{2}\|x - b_{t}\|^{2} + \lambda \|z\|_{1} \\[4pt]
			\textup{s.t.} & Fx - z = 0,
		\end{array}
	\right.
\end{equation}
where $b_t$ is generated by the standard normal distribution. The specific generating code can be found in Boyd's website (http://www.stanford.edu/$\thicksim$boyd/papers/admm/total$\_$variation/total$\_$variation$\_$ex-\\ample.html), although we modified the codes to accommodate our setup.


Then the Lasso is reformulated as Problem (\ref{OCO-TV}). The augmented Lagrangian is defined as
\begin{equation*}
	\begin{aligned}
		\mathcal{L}_{\sigma}^{t} (x, z; y) = & \frac{1}{2}\|x - b_{t}\|^{2} + \lambda \|z\|_{1} + \langle y, Fx - z \rangle + \frac{\sigma}{2} \|Fx - z\|^{2}.
	\end{aligned}
\end{equation*}
Then Online-spADMM may be described as follows.\\
{\bf Online-spADMM} for TV problem (\ref{TV}).

\begin{description}
\item[Step 0 ] Input $(x^1,z^1,\mu^1,\lambda^1)\in \Re^n \times \Re^n \times \Re^m\times \Re^n .$ Let  $\tau\in (0, +\infty)$ be a positive parameter   (e.g., $\tau \in ( 0, (1+\sqrt{5} )/2)$\,),  $S_1\in \mathbb S^n_+$ be a symmetric positive semidefinite matrix.  Set $k:=1$.
\item[ Step 1]
Set \begin{equation}\label{xnaq_TV}
\begin{array}{l}
x^{k+1}\in \hbox{arg}\min \, \frac{1}{2}\|x - b_k\|^2 +
\langle F^T y^k, x \rangle + \frac{\sigma}{2} \|Fx - z^k\|^2 + \frac{\sigma}{2}\|x-x^k\|^2_{S_k}\, ,\\[2mm]

z^{k+1}\in \hbox{arg}\min\, \lambda \|z\|_1 + \frac{\sigma}{2}\|z - (Fx^{k+1}+\frac{y^k}{\sigma})\|^2 ,\\[2mm]

y^{k+1}=y^k+\tau\sigma (Fx^{k+1}-z^{k+1}).
\end{array}
\end{equation}
\item[ Step 2] Receive a cost function $f_{k+1}$ and incur loss $f_{k+1}(x^{k+1})$ and constraint violation
$\|Fx^{k+1} - z^{k+1}\|$.
\item[ Step 3] Choose a symmetric positive semidefinite matrix $S_{k+1}\in \mathbb S^n_+$.
\item[ Step 4] Set $k:=k+1$  and go to Step 1.
\end{description}
For integer $N$, choose $\alpha >0$ large enough such that
$$
(\alpha-\displaystyle\frac{1}{\sigma}) I - F^{T}F \succeq 0, \forall t=1,\ldots, N.
$$
Define
\begin{equation}\label{k:Sk}
S_t=(\alpha-\displaystyle\frac{1}{\sigma}) I - F^{T}F \succeq 0, \forall t=1,\ldots, N.
\end{equation}
 Let $\sigma=a\sqrt{N}$, where $a$ is a scalar. Then subproblems for $x^{k+1}$ and $z^{k+1}$ in (\ref{xnaq_TV}) have the following explicit solutions:
\begin{equation}\label{eq:Solq}
\begin{array}{rl}
x^{k+1}= &\displaystyle \frac{1}{\alpha}(S_kx^k + F^T z^k)-\displaystyle \frac{1}{\alpha a \sqrt{N}}(F^T y^k - b_k),\\[10pt]
z^{k+1}= &\hbox{arg}\min \displaystyle\frac{\lambda}{\sqrt{N}}\|z\|_{1} + \frac{1}{2}\|z-(Fx^{k+1}+\displaystyle\frac{y^k}{\sqrt{N}})\|^2\\[10pt]
= &\textup{sgn}(Fx^{k + 1} + y^{k} / \sqrt{N}) \cdot \textup{max}(|Fx^{k + 1} + y^{k} / \sqrt{N}| - \lambda / \sqrt{N}, 0).
\end{array}
\end{equation}

In our numerical tests, we choose $n = [10, 20, 50, 100]$,  total number of iterations $T = 5000$, $\tau = [0.1, 0.3, 1.618]$ and $\sigma = \sqrt{N}$. We intend to compare the performance
of Online-spADMM with the Online Alternating Direction Method (OADM) algorithm proposed in \cite{WangBan2013} is in the following form:

\begin{equation*}
	\left\{
		\begin{array}{rl}
			x^{k + 1} = & ((1 + \eta_2)I + \eta_1 F^T F)^{-1} [(\eta_2 x_k + \eta_1 F^T z^k + b_k) - F^T y^k ], \\[6pt]
			z^{k + 1} = & \textup{argmin} ~ \lambda \|z\|_{1} + \frac{\eta_1}{2} \|z - (Fx^{k + 1} + \frac{y^{k}}{\eta_1})\|^{2} \\[4pt]
			= & \textup{sgn}(Fx^{k + 1} + y^{k} / \eta_1) \cdot \textup{max}(|Fx^{k + 1} + y^{k} / \eta_1| - \lambda / \eta_1, 0), \\[6pt]
			y^{k + 1} = & y^{k} + \eta_1 (x^{k + 1} - z^{k + 1}).
		\end{array}
	\right.
\end{equation*}
We set the parameters to $\eta_1 = \sqrt{N}$ and $\eta_2 = T / 2$.

\begin{figure}[htbp]
	\centering
	\subfigure[time-averaged objective regret (n = 10)]{
	\includegraphics[height =4cm, width=5.5cm]{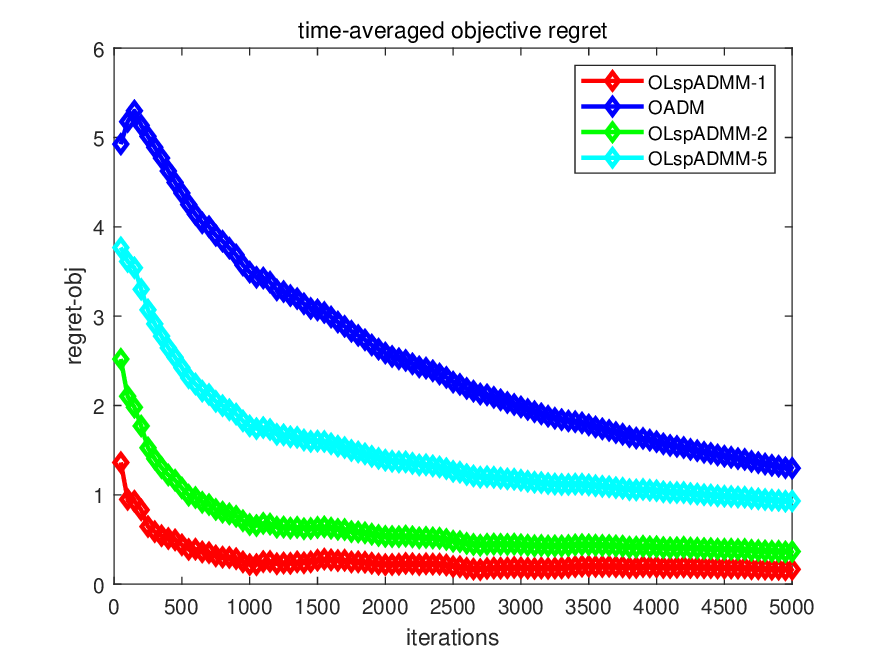}\label{figure1a}
	}
	\quad
	\subfigure[time-averaged constraint violation (n = 10)]{
	\includegraphics[height =4cm, width=5.5cm]{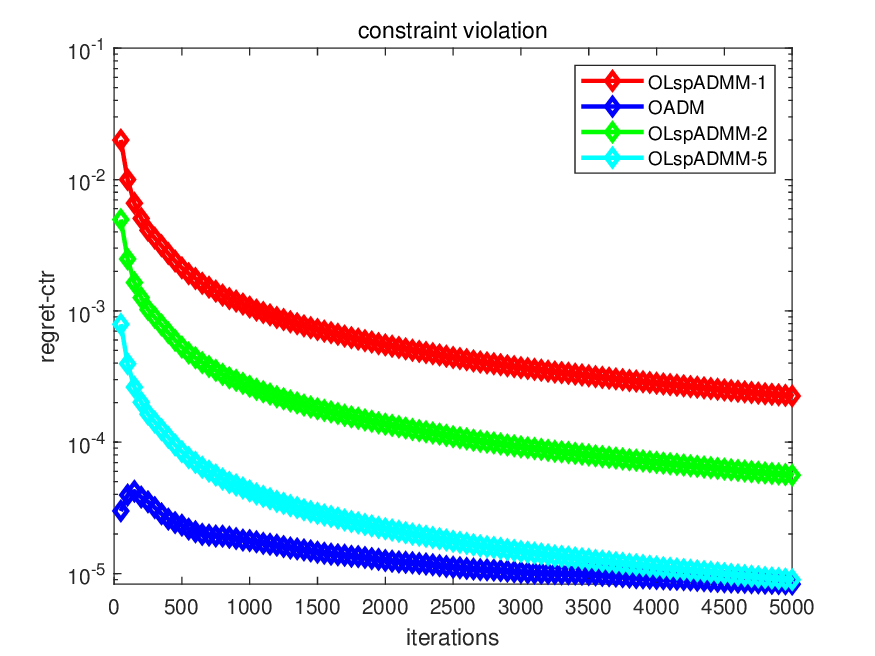} \label{figure1b}
	}
    \subfigure[time-averaged objective regret (n = 20)]{
	\includegraphics[height =4cm, width=5.5cm]{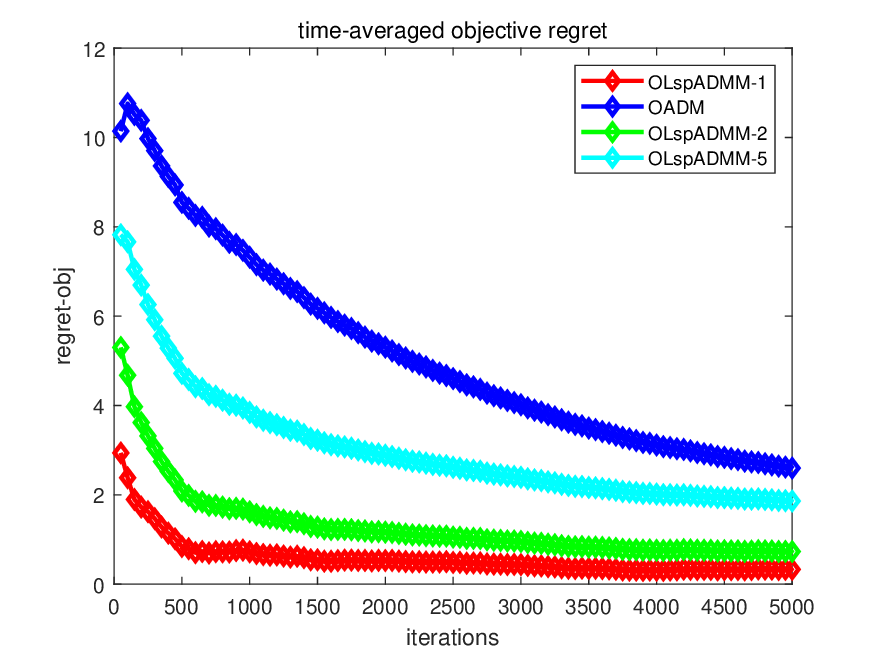} \label{figure1b}
	}
    \quad
    \subfigure[time-averaged constraint violation (n = 20)]{
	\includegraphics[height =4cm, width=5.5cm]{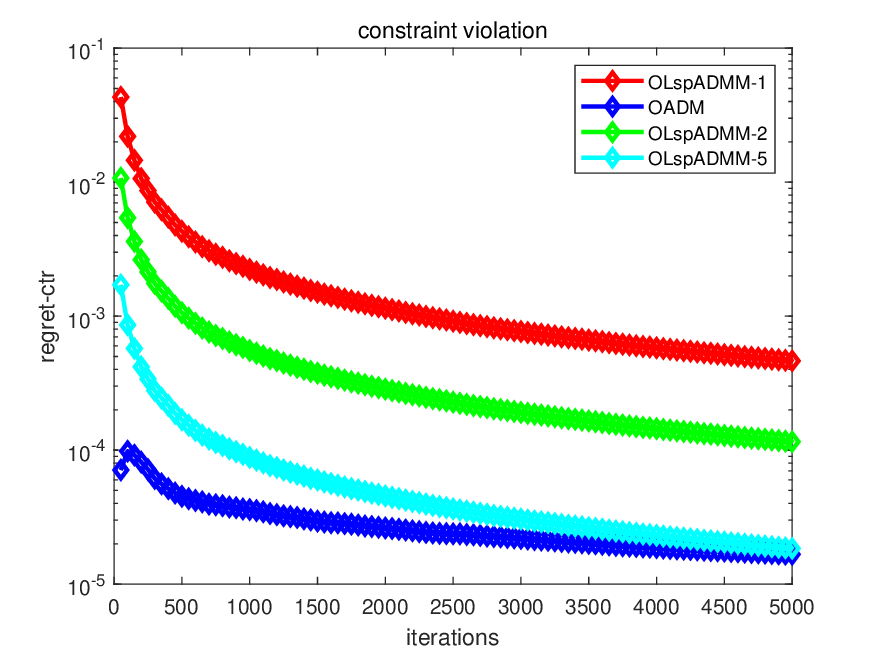} \label{figure1b}
	}
    \subfigure[time-averaged objective regret (n = 50)]{
	\includegraphics[height =4cm, width=5.5cm]{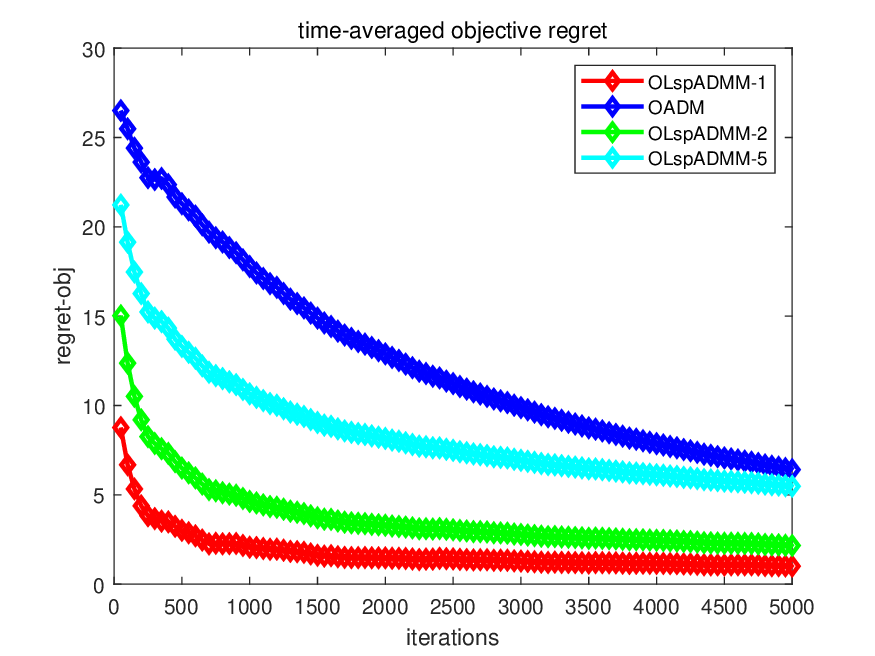} \label{figure1b}
	}
    \quad
    \subfigure[time-averaged constraint violation (n = 50)]{
	\includegraphics[height =4cm, width=5.5cm]{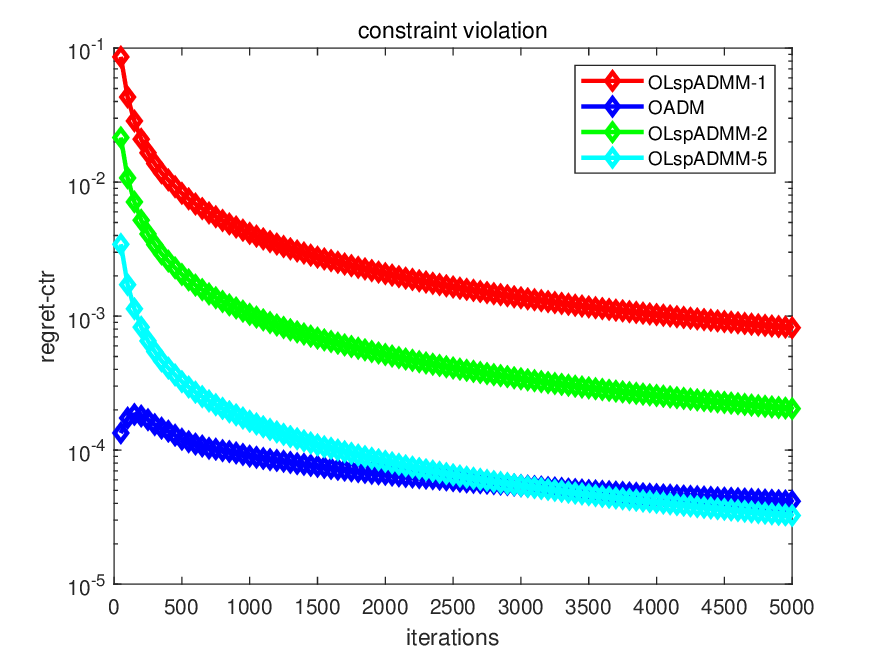} \label{figure1b}
	}
    \subfigure[time-averaged objective regret (n = 100)]{
	\includegraphics[height =4cm, width=5.5cm]{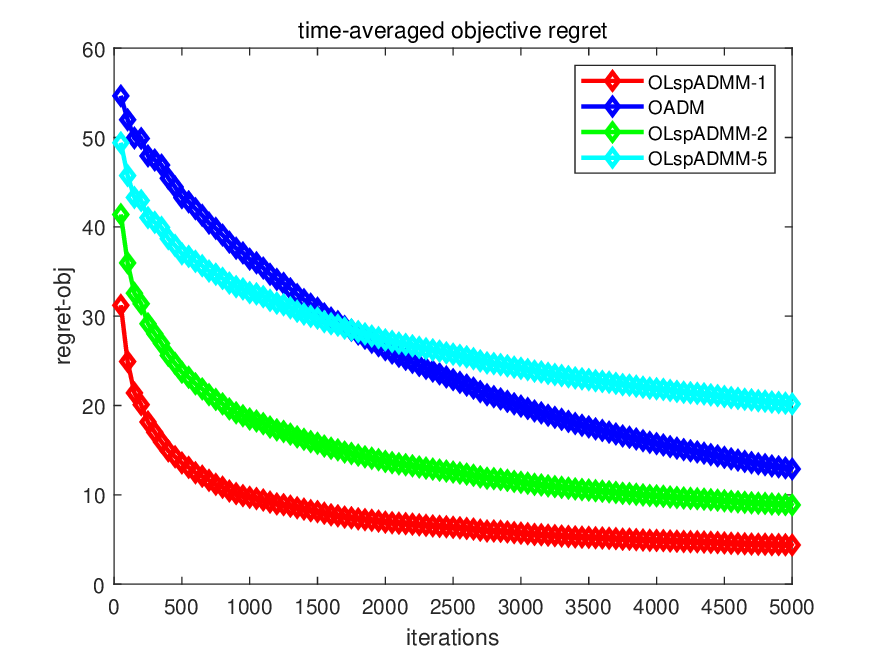} \label{figure1b}
	}
    \quad
    \subfigure[time-averaged constraint violation (n = 100)]{
	\includegraphics[height =4cm, width=5.5cm]{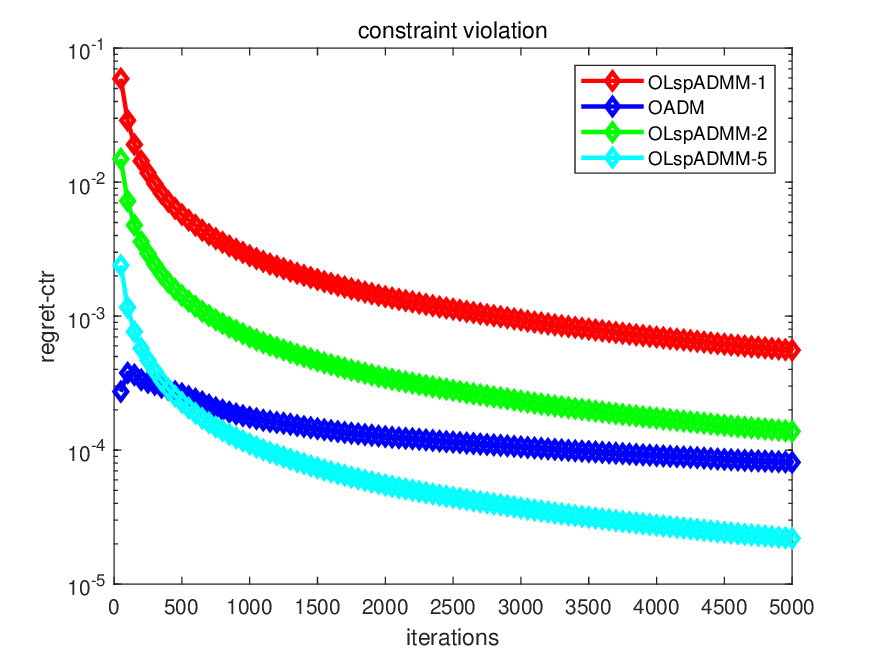} \label{figure1b}
	}
	\caption{Comparison of algorithms with respect to the time-average objective regret and constraint violation for TV.} \label{figure3}
\end{figure}

The results are shown in Figure \ref{figure3}. Since parameter $\tau$ changes have little effect on the performance of Online-spADMM, we consider the effect of $\sigma$ on this algorithm. Combine Figure \ref{figure3} and Table \ref{table5}, we can see that Lasso and OQO have the similar numerical performance about the time-averaged objective regret and constraint violation. It's easy to observe that the larger $a$, the larger time-averaged objective regret and the time-averaged objective regret is the opposite. From table \ref{table6}, we can observe that Online-spADMM has a shorter running time.

\begin{table}[h]
\caption{Comparison between the time-averaged objective regret of Online semi-proximal ADMM with the different selections of parameter $a$ (OLspADMM-$a$) and Online Alternating Direction Method(OADM), where the total number of iterations $T = 5000$.} \label{table5}
\begin{tabular*}{\textwidth}{@{\extracolsep\fill}lcccc}
\toprule%
\multicolumn{1}{@{}c@{}}{Dimension} &
\multicolumn{1}{@{}c@{}}{OLspADMM-1} & \multicolumn{1}{@{}c@{}}{OLspADMM-2} & \multicolumn{1}{@{}c@{}}{OLspADMM-5} &
\multicolumn{1}{@{}c@{}}{OADM}\\
\midrule
10  & 0.166 & 0.366 & 0.930  & 1.298  \\
20  & 0.330 & 0.730 & 1.860  & 2.598  \\
50  & 1.005 & 2.161 & 5.477  & 6.408  \\
100 & 4.382 & 8.863 & 20.180 & 12.890 \\
\bottomrule
\end{tabular*}
\end{table}

\begin{table}[h]
\caption{Same as Table \ref{table2} but for TV.} \label{table6}
\begin{tabular*}{\textwidth}{@{\extracolsep\fill}lcc}
\toprule%
&
 \multicolumn{1}{@{}c@{}}{Online-spADMM} & \multicolumn{1}{@{}c@{}}{OADM} \\\cmidrule{2-3}%
IterNum & n = 10$|$n = 20$|$n = 50$|$n = 100 & n = 10$|$n = 20$|$n = 50$|$n = 100 \\
\midrule
5000  & 0.028 $|$ 0.029 $|$ 0.046 $|$ 0.109 & 0.061 $|$ 0.071 $|$ 0.146 $|$ 0.382 \\
10000 & 0.053 $|$ 0.057 $|$ 0.091 $|$ 0.218 & 0.118 $|$ 0.141 $|$ 0.291 $|$ 0.763 \\
20000 & 0.103 $|$ 0.114 $|$ 0.180 $|$ 0.442 & 0.237 $|$ 0.282 $|$ 0.578 $|$ 1.521 \\
50000 & 0.257 $|$ 0.287 $|$ 0.451 $|$ 1.112 & 0.587 $|$ 0.702 $|$ 1.459 $|$ 3.822 \\
\bottomrule
\end{tabular*}
\end{table}

\section{Conclusions}
 In this paper, we have established regrets for objective and constraint violation of Online-spADMM for solving online linearly constrained convex composite optimization problems 
 {and presented numerical experiments to verify our theoretical results}. One significant feature of our approach is that the bounds for objective and constraint violation are obtained under weak  assumptions for objective functions.
  As the  bound for solution regret in Theorem \ref{iter-complexityanother} is not satisfactory, an  important issue left unanswered is to
  find  sufficient conditions for ensuring ${\rm O}(\sqrt{N})$ regret bound for  solution errors. For spADMM, whether we can obtain ${\rm O}(\sqrt{N})$ constraint violation regret bound  when $\sigma$ is a fixed constant is another topic worth studying.
This paper only discusses the case for constant constraint set $\Phi=\{(x,z): Ax+Bz=c\}$, a difficult problem is left to study is the case when $A,B,c$ is changing with time $t$, just like the online linear optimization considered by \cite{AWye2014}, or even a more complicated case where constraints are time-varying inequalities considered in \cite{YMNeely2017}. These  online optimization models are worth studying as they cover a large number of important practical problems.\\[20pt]

{\bf Appendix A. Proof of Theorem \ref{inequality}}. Since
$$
\begin{array}{l}
x^{k+1}\in \hbox{arg}\min \, \L^k_{\sigma }(x,z^k; y^k) +\frac{1}{2}\|x-x^k\|^2_{\sigma\S_k}\, ,\\[2mm]
z^{k+1}\in \hbox{arg}\min\, \L^k_{\sigma  }(x^{k+1},z; y^k) +\frac{1}{2}\|z-z^k\|^2_\T,
\end{array}
$$
we have, from the optimality for convex programming, that
$$
\begin{array}{l}
A^*(-y^k-\sigma(Ax^{k+1}+Bz^k-c))-\sigma\S_k(x^{k+1}-x^k)\in \partial f_k(x^{k+1}),\\[6pt]
B^*(-y^k-\sigma(Ax^{k+1}+Bz^{k+1}-c))-\T(z^{k+1}-z^k)\in \partial g(z^{k+1}),
\end{array}
$$
which, combing $y^{k+1}=y^k+\tau\sigma (Ax^{k+1}+Bz^{k+1}-c)$, implies
\begin{equation}\label{eq:a1}
\begin{array}{r}
A^*[-y^{k+1}-(1-\tau)\sigma(Ax^{k+1}+Bz^{k+1}-c))-\sigma B(z^k-z^{k+1})] \quad \quad \quad \quad  \quad \,\,\\[6pt]
-\sigma\S_k(x^{k+1}-x^k)\in \partial f_k(x^{k+1}),\\[6pt]
B^*[-y^{k+1}-(1-\tau)\sigma(Ax^{k+1}+Bz^{k+1}-c)]-\T(z^{k+1}-z^k)\in \partial g(z^{k+1}).
\end{array}
\end{equation}
It follows from the convexity of $f_k$ that
\begin{equation}\label{eq:a2}
\begin{array}{ll}
f_k(x^{k+1})-f_k(\widehat x) &\leq \langle -y^{k+1},A(x^{k+1}-\widehat x)\rangle+\sigma \langle B(z^{k+1}-z^k),A(x^{k+1}-\widehat x)\rangle \\[6pt]
& \quad \quad -(1-\tau)\sigma \langle Ax^{k+1}+Bz^{k+1}-c, A(x^{k+1}-\widehat x)\rangle\\[6pt]
& \quad \quad -\langle \sigma\S_k(x^{k+1}-x^k),x^{k+1}-\widehat x \rangle-\displaystyle \frac{1}{2}\|x^{k+1}-\widehat x\|_{\Sigma_{f_k}}^2\\[6pt]
& =\langle -y^{k+1},A(x^{k+1}-\widehat x)\rangle+\sigma \langle B(z^{k+1}-z^k),Ax^{k+1}+B\widehat z-c\rangle \\[6pt]
& \quad \quad -(1-\tau)\sigma \langle Ax^{k+1}+Bz^{k+1}-c, A(x^{k+1}-\widehat x)\rangle\\[6pt]
& \quad \quad -\langle \sigma\S_k(x^{k+1}-x^k),x^{k+1}-\widehat x \rangle-\displaystyle \frac{1}{2}\|x^{k+1}-\widehat x\|_{\Sigma_{f_k}}^2\\[6pt]
\end{array}
\end{equation}
From the convexity of $g$, we have
\begin{equation}\label{eq:a3}
\begin{array}{ll}
g(z^{k+1})-g(\widehat z) & \leq  \langle -y^{k+1}-(1-\tau)\sigma(Ax^{k+1}+Bz^{k+1}-c), B(z^{k+1}-\widehat z)
\rangle\\[6pt]
& \quad \quad - \langle \T(z^{k+1}-z^k),z^{k+1}-\widehat z \rangle -\displaystyle \frac{1}{2}\|z^{k+1}-\widehat z\|_{\Sigma_g}^2.
\end{array}
\end{equation}
Adding both sides of (\ref{eq:a2}) and (\ref{eq:a3}), we obtain
$$
\begin{array}{l}
[f_k(x^{k+1})+g(z^{k+1})]-[f_k(\widehat x)+g(\widehat z)] \\[6pt]
\quad \leq  \left\langle -y^{k+1},Ax^{k+1}+Bz^{k+1}-c \right\rangle-
(1-\tau)\sigma \|Ax^{k+1}+Bz^{k+1}-c\|^2\\[6pt]
\quad \quad +\sigma \left \langle B(z^{k+1}-z^k),Ax^{k+1}+B\widehat z-c\right \rangle-\left \langle \sigma\S_k(x^{k+1}-x^k),x^{k+1}-\widehat x \right \rangle \\[6pt]
 \quad \quad -\left \langle \T(z^{k+1}-z^k),z^{k+1}-\widehat z\right \rangle-\displaystyle \frac{1}{2}\|x^{k+1}-\widehat x\|^2_{\Sigma_{f_k}}-\displaystyle \frac{1}{2}\|z^{k+1}-\widehat z\|^2_{\Sigma_g}\\[10pt]
\quad = (\tau\sigma)^{-1}\left \langle y^{k+1},y^k-y^{k+1} \right\rangle-(1-\tau)\sigma \|Ax^{k+1}+Bz^{k+1}-c\|^2\\[6pt]
 \quad \quad +\displaystyle \frac{\sigma}{2}\left\{ \|Bz^k-B\widehat z\|^2-\|Bz^{k+1}-B\widehat z\|^2+\|Ax^{k+1}+Bz^{k+1}-c\|^2-
\|Ax^{k+1}+Bz^k-c\|^2\right\}\\[6pt]
\quad \quad -\left \langle \sigma\S_k(x^{k+1}-x^k),x^{k+1}-\widehat x\right \rangle-\left \langle \T(z^{k+1}-z^k),z^{k+1}-\widehat z\right \rangle
 -\displaystyle \frac{1}{2}\|x^{k+1}-\widehat x\|^2_{\Sigma_{f_k}}-\displaystyle \frac{1}{2}\|z^{k+1}-\widehat z\|^2_{\Sigma_g}.
 \end{array}
$$
From this and the identities
$$
\begin{array}{ll}
\left \langle y^{k+1},y^k-y^{k+1} \right \rangle &= \displaystyle \frac{1}{2}\left[ \|y^k\|^2-\|y^k-y^{k+1}\|^2-\|y^{k+1}\|^2\right]\\[6pt]
&= \displaystyle \frac{1}{2}\left[ \|y^k\|^2-\|y^{k+1}\|^2\right]-\displaystyle \frac{1}{2}[\sigma \tau]^2\|Ax^{k+1}+Bz^{k+1}-c\|^2,
\end{array}
$$
we obtain
\begin{equation}\label{eq:a4}
\begin{array}{l}
[f_k(x^{k+1})+g(z^{k+1})]-[f_k(\widehat x)+g(\widehat z)] \\[10pt]
\quad \leq  \displaystyle \frac{1}{2\sigma\tau}\left[ \|y^k\|^2-\|y^{k+1}\|^2\right]-\displaystyle \frac{1}{2}(2-\tau)\sigma
\|Ax^{k+1}+Bz^{k+1}-c \|^2\\[10pt]
\quad \quad +\displaystyle \frac{\sigma}{2}\left\{ \|Bz^k-B\widehat z\|^2-\|Bz^{k+1}-B\widehat z\|^2+\|Ax^{k+1}+Bz^{k+1}-c\|^2\right.\\[10pt]
\quad \quad \quad \quad \quad\displaystyle \left.-
\|Ax^{k+1}+Bz^k-c\|^2\right\}
 -\left \langle \sigma\S_k(x^{k+1}-x^k),x^{k+1}-\widehat x\right \rangle\\[10pt]
\quad \quad \quad \quad \quad-\left \langle \T(z^{k+1}-z^k),z^{k+1}-\widehat z\right \rangle
 -\displaystyle \frac{1}{2}\|x^{k+1}-\widehat x\|^2_{\Sigma_{f_k}}-\displaystyle \frac{1}{2}\|z^{k+1}-\widehat z\|^2_{\Sigma_g}.
 \end{array}
\end{equation}
Obviously we have the following equalities
\begin{equation}\label{eq:a5}
\begin{array}{l}
\|Ax^{k+1}+Bz^{k+1}-c\|^2-
\|Ax^{k+1}+Bz^k-c\|^2\\[6pt]
\quad =\|Ax^{k+1}+Bz^{k+1}-c\|^2-\|(Ax^{k+1}+Bz^{k+1}-c)+B(z^k-z^{k+1})\|^2\\[6pt]
\quad =-2 \left \langle Ax^{k+1}+Bz^{k+1}-c,B(z^k-z^{k+1}) \right \rangle-\|B(z^{k+1}-z^k)\|^2
\end{array}
\end{equation}
and for $R^k=Ax^k+Bz^k-c$, we have
\begin{equation}\label{eq:a6}
\begin{array}{l}
\sigma \left \langle B(z^{k+1}-z^k),Ax^{k+1}+Bz^{k+1}-c \right \rangle\\[6pt]
=(1-\tau)\sigma \left \langle B(z^{k+1}-z^k),Ax^{k+1}+Bz^{k+1}-c \right \rangle+\left \langle B(z^{k+1}-z^k),y^k-y^{k+1} \right \rangle\\[6pt]
=(1-\tau)\sigma \left \langle B(z^{k+1}-z^k),Ax^k+Bz^k-c \right \rangle+\left \langle B(z^{k+1}-z^k),y^k-y^{k+1} \right \rangle\\[6pt]
\quad \quad +(1-\tau)\sigma \left \langle B(z^{k+1}-z^k),[R^{k+1}-R^k]\right \rangle\\[6pt]
=(1-\tau)\sigma \left \langle B(z^{k+1}-z^k),Ax^k+Bz^k-c \right \rangle\\[6pt]
\quad \quad + \left \langle B(z^{k+1}-z^k),y^k-(1-\tau)\sigma R^k-(y^{k+1}-(1-\tau)\sigma R^{k+1})\right \rangle\\[6pt]
=(1-\tau)\sigma \left \langle B(z^{k+1}-z^k),Ax^k+Bz^k-c \right \rangle\\[6pt]
\quad \quad + \left \langle (z^{k+1}-z^k),B^*(y^k-(1-\tau)\sigma R^k)-\T(z^k-z^{k-1})\right.\\[6pt]
\quad \quad\quad \quad \quad\quad \quad\quad \quad\quad \quad\left. -B^*(y^{k+1}-(1-\tau)\sigma R^{k+1})-\T(z^{k+1}-z^k)\right \rangle\\[6pt]
\quad \quad\quad \quad + \left \langle z^{k+1}-z^k,\T(z^k-z^{k-1}) - \T(z^{k+1}-z^k)\right \rangle.
\end{array}
\end{equation}
Since
$$
B^*(y^k-(1-\tau)\sigma R^k)-\T(z^k-z^{k-1})\in \partial g(z^k), \,B^*(y^{k+1}-(1-\tau)\sigma R^{k+1})-\T(z^{k+1}-z^k)
\in\partial g(z^{k+1}),
$$
one has that
$$
\begin{array}{l}
\left \langle (z^{k+1}-z^k),B^*(y^k-(1-\tau)\sigma R^k)-\T(z^k-z^{k-1})\right.\\[6pt]
\quad \quad\quad \quad \quad\quad \quad\quad \quad\quad \quad\left. -B^*(y^{k+1}-(1-\tau)\sigma R^{k+1})-\T(z^{k+1}-z^k)\right \rangle \leq 0.\\[6pt]
\end{array}
$$
Thus we have from (\ref{eq:a6}) that
\begin{equation}\label{eq:a7}
\begin{array}{l}
\sigma \left \langle B(z^{k+1}-z^k),Ax^{k+1}+Bz^{k+1}-c \right \rangle\\[6pt]
\leq
(1-\tau)\sigma \left \langle B(z^{k+1}-z^k),Ax^k+Bz^k-c \right \rangle\\[6pt]
\quad \quad\quad \quad + \left \langle z^{k+1}-z^k,\T(z^k-z^{k-1}) - \T(z^{k+1}-z^k)\right \rangle\\[6pt]
\leq (1-\tau)\sigma \left \langle B(z^{k+1}-z^k),Ax^k+Bz^k-c \right \rangle\\[6pt]
\quad \quad\quad \quad -\displaystyle \frac{1}{2}\|z^{k+1}-z^k\|_{\T}^2+\displaystyle \frac{1}{2}\|z^{k}-z^{k-1}\|_{\T}^2.
\end{array}
\end{equation}
Thus we have from (\ref{eq:a5}) and  (\ref{eq:a7}) that
\begin{equation}\label{eq:a8}
\begin{array}{l}
\displaystyle \frac{\sigma}{2}\left[\|Ax^{k+1}+Bz^{k+1}-c\|^2-\|Ax^{k+1}+Bz^{k}-c \|^2  \right]\\[6pt]
\quad =\displaystyle \frac{\sigma}{2}\left[-2 \langle Ax^{k+1}+Bz^{k+1}-c,B(z^k-z^{k+1})\rangle-\|B(z^{k+1}-z^{k})\|^2  \right]\\[6pt]
\quad =\sigma \langle B(z^k-z^{k+1}),  Ax^{k+1}+Bz^{k+1}-c\rangle-\displaystyle \frac{\sigma}{2}|B(z^{k+1}-z^{k})\|^2 \\[6pt]
\quad \leq (1-\tau)\sigma \left \langle B(z^{k+1}-z^k),Ax^k+Bz^k-c \right \rangle-\displaystyle \frac{\sigma}{2}|B(z^{k+1}-z^{k})\|^2\\[6pt]
\quad \quad\quad \quad -\displaystyle \frac{1}{2}\|z^{k+1}-z^k\|_{\T}^2+\displaystyle \frac{1}{2}\|z^{k}-z^{k-1}\|_{\T}^2.
\end{array}
\end{equation}
Combining (\ref{eq:a4}) and (\ref{eq:a8}), we obtain
\begin{equation}\label{eq:a9}
\begin{array}{l}
[f_k(x^{k+1})+g(z^{k+1})]-[f_k(\widehat x)+g(\widehat z)] \\[10pt]
\quad \leq  \displaystyle \frac{1}{2\sigma\tau}\left[ \|y^k\|^2-\|y^{k+1}\|^2\right]-\displaystyle \frac{1}{2}(2-\tau)\sigma
\|Ax^{k+1}+Bz^{k+1}-c \|^2\\[10pt]
\quad \quad +\displaystyle \frac{\sigma}{2}\left\{ \|Bz^k-B\widehat z\|^2-\|Bz^{k+1}-B\widehat z\|^2\right\}\\[10pt]
\quad \quad +(1-\tau)\sigma \left \langle B(z^{k+1}-z^k),Ax^k+Bz^k-c \right \rangle-\displaystyle \frac{\sigma}{2}\|B(z^{k+1}-z^{k})\|^2\\[6pt]
\quad \quad\quad \quad -\displaystyle \frac{1}{2}\|z^{k+1}-z^k\|_{\T}^2+\displaystyle \frac{1}{2}\|z^{k}-z^{k-1}\|_{\T}^2\\[6pt]
\quad \quad +\displaystyle \frac{1}{2}
 \left [\|x^k-\widehat x \|_{\sigma\S_k}^2-\|x^{k+1}-x^k\|_{\sigma\S_k}^2-\|x^{k+1}-\widehat x\|_{\sigma\S_k}^2\right]-\displaystyle \frac{1}{2}\|x^{k+1}-\widehat x\|^2_{\Sigma_{f_k}}\\[10pt]
\quad \quad  +\displaystyle \frac{1}{2}
 \left [\|z^k-\widehat z \|_{\T}^2-\|z^{k+1}-z^k\|_{\T}^2-\|z^{k+1}-\widehat z\|_{\T}^2\right]
 -\displaystyle \frac{1}{2}\|z^{k+1}-\widehat z\|^2_{\Sigma_g}.
 \end{array}
\end{equation}
We consider two cases when  $\tau \in (0,1]$ and $(1,(1+\sqrt{5})/2)$, respectively.\\
Case (i): when $\tau \in (0,1]$.  We have
\begin{equation}\label{eq:a9}
\begin{array}{l}
(1-\tau)\sigma \left \langle B(z^{k+1}-z^k),Ax^k+Bz^k-c \right \rangle\\[6pt]
\quad \, \leq \displaystyle \frac{1}{2}(1-\tau)\sigma [\|B(z^{k+1}-z^k)\|^2+\|Ax^k+Bz^k-c\|^2].
\end{array}
\end{equation}
It follows from (\ref{eq:a9}), in this case, that
$$
\begin{array}{l}
[f_k(x^{k+1})+g(z^{k+1})]-[f_k(\widehat x)+g(\widehat z)] \\[10pt]
\quad \leq  \displaystyle \frac{1}{2\sigma\tau}\left[ \|y^k\|^2-\|y^{k+1}\|^2\right]\\[10pt]
\quad \quad +\displaystyle \frac{\sigma}{2}\left\{ \|Bz^k-B\widehat z\|^2-\|Bz^{k+1}-B\widehat z\|^2\right\}\\[10pt]
\quad \quad\quad \quad +\displaystyle \frac{1}{2}\left[\|z^{k}-z^{k-1}\|_{\T}^2-\|z^{k+1}-z^k\|_{\T}^2\right] \\[6pt]
\quad \quad +\displaystyle \frac{1}{2}
 \left [\|x^k-\widehat x \|_{\sigma\S_k}^2-\|x^{k+1}-x^k\|_{\sigma\S_k}^2-\|x^{k+1}-\widehat x\|_{\sigma\S_k}^2\right]-\displaystyle \frac{1}{2}\|x^{k+1}-\widehat x\|^2_{\Sigma_{f_k}}\\[10pt]
 \end{array}
 $$
 \begin{equation}\label{eq:a10}
\begin{array}{l}
\quad \quad  +\displaystyle \frac{1}{2}
 \left [\|z^k-\widehat z \|_{\T}^2-\|z^{k+1}-z^k\|_{\T}^2-\|z^{k+1}-\widehat z\|_{\T}^2\right]
 -\displaystyle \frac{1}{2}\|z^{k+1}-\widehat z\|^2_{\Sigma_g}.\\[10pt]
 \quad \quad +\displaystyle \frac{1}{2}(1-\tau)\sigma \left [\|Ax^{k}+Bz^{k}-c \|^2-\|Ax^{k+1}+Bz^{k+1}-c \|^2\right]\\[10pt]
 \quad \quad -\displaystyle \frac{1}{2}\sigma \|Ax^{k+1}+Bz^{k+1}-c \|^2-\displaystyle \frac{\tau \sigma}{2}\|B(z^{k+1}-z^k)\|^2.
 \end{array}
\end{equation}
Case (ii): when $(1,(1+\sqrt{5})/2)$.  We have
\begin{equation}\label{eq:a11}
\begin{array}{l}
(1-\tau)\sigma \left \langle B(z^{k+1}-z^k),Ax^k+Bz^k-c \right \rangle\\[6pt]
\quad \, \leq \displaystyle \frac{1}{2}(\tau-1)\sigma [\tau\|B(z^{k+1}-z^k)\|^2+\tau^{-1}\|Ax^k+Bz^k-c\|^2].
\end{array}
\end{equation}
It follows from (\ref{eq:a9}), in this case, that
\begin{equation}\label{eq:a12}
\begin{array}{l}
[f_k(x^{k+1})+g(z^{k+1})]-[f_k(\widehat x)+g(\widehat z)] \\[10pt]
\quad \leq  \displaystyle \frac{1}{2\sigma\tau}\left[ \|y^k\|^2-\|y^{k+1}\|^2\right]\\[10pt]
\quad \quad +\displaystyle \frac{\sigma}{2}\left\{ \|Bz^k-B\widehat z\|^2-\|Bz^{k+1}-B\widehat z\|^2\right\}\\[10pt]
\quad \quad\quad \quad +\displaystyle \frac{1}{2}\left[\|z^{k}-z^{k-1}\|_{\T}^2-\|z^{k+1}-z^k\|_{\T}^2\right] \\[6pt]
\quad \quad +\displaystyle \frac{1}{2}
 \left [\|x^k-\widehat x \|_{\sigma\S_k}^2-\|x^{k+1}-x^k\|_{\sigma\S_k}^2-\|x^{k+1}-\widehat x\|_{\sigma\S_k}^2\right]-\displaystyle \frac{1}{2}\|x^{k+1}-\widehat x\|^2_{\Sigma_{f_k}}\\[10pt]
\quad \quad  +\displaystyle \frac{1}{2}
 \left [\|z^k-\widehat z \|_{\T}^2-\|z^{k+1}-z^k\|_{\T}^2-\|z^{k+1}-\widehat z\|_{\T}^2\right]
 -\displaystyle \frac{1}{2}\|z^{k+1}-\widehat z\|^2_{\Sigma_g}.\\[10pt]
 \quad \quad +\displaystyle \frac{1}{2}(1-\tau^{-1})\sigma \left [\|Ax^{k}+Bz^{k}-c \|^2-\|Ax^{k+1}+Bz^{k+1}-c \|^2\right]\\[10pt]
 \quad \quad -\displaystyle \frac{\sigma}{2\tau}(1+\tau-\tau^2) \|Ax^{k+1}+Bz^{k+1}-c \|^2-\displaystyle \frac{ \sigma}{2}(1+\tau-\tau^2)\|B(z^{k+1}-z^k)\|^2.
 \end{array}
\end{equation}
In view of (\ref{eq:a10}) and (\ref{eq:a12}), we obtain   formula
(\ref{eq:mainineq}).\hfill $\Box$

\section*{Declarations}
{\bf Acknowledgements.} We are grateful to the reviewer for his/her constructive feedback.\\

{\bf Authors' contributions.} All authors contributed equally to this study.\\

{\bf Availability of supporting data.} The datasets used or analyzed in the study are available from the author upon reasonable request.\\

{\bf Funding.} This work was supported by National Key R\&D Program of China (2022YFA1004000),
 the National Natural Science Foundation of China (12201097,12071055,12371298), and  partially supported by Dalian High-level Talent Innovation Project (2020RD09).\\

{\bf Ethical Approval.} Not applicable.\\

{\bf Competing interests.} The authors declare no competing interests.\\


\begin{thebibliography}{99}
\bibitem{AWye2014}
{S. Agrawal, Z. Z. Wang and Y. Y. Ye,
A dynamic near-optimal algorithm for online linear programming, Operations
Research, 62(4)(2014), 876-890.}

\bibitem{boyd11}
S. Boyd, N.  Parikh, E. Chu,  B. Peleato and J. Eckstein,
  Distributed optimization and statistical learning via the alternating direction method of multipliers,
 {Found. Trends Mach. Learn.}, 3 (2011), 1-122.

\bibitem{chaudhary2021safe}
{S. Chaudhary, D. Kalathil,
  Safe Online Convex Optimization with Unknown Linear Safety Constraints,  Proceedings of the AAAI Conference on Artificial Intelligence, 36(6)(2022), 6175-6182.}

\bibitem{CSToh2015}
L. Chen, D. F.  Sun and K.-C. Toh,
  An effcient inexact symmetric Gauss-Seidel based majorized ADMM for high-dimensional convex composite conic programming,  Mathematical Programming,   161(1-2) (2017), 237-270

%





\bibitem{CLSToh2015}
Y.~Cui, X. D.~Li, D. F.~Sun  and K.-C.~Toh,
  On the convergence properties of a majorized ADMM for linearly constrained convex optimization problems with coupled objective functions, Journal of Optimization Theory and Applications, 169 (2016),1013-1041.


\bibitem{JohnSinger2009}
{J. Duchi and  Y. Singer, Efficient online and batch learning using forward backward splitting, Journal of Machine Learning Research, 10 (2009), 2899-2934.}

\bibitem{Fazel} M. Fazel, T. K. Pong, D. F. Sun and P. Tseng,
Hankel matrix rank minimization with applications to system identification and realization, SIAM Journal on Matrix Analysis and  Applications,  34 (2013),946-977.

\bibitem{Gabay:1976ff}
D. Gabay and B. Mercier,
 A dual algorithm for the solution of nonlinear variational problems via finite element approximation,
{ Computational Mathematics and Applications}, 2 (1976),17-40.


\bibitem{glo75}
R. Glowinski and A. Marroco,
  Sur l'approximation, par \'el\'ements finis d'ordre un, et la r\'esolution, par p\'enalisation-dualit\'e d'une classe de probl\`emes de Dirichlet non lin\'eaires,
  {Revue fran\c{c}aise d'atomatique, Informatique Recherche Op\'erationelle. Analyse Num\'erique},  9 (1975),41-76.

\bibitem{HSZhang2018} D. R. Han, D. F. Sun, and L. W. Zhang, Linear Rate Convergence of the
Alternating Direction Method of Multipliers for Convex Composite Programming, Mathematics of Operations Research,
 43(2)(2018), 622-637.

\bibitem{Hazan2015}
{E. Hazan, Introduction to Online Convex Optimization,Foundations and Trends
in Optimization,  2(3-4)(2015),157-325.}



\bibitem{hosseini2014online}
S. Hosseini, A. Chapman and M. Mesbahi, Introduction to Online Convex Optimization,Foundations and Trends
in Optimization,  2(3-4)(2015),157-325.




\bibitem{JJArch2015}
{R. Jenatton, J. Huang, C. Archambeau, Adaptive algorithms for online convex optimization with long-term constraints, International Conference on Machine Learning, 2016, 402-411.}





\bibitem{Kalai2005}
A. Kalai and S. Vempala, Efficient algorithms for online decision problems, Journal of Computer and System Sciences, 71(3) (2005), 291-307.

\bibitem{Kiv1997} J. Kivinen and M. Warmuth, Exponentiated gradient versus gradient descent for linear predictors, Information and Computation,  132(1)(1997),1-64.

\bibitem{LSToh2014}
X. D.~Li, D. F.~Sun and K. -C.~Toh,
A Schur complement based semi-proximal ADMM for convex quadratic conic programming and extensions,  Mathematical Programming,  155 (2016),333-373.

\bibitem{LMercier1979}
P. L. Lions and B. Mercier,   Splitting algorithms for the sum of two nonlinear operators, SIAM Journal on  Numerical
Analysis, 16 (1979), 964-979.
%
\bibitem{Little1998} N. Littlestone, Learning quickly when irrelevant attributes abound: A new linear-threshold algorithm, Machine Learning, 2(1988), 285-318.

\bibitem{liu2018online}
{B. B. Liu, J. D. Li, Y. Q. Liu, X. J. Liang, L. Jian and H. Huan, Online newton step algorithm with estimated gradient, arXiv preprint arXiv:1811.09955, 2018.}



\bibitem{liu2018zeroth}
{S. J. Liu, J. Chen, P. Y. Chen and A. O. Hero, Zeroth-order online alternating direction method of multipliers: Convergence analysis and applications, International Conference on Artificial Intelligence and Statistics, 2018, 288-297.}



\bibitem{MRTal2012}
{M. Mohri, A. Rostamizadeh, and A. Talwalkar,
Foundations of Machine Learning,The MIT Press Cambridge, Massachusetts London, England,2012.}

\bibitem{MRYang2012}
{M. Mahdavi, R. Jin and T. B. Yang, Trading regret for efficiency: online convex optimization with long term constraints,Journal of Machine Learning Research, 13 (2012), 2503-2528.}

%
\bibitem{rockafellar76A}
R. T.~Rockafellar,
 Augmented Lagrangians and applications of the proximal point algorithm in convex programming,
  {Mathematics of Operations Research}, 1 (1976),97-116.

\bibitem{rockafellar76B}
R. T.~Rockafellar,
  Monotone operators and the proximal point algorithm,
 {SIAM Journal on Control and Optimization}, 14 (1976),877-898.


 \bibitem{Rosenblatt1988}
  F. Rosenblatt, The perceptron, A probabilistic model for information storage andorganizationinthebrain,Psychological Review, 65(1958),386-407. (Reprinted in Neurocomputing, (MIT Press, 1988)).

 \bibitem{rudin1992nonlinear}
  L. I. Rudin, S. Osher and E. Fatemi, Nonlinear total variation based noise removal algorithms, Physica D: nonlinear phenomena, 60(1-4)(1992),259-268.


\bibitem{Shai2007a}
S. Shalev-Shwartz, Online learning: Theory, algorithms, and applications, Ph.D. Thesis, The Hebrew University, 2007.

\bibitem{Shai2011}
{S. Shalev-Shwartz, Online Learning and Online
Convex Optimization, Foundations and Trends in
Machine Learning,
 4(2)(2011), 107-194.}


\bibitem{SSSBen2014}
{S. Shalev-Shwartz and S. Ben-David, Understanding Machine Learning,From Theory to Algorithms,Cambridge University Press,2014.}


 \bibitem{Shai2007b}S. Shalev-Shwartz and Y. Singer, A primal-dual perspective of online learning algorithms, Machine Learning Journal,  69(2)2(2007), 115-142.


\bibitem{STYang2015}
D. F. Sun,  K. -C.  Toh and L. Q. Yang,
 A convergent 3-block semi-proximal alternating direction method of multipliers for conic programming with 4-type constraints,
 SIAM Journal on  Optimization, 25 (2015),882-915.

\bibitem{tibshirani1996regression}
{R. Tibshirani, Regression shrinkage and selection via the lasso, Journal of the Royal Statistical Society, Series B, 58(1)(1996), 267-288.}

\bibitem{WangBan2013}
{H. H. Wang and  A. Banerjee, Online alternating direction method (longer version), the 29th International Conference on Machine Learning, 2012.}


\bibitem{STYang2015}
D. F. Sun,  K. -C.  Toh and L. Q. Yang,
 A convergent 3-block semi-proximal alternating direction method of multipliers for conic programming with 4-type constraints,
 SIAM Journal on  Optimization, 25 (2015),882-915.



\bibitem{Xiao2010}
{L. Xiao, Dual averaging methods for regularized stochastic learning and online optimization,Journal of Machine Learning Research, 11 (2010), 2543-2596.}

\bibitem{YNeedly2016}
{H. Yu and M. J. Neely,  A low complexity algorithm with ${\rm O}(\sqrt{T})$ regret and finite constraint violations for online convex optimization with long term constraints,The Journal of Machine Learning Research, 8(5)(2020), 1-25.}

\bibitem{YMNeely2017}
{Hao Yu, Michael J. Neely, Xiaohan Wei, Online convex optimization with stochastic constraints, Adv Neural Inf Process Syst, 30(2017).}
\end{thebibliography}
\end{document}